
%


%
%
%
\font\fourteenrm=cmr12  scaled \magstep1
\font\fourteenbf=cmbx12 scaled \magstep1
\font\fourteentt=cmtt12 scaled \magstep1
\font\fourteensl=cmsl12 scaled \magstep1
\font\fourteensy=cmsy10 scaled \magstep2
\font\fourteeni=cmmi12  scaled \magstep1
\font\fourteenit=cmti12 scaled \magstep1
\font\fourteensc=cmcsc10 scaled \magstep2
\font\fourteenbit=cmssi12 scaled \magstep1
\font\fourteenbbb=msym10 scaled \magstep2
\font\twelverm=cmr12
\font\twelvebf=cmbx12
\font\twelvett=cmtt12
\font\twelvesl=cmsl12
\font\twelvesy=cmsy10 scaled \magstep1
\font\twelvei=cmmi12
\font\twelveit=cmti12
\font\twelvesc=cmcsc10 scaled \magstep1
\font\twelvebit=cmssi12
\font\twelvebbb=msym10 scaled \magstep1
\font\tenrm=cmr10
\font\tenbf=cmb10 
\font\tentt=cmtt10 
\font\tensl=cmsl10 
\font\tensy=cmsy10 
\font\teni=cmmi10 
\font\tenit=cmti10 
\font\tensc=cmcsc10 
\font\tenbit=cmssi10
\font\tenbbb=msym10
\font\ninei=cmmi9 
\font\ninerm=cmr9

\font\ninesy=cmsy9 
\font\ninebbb=msym10 at 9 pt
\font\eighti=cmmi8
\font\eightrm=cmr8

\font\eightsy=cmsy8
\font\eightbbb=msym10 at 8 pt
\font\seveni=cmmi7 
\font\sevenrm=cmr7
\font\sevensy=cmsy7
%
%
%

%

%
\newfam\bbbfam
%
\def\tenpoint{
\def\rm{\fam0\tenrm}
\textfont0=\tenrm \scriptfont0=\eightrm \scriptscriptfont0=\sevenrm
\textfont1=\teni \scriptfont1=\eighti \scriptscriptfont1=\seveni
\textfont2=\tensy \scriptfont2=\eightsy \scriptscriptfont2=\sevensy
\textfont3=\tenex \scriptfont3=\tenex \scriptscriptfont3=\tenex
\textfont\itfam=\tenit
\def\it{\fam\itfam\tenit}
\textfont\slfam=\tensl
\def\sl{\fam\slfam\tensl}
\textfont\bffam=\tenbf
\def\bf{\fam\bffam\tenbf}
\textfont\ttfam=\tentt
\def\tt{\fam\ttfam\tentt}
\def\sc{\tensc}
\def\bit{\tenbit}
\def\bbb{\fam\bbbfam\twelvebbb}
\textfont\bbbfam=\tenbbb
\scriptfont\bbbfam=\eightbbb
\scriptscriptfont\bbbfam=\sevenbf
\normalbaselineskip=12pt
\setbox\strutbox=\hbox{\vrule height10pt depth4pt width0pt}%
\normalbaselines\rm}
%
%
\def\twelvepoint{
\def\rm{\fam0\twelverm}
\textfont0=\twelverm \scriptfont0=\ninerm \scriptscriptfont0=\sevenrm
\textfont1=\twelvei \scriptfont1=\ninei \scriptscriptfont1=\seveni
\textfont2=\twelvesy \scriptfont2=\ninesy \scriptscriptfont2=\sevensy
\textfont3=\tenex \scriptfont3=\tenex \scriptscriptfont3=\tenex
\textfont\itfam=\twelveit
\def\it{\fam\itfam\twelveit}
\textfont\slfam=\twelvesl
\def\sl{\fam\slfam\twelvesl}
\textfont\bffam=\twelvebf
\def\bf{\fam\bffam\twelvebf}
\textfont\ttfam=\twelvett
\def\tt{\fam\ttfam\twelvett}
\def\sc{\twelvesc}
\def\bit{\twelvebit}
\def\bbb{\fam\bbbfam\twelvebbb}
\textfont\bbbfam=\twelvebbb
\scriptfont\bbbfam=\ninebbb 
\scriptscriptfont\bbbfam=\sevenbf
\normalbaselineskip=14pt
\setbox\strutbox=\hbox{\vrule height10pt depth4pt width0pt}%
\normalbaselines\rm}
%
%
\def\fourteenpoint{
\def\rm{\fam0\fourteenrm}
\textfont0=\fourteenrm \scriptfont0=\twelverm \scriptscriptfont0=\tenrm
\textfont1=\fourteeni \scriptfont1=\twelvei \scriptscriptfont1=\teni
\textfont2=\fourteensy \scriptfont2=\twelvesy \scriptscriptfont2=\tensy
\textfont3=\tenex \scriptfont3=\tenex \scriptscriptfont3=\tenex
\textfont\itfam=\fourteenit
\def\it{\fam\itfam\fourteenit}
\textfont\slfam=\fourteensl
\def\sl{\fam\slfam\fourteensl}
\textfont\bffam=\fourteenbf
\def\bf{\fam\bffam\fourteenbf}
\textfont\ttfam=\fourteentt
\def\tt{\fam\ttfam\fourteentt}
\def\sc{\fourteensc}
\def\bit{\fourteenbit}
\def\bbb{\fam\bbbfam\twelvebbb}
\textfont\bbbfam=\fourteenbbb
\scriptfont\bbbfam=\tenrm 
\scriptscriptfont\bbbfam=\eightrm
\normalbaselineskip=16pt
\setbox\strutbox=\hbox{\vrule height10pt depth4pt width0pt}%
\normalbaselines\rm}
%
%
%
\twelvepoint
\abovedisplayskip 14pt plus 3pt minus 10pt%
\belowdisplayskip 14pt plus 3pt minus 10pt%
\abovedisplayshortskip 0pt plus 3pt%
\belowdisplayshortskip 8pt plus 3pt minus 5pt%
\parskip 3pt plus 1.5pt
\hsize=6.5in
\vsize=8.9in
%

%

%

\input psfig
\input amssym.def
\input amssym
\def\widehat{\mathaccent"0362 }
\def\ss{\smallskip}
\def\ms{\medskip}
\def\bs{\bigskip}
\def\ni{\noindent}
\def\cl{\centerline}

\def\QED {$\quad\square$}
\def\QP{\smallskip\leftskip=.4in\rightskip=.4in\noindent}
\def\ref{\hangindent=1pc \hangafter=1 \noindent}

\def\ssm{\smallsetminus}
\def\dist{{\rm dist}}

\def\[{$\,}
\def\]{\,$}
\def\C{{\bf C}}
\def\R{{\bf R}}
\def\Z{{\bf Z}}
\def\N{{\bf N}}

\def\E{{\cal E}}

\def\>{\succ}

\def\={~=~ }

\font\bit=cmssi12 at 12truept
\font\title=cmbx12

\mathsurround = 1pt
\abovedisplayskip=6pt
\belowdisplayskip=6pt
\parskip=2pt

\def\IMSmarkvadjust{0 pt}
\def\IMSmarkhadjust{0 pt}
\def\IMSmarkhpadding{0 pt}
\def\SBIMSMark#1#2#3{
 \font\SBF=cmss10 at 10 true pt
 \font\SBI=cmssi10 at 10 true pt
 \setbox0=\hbox{\SBF \hbox to \IMSmarkhpadding{\relax}
                Stony Brook IMS Preprint \##1}
 \setbox2=\hbox to \wd0{\hfil \SBI #2}
 \setbox4=\hbox to \wd0{\hfil \SBI #3}
 \setbox6=\hbox to \wd0{\hss
             \vbox{\hsize=\wd0 \parskip=0pt \baselineskip=10 true pt
                   \copy0 \break%
                   \copy2 \break%
                   \copy4 \break}}
 \dimen0=\ht6   \advance\dimen0 by \vsize \advance\dimen0 by 8 true pt
                \advance\dimen0 by -\pagetotal
	        \advance\dimen0 by \IMSmarkvadjust
 \dimen2=\hsize \advance\dimen2 by .25 true in
	        \advance\dimen2 by \IMSmarkhadjust

%
%
  \openin2=publishd.tex
  \ifeof2\setbox0=\hbox to 0pt{}
  \else 
     \setbox0=\hbox to 3.1 true in{
                \vbox to \ht6{\hsize=3 true in \parskip=0pt  \noindent  
                {\SBI Published in modified form:}\hfil\break
                \input publishd.tex 
                \vfill}}
  \fi
  \closein2
  \ht0=0pt \dp0=0pt
 \ht6=0pt \dp6=0pt
 \setbox8=\vbox to \dimen0{\vfill \hbox to \dimen2{\copy0 \hss \copy6}}
 \ht8=0pt \dp8=0pt \wd8=0pt
 \copy8
 \message{*** Stony Brook IMS Preprint #1, #2. #3 ***}
}

\def\IMSmarkhpadding{5em}
\SBIMSMark{1997/10}{originally ``On bicritical rational maps'', Sept. 1997}
{Revised and retitled, April 1999}

\def\bic{{\bf bicrit}}
\def\per{{\rm Per}}
\def\D{{\bf D}}

\def\CL{{\cal C}}
\def\SL{{\cal S}}
\def\m{{\cal M}}
\def\M{\widehat{\cal M}} 
\def\E{\widehat{\cal M}}
\def \iso{{\,\buildrel \simeq\over\longrightarrow\,}}
\def\bic {{\bf Bicrit}}
\def\C{{\bbb C}}
\def\R{{\textstyle\bbb R}}
\def\Z{{\bbb Z}}
\def\D{{\bbb D}}

\headline={%
    \ifnum\pageno > 1 {%
        \hss\tenrm\ifodd\pageno 1. CONJUGACY INVARIANTS\hss\folio
        \else\folio\hfill MAPS WITH TWO CRITICAL POINTS\hfill\fi}%
    \else \hss\vbox{}\hss
    \fi}
\footline={}

\cl{\bf  On Rational Maps with Two Critical Points}\ss

\cl{J. Milnor}\ms

{\QP\bit{\bf Abstract.} This is a preliminary investigation of the geometry and
dynamics of rational maps with only two critical points.\ms}

\cl{\bf Contents:}

\hskip .8in 0. Introduction

\hskip .8in 1. Conjugacy invariants and the moduli space \[\m\]

\hskip .8in 2. Fixed points and the curves \[\per_1(\lambda)\]

\hskip .8in 3. Shift locus or connectedness locus

\hskip .8in 4. The parabolic shift locus \[\SL_{\rm par}\cong\C\ssm
\overline\D\]

\hskip .8in 5. Real maps

\hskip .8in 6. The extended moduli space \[\E= \m\cup L_\infty\]

\hskip .8in 7. The extended hyperbolic shift locus  \[\widehat\SL_{\rm hyp}
\subset\E\]

\hskip .8in 8. The curves \[\overline\per_p(\lambda)\subset \E\]

\hskip .8in Appendix A: No Herman rings

\hskip .8in Appendix B: Totally disconnected Julia sets

\hskip .8in Appendix C: Cross-ratio formulas

\hskip .8in Appendix D: Entire and meromorphic maps

\hskip .8in References
\ms

{\bf \S0. Introduction.} We consider rational maps \[f:\widehat\C\to\widehat\C\]
of degree \[n\ge 2\] which are {\bit bicritical\/}, in the sense that they have
only two critical points. Note that every rational map of degree two
is bicritical. (See [M2], [R2], [R3], [Si], [Sti] for discussion of this case.)
For \[n>2\], bicriticality is a very strong restriction. In fact bicritical
maps seem to behave much more like quadratic rational maps than like general
rational maps of degree \[n\].

It is shown that
the moduli space $\m=\m_n$, consisting of all holomorphic conjugacy classes of
bicritical maps of degree \[n\],
is biholomorphic to \[\C^2\]. Furthermore, the Julia set of a bicritical
map is either
connected, or totally disconnected and isomorphic to the one sided shift on
\[n\] symbols. In the latter case this Julia set can be either hyperbolic
or parabolic. Correspondingly the moduli space splits as the disjoint union
$$ \m\=\CL\,\cup\,\SL_{\rm hyp}\,\cup\,\SL_{\rm par} $$
of the {\bit connectedness locus\/}, the
{\bit hyperbolic shift locus\/}, and the {\bit parabolic shift locus\/}.
Here \[\SL_{\rm hyp}\] is a connected open subset of \[\m\]
with free cyclic fundamental group,
while \[\SL_{\rm par}\] is a codimension one subset, conformally isomorphic
to \[\C\ssm\overline\D\].\ss

{\bf Remark 0.1.} There is another interesting trichotomy obtained by
considering the multipliers \[\lambda_1\,,\,\ldots\,,\,\lambda_{n+1}\]
at the various fixed points (see \S2). If we assume that \[|\lambda_j|\ne 1\]
for all of these multipliers, then there are three possibilities as follows.
If two of the fixed points are attracting, then we are in the {\bit principal
hyperbolic component\/}, and the Julia set is a quasicircle. If there is
just one attracting fixed point, then we are in the {\bit polynomial-like\/}
case, and can reduce to the polynomial case by a quasiconformal surgery.
(If \[N\] is a compact neighborhood of the attracting point which contains
exactly one critical value, with \[f(N)\] compactly contained in \[ N\],
then \[f\] carries
\[\widehat\C\ssm f^{-1}(N)\] onto \[\widehat\C\ssm N\] by a map which is
polynomial-like in the sense of [DH2].)
Note that the hyperbolic shift locus is included here.
Finally, it may happen that all \[n+1\] fixed points are strictly repelling.
This {\bit essentially non-polynomial-like\/}
case is the most interesting, since this is where we must look for any
new or exotic behavior. (Compare 3.2.)\ss

In order to understand limiting behavior
as the rational map becomes degenerate, it is convenient to introduce a partial
compactification of moduli space by adding a line \[L_\infty\cong\C\] of
``points at infinity''. The resulting {\bit extended moduli space\/}
\[\E=\m\cup L_\infty\] fibers as a complex line bundle
$$	\C~\hookrightarrow~\E~\buildrel X\over\longrightarrow~\widehat \C $$
over the Riemann sphere, with Chern number equal to \[n-1\].
Here \[X:\m\to\C\] is a certain conjugacy class invariant which can be
described up to sign as a cross-ratio (see 1.7), and \[X(L_\infty)=\infty\].
The connectedness locus \[\CL\subset\m\] has compact closure within \[\E\].

For each \[\lambda\in\C\ssm \{0\}\]
the curve \[\per_1(\lambda)\] consisting of conjugacy
classes of maps with a fixed point of multiplier \[\lambda\] forms a
holomorphic section of the line bundle \[\E\to\widehat\C\]. Any two such
sections have exactly \[n-1\] intersections,
counted with multiplicity. On the other hand, for \[\lambda=0\] the locus
\[\per_1(0)\] is not a section, but rather coincides with a fiber
$$	L_0\=\{\,(f)\in \E\,~;~X(f)=0\,\}~.$$
This fiber can be identified with the set of conjugacy classes of
unicritical polynomial maps \[z\mapsto z^n+{\rm constant}\].\ss

The moduli space \[\m=\m_n\cong\C^2\] contains a {\bit real subspace\/} \[\m_\R\cong
\R^2\]. This consists not only of conjugacy classes of maps with real
coefficients but also, when the degree \[n\] is odd, of a more exotic region
consisting of conjugacy classes of maps \[f\] which commute with the antipodal
map \[z\mapsto -1/\overline z\] of the Riemann sphere. Such \[f\] give rise to
dynamical systems on the nonorientable surface which is
obtained by identifying \[z\] with \[-1/\overline z\].
Similarly the extended moduli space \[\E\supset \m\] contains a real subset
\[\E_\R\] which fibers as a real line bundle
$$	\R~\hookrightarrow~\E_\R~\buildrel X\over\longrightarrow~\R\cup
\{\infty\} $$
over the circle \[\R\cup\{\infty\}\].
Topologically, \[\E_\R\] is either a cylinder or M\"obius band according as
\[n\] is odd or even.\ss

{\bf Erratum.} The preliminary version ``On bicritical rational maps''
of this note contained a wrong computation of the degree of
\[\per_p(\lambda)\], and also an incorrect proof in Appendix A.\ss

{\bf Acknowledgements.}
I am indebted to A. Epstein for suggesting the problem, to D. Schleicher
and M. Shishikura
for valuable assistance, particularly in the proof of 4.2, and to A. Douady
for suggestions leading to Appendices A and D. Furthermore, I am indebted to the
Gabriella and Paul Rosenbaum Foundation and to the National Science
Foundation (Grant DMS-9505833)
for their support of mathematical activity at the Institute for Mathematical
Sciences at Stony Brook.
\bs\bs\eject


\cl{\bf \S1. Conjugacy Invariants and the Moduli Space \[\m\]}\ms

Let \[\bic_n\] be the space of all bicritical maps of degree \[n\ge 2\].
(It is not hard to check that \[\bic_n\] is a smooth 5-dimensional complex
manifold.) By definition, two rational maps \[f\] and \[g\] are
(holomorphically) {\bit conjugate\/} if there exists a M\"obius automorphism
\[\phi\] of the Riemann sphere so that \[g=\phi\circ f\circ\phi^{-1}\].
We are interested in the {\bit moduli space\/} \[\m=\m_n\] consisting of
all conjugacy classes \[(f)\] of degree \[n\] bicritical maps. The first two
sections will provide a rather formal algebraic description of this space.


First consider
the {\bit marked\/} moduli space \[\m'\], consisting of conjugacy classes of
\[f\] with numbered critical points \[ c_1\,,\, c_2\]. Here, by definition,
the {\bit conjugacy class\/} of \[(f\,,\, c_1\,,\, c_2)\]
consists of all \[\big(\phi\circ f\circ\phi^{-1}\,,\,\phi( c_1)
\,,\,\phi( c_2)\big)\], where \[\phi\] ranges over
M\"obius automorphisms of the Riemann sphere \[\widehat\C\].
In order to construct a complete set of invariants for such an
\[(f\,,\, c_1\,,\, c_2)\], we proceed as follows.

{\QP{\bf Lemma 1.1.} \it If we put the critical points
\[ c_1\] at infinity and \[ c_2\] at zero, then \[f\] must have the form
$$	f(z)\= {az^n+b\over cz^n+d}~,\qquad{\rm with~derivative}\quad
	f'(z)\={nz^{n-1}(ad-bc)\over (cz^n+d)^2}~. \eqno (1)$$\ss}

(There should be no confusion between the
coefficient \[c\] and the critical points \[c_j\].)
Here the determinant \[ad-bc\] must be non-zero.
Note that this transformation depends only on the ratios \[(a:b:c:d)\].\ss

{\bf Proof of 1.1.} Write \[f(z)\]
as the quotient \[p(z)/q(z)\] and look at the equation \[f(z)=v\] or
\[p(z)-vq(z)=0\]. Suppose, to fix our ideas, that the two critical values
\[v_1=f(0)\] and \[v_2=f(\infty)\] are finite. Then the polynomial
\[p(z)-v_2q(z)\] has no finite roots, and hence must be constant.
Similarly, \[p(z)-v_1q(z)\] has no non-zero roots, hence must have the form
\[kz^m\]. Solving the resulting linear equations for \[p(z)\] and \[q(z)\],
the conclusion follows easily. The case where \[v_1\] or \[v_2\] is infinite
can be handled by a similar argument.\QED\ss

{\bf Remark.} The Julia set of a bicritical map of degree \[n\] always has
an \[n$-fold rotational symmetry about its critical points. If we use the
normal form (1), then this symmetry is expressed by the equation \[f(\omega z)
=f(z)\] and hence \[J=\omega J\], whenever \[\omega^n=1\].

{\QP{\bf Theorem 1.2.} \it A complete set of conjugacy invariants for a map
\[f\] in the normal form $(1)$, with
marked critical points at zero and infinity, is given by the expressions
$$	X\={bc\over ad-bc}\,,\qquad Y_1\={a^{n+1}b^{n-1}\over(ad-bc)^n}\,,\qquad
 Y_2\={c^{n-1}d^{n+1}\over(ad-bc)^n}~. \eqno (2)$$
These are subject to the relation
$$	Y_1Y_2\=X^{n-1}(X+1)^{n+1}~, \eqno (3)$$
but to no other relations. Hence the moduli space \[\m'\] consisting of all
such marked conjugacy classes is homeomorphic to the affine algebraic variety
consisting of all \[~(X,Y_1,Y_2)\,\in\,\C^3~\] satisfying equation \[(3)\].\ms}

As an example, the polynomial function \[f(z) =z^n+b\] corresponds to
the matrix
$$	\left[\matrix{ a & b \cr c & d\cr}\right ]\=
	\left[\matrix{ 1 & b \cr 0 & 1\cr}\right ]\qquad{\rm with~ invariants}
\qquad	X\=Y_2\=0\,,\quad Y_1=b^{n-1}~.$$
Compare Figures 1a and 2 for pictures of the \[Y_1$-coordinate plane
and the \[b$-coordinate plane for degree \[n=4\], in this locus \[X=Y_2=0\]
of polynomial maps.\ss

{\bf Proof of 1.2.}
Multiplying the four coefficients in (1) by a common factor,
we can normalize so that \[ad-bc=1\]. The coefficients are then uniquely
determined up to a common change of sign. Note that the expressions (2)
are all invariant under this transformation, since numerator and denominator
are homogeneous of the same degree. With this normalization, we can write
(2) in the simpler form
$$	X+1\=ad\,,\qquad X\=bc\,,\qquad Y_1\=a^{n+1}b^{n-1}\,,\qquad Y_2\=c^{n-1}d^{n+1}~.
\eqno (2') $$
Such a normal form with critical points at zero and infinity is not unique,
since we are
still free to conjugate by a M\"obius automorphism which fixes
both zero and infinity. If we write such an automorphism as \[\phi(z)=z/ t^2\],
then we must replace \[f\] by \[\phi\circ f\circ\phi^{-1}(z)=
f(t^2z)/t^2\]. A brief computation shows
that the four coefficients, normalized so that the determinant remains \[+1\],
are then transformed by the rule
$$	(a\,,\,b\,,\,c\,,\,d)~~\mapsto~~( t^{n-1}a\,,~~ b/t^{n+1}\,,~~
 t^{n+1}c\,,~~ d/t^{n-1})~. \eqno (4)$$
It is clear that the three expressions $(2')$ are
invariant under this transformation (4), and also under a simultaneous change
of sign for \[a,b,c,d\], and that they satisfy the required relation (3).

Conversely, given \[(X,Y_1,Y_2)\] satisfying (3), we must show that there is
one and only one corresponding choice of \[\pm(a,b,c,d)\], up to the
transformation (4).\ss

{\bf Case 1.} Suppose that \[Y_2\ne 0\], hence \[c\ne 0\] and \[d\ne 0\].
Then using (4) we can
make a linear change of variables so that \[d=1\]. It follows from $(2')$ that
$$	a\=X+1\,,\qquad b=X/c\,,\qquad c^{n-1}\=Y_2~.$$
Thus we obtain a normal form which is uniquely determined by \[X\] and \[Y_2\],
up to a choice of \[(n-1)$-st root for \[Y_2\]. With this choice of \[a,b,c,d\],
note that the relation
\[Y_1=a^{n+1}b^{n-1}\] follows from (3). But applying (4) again
with \[t\] equal to any \[(2n-2)$-nd root of unity, since \[t^{n-1}=
\pm 1\], we see that
$$	(a\,,\,b\,,\,c\,,\,d)~\mapsto~\pm(a\,,\,b/t^2\,,\,t^2c\,,\,d)~,$$
where \[t^2\] can be an arbitrary \[(n-1)$-st root of unity. Therefore,
the conjugacy class does not depend on a particular choice of \[(n-1)$-st
root \[c\].\ss

{\bf Case 2.} If \[Y_1\ne 0\] the argument is similar.\ss

{\bf Case 3.} If \[Y_1=Y_2=0\], then \[X\] must be either \[0\] or \[-1\] by
equation (3). If \[X=Y_1=Y_2=0\], then making use of the hypothesis that
\[ad-bc=1\] we see that \[b=c=0\] and that
 \[f\] is conjugate to $$z~\mapsto~ z^n~.$$
On the other hand, if \[X+1=Y_1=Y_2=0\], then it follows similarly
that \[a=d=0\] and that
\[f\] is conjugate to the map $$z~\mapsto~ 1/z^n~.$$
(Remark: These two exceptional points in moduli space will
often require special attention.) This completes the proof.\QED\ss

Now consider the quotient space of \[\m'\] under the involution 
$$	(f\,,\, c_1\,,\, c_2)~\leftrightarrow~
	(f\,,\, c_2\,,\, c_1) $$
which interchanges the two critical points.

{\QP{\bf Corollary 1.3.} \it This quotient space \[{\m}\],
consisting of all holomorphic conjugacy
classes of degree \[n\] bicritical maps, is biholomorphic
to \[\C^2\], with coordinates \[X\] and \[Y=Y_1+Y_2\].\ms}

We will use the notation \[L_{X_0}\] 
for the complex line consisting of all \[(f)\in \m\] with \[X(f)=X_0\].\ss

{\bf Proof of 1.3.}
If \[f\] is given by (1), then a holomorphically conjugate map
with critical points interchanged is given by
$$	{1\over f(1/z)}\={dz^n+c\over bz^n+a}~. $$
Thus \[(a,b,c,d)\leftrightarrow(d,c,b,a)\] and
\[Y_1\leftrightarrow Y_2\], with \[X\] fixed.
We must form the quotient of \[\m'\]
under this involution. Since the product
\[Y_1Y_2\] can be expressed as a smooth function of \[X\], it is easy to
check that the two quantities \[X\] and \[Y=Y_1+Y_2\] form a complete and
independent set of invariants for the quotient variety \[\m\].\QED\ss

{\bf Remark.} Note that the algebraic variety (3) has a singular point
at \[X=-1\], \[Y_1=Y_2=0\], and (if \[n\ge 3\]) another singular point
at \[X=Y_1=Y_2=0\]. However, by passing to the quotient
variety in which we unmark the critical points, these two singular points
miraculously disappear.\ss

{\QP{\bf Corollary 1.4.} \it The {\bit symmetry locus\/} \[\Sigma\subset \m\],
consisting of all conjugacy classes  of \[f\] which commute with some
M\"obius automorphism, is the variety defined by the equation
$$	Y^2\=4\,X^{n-1}(X+1)^{n+1} ~. $$\ss}

{\bf Proof.} First suppose that there exists a non-trivial automorphism
which fixes the two critical points. Then it must also fix the critical
values, hence the set of critical points must coincide with the set of
critical values. There are only two possibilities: Either \[f\] fixes both
critical points hence \[(f)\] is the conjugacy class of \[z\mapsto z^n\]
with \[X=Y=0\], or else \[f\] interchanges the two critical points, hence
\[(f)\] is the class of \[z\mapsto 1/z^n\] with \[X+1=Y_1=Y_2=0\].
(In these two exceptional cases, the group of
automorphims fixing the critical points is cyclic of order \[n-1\] or \[n+1\]
respectively, and the full group of automorphisms is dihedral of order
\[2(n-1)\] or \[2(n+1)\] respectively.)

If we exclude these two cases, then a non-trivial automorphism \[\iota\]
commuting with \[f\] must be an involution which
interchanges the two critical points, and must be unique. It is easy
to see that such an involution exists if and only if \[Y_1=Y_2=Y/2\].
Since the two exceptional cases also satisfy this equation, the
conclusion then follows from (3).\QED\ss

{\bf Remark 1.5.} When \[n\] is odd, this
symmetry locus is a reducible variety, splitting as \[\Sigma_+\cup\Sigma_-\]
where \[\Sigma_\pm\] is defined by
$$	Y\=\pm 2\,X^{(n-1)/2}(X+1)^{(n+1)/2}~.  $$
In fact a class \[(f)\in{\m}\] can be symmetric in two essentially
different ways when \[n\] is odd. Let \[\iota\] be the (usually unique)
involution which commutes with \[f\]. Then the two fixed points of \[\iota\]
are also fixed by \[f\] when \[(f)\in\Sigma_+\], but are interchanged by
\[f\] when \[f\in\Sigma_-\]. This can be proved by using the normal form
(1), taking \[\iota(z)\] to be \[1/z\], so that
$$	f(z)\=\pm\,{az^n+b\over bz^n+a}~. \eqno (5_\pm) $$
(Note that the two exceptional conjugacy classes,
where the involution \[\iota\] is not uniquely determined, constitute the
intersection \[\Sigma_+\cap\Sigma_-\].) For \[(f)\in\Sigma_+\], computation
shows that the two invariant fixed points have multipliers \[\lambda_1\]
and \[\lambda_2\] with sum \[\lambda_1+\lambda_2=2n(2X+1)\]
and with product \[\lambda_1\lambda_2=n^2\]. On the other hand, for \[(f)\in
\Sigma_-\] the two fixed points of \[\iota\] constitute a period two orbit
for \[f\] with multiplier \[n^2\]. It follows that \[\Sigma_-\] is contained
as one irreducible component in the curve \[\per_2(n^2)\] of \S8.

Evidently, the two halves of the symmetry locus
represent quite different dynamic behavior. For example if \[(f)\in\Sigma_-\]
then there are either two non-repelling fixed points or none. In either case,
it follows that the Julia set is connected. On the other hand, if we use the
normal form $(5_+)$ then a brief computation shows that the multiplier
at the fixed point \[1=f(1)=\iota(1)\] equals \[n(a-b)/(a+b)\]. Whenever this
fixed point is attracting, it follows by symmetry that both critical points
must lie in its immediate basin. Using Theorem B.5 in the Appendix,
it follows that the Julia set is totally disconnected.

For \[n\] even, the
symmetry locus is irreducible, conformally isomorphic to a punctured plane.
In fact each \[(f)\in\Sigma\] has a unique fixed point which is invariant
under the involution. The multiplier \[\lambda\] at this fixed point can
take any non-zero value, and computation shows that
\[X=(\lambda/n+n/\lambda-2)/4\] is then uniquely determined. (The
correspondence \[\lambda\mapsto X\] is two-to-one, since a generic fiber
\[L_X\] intersects \[\Sigma\] in two different points.)\ss

{\bf Remark 1.6.} The following observation is due to Adam Epstein. Again let
\[\iota\] be the (usually unique) involution commuting with \[f\].
Then there is a natural involution of the symmetry locus given
by \[(f)\mapsto(f^\star)\], where \[f^\star=f\circ \iota=\iota\circ f\]
has the same Julia set as \[f\]. The map \[f^\star\]
has invariants \[ X^\star=-1-X\] and \[ Y^\star=(-1)^nYX/(1+X)\].
When \[n\] is odd, note that this involution maps each irreducible
component \[\Sigma_\pm\] to itself. Using the normal form \[(5_\pm)\],
this involution interchanges the coefficients \[a\] and \[b\].

An interesting involution of the entire moduli space \[\m\] is given by the
correspondence \[(f)\mapsto (J\circ f)\], where \[J=J_f\] is the unique
involution of the Riemann sphere which fixes the two critical values. If \[f\]
is given by (1), then
$$	J\circ f(z)\= {az^n-b\over cz^n-d} ~, $$
with \[X(J\circ f)=X(f)\] and \[Y_i(J\circ f)=-Y_i(f)\]. Thus this involution
maps the symmetry locus to itself, interchanging \[\Sigma_+\] and \[\Sigma_-\]
in the odd degree case.\ms

The correspondence
\[(f)\mapsto X(f)\] is rather natural, and can be defined in several different
ways. For example we will see in \S2 that \[X(f)\] is linearly related to the
sum of the multipliers at the various fixed points of \[f\]. As another
example, note the following cross-ratio formula. (Compare Appendix C.)

{\QP{\bf Lemma 1.7.} \it If \[f\] is a rational map with critical points
\[ c_1\,,\, c_2\] and with critical values \[v_j=f( c_j)\],
then the invariant \[X=X(f)\] is equal to the negative of the cross-ratio
$$	{( c_1-v_1)\,(c_2- v_2)\over( c_1- c_2)\,
	(v_1-v_2)}~.$$\ss}


{\bf Proof.} This cross-ratio is clearly well defined and invariant under
conjugation. (Note that the denominator never vanishes.) Putting the
critical points at \[c_1=\infty\] and \[c_2=0\], the left hand factors
cancel and the cross-ratio reduces to
$$	{0-v_2\over v_1-v_2}\={-b/d\over a/c-b/d}\={-bc\over ad-bc}~, $$
as required.\QED\ss

For further cross-ratio formulas, see Appendix C.\ms

{\QP{\bf Corollary 1.8.} \it Denote
the modulus of an annulus \[{\rm A}\subset\widehat\C\] by \[\rm mod(A)\],
and let \[{\rm mod}(f)\ge 0\] be the largest possible
modulus of an annulus in \[\widehat\C\] which separates the critical values
of \[f\] from the critical points of \[f\] (taking \[{\rm mod}(f)=0\]
when there is no such annulus). Given a sequence of conjugacy
classes \[(f_i)\in \m\], the invariants \[|X(f_i)|\] tend to infinity if and
only if the invariants \[{\rm mod}(f_i)\] tend to infinity.\ss}

{\bf Proof.} In fact we will show that
$$	{\log(r)\over 2\pi}~\le~{\rm mod}(f)~\le~
	{\rm mod}\Big(\widehat\C\ssm\big([-1,0]\cup[r\,,\,+\infty]\big)\Big)
        ~, $$
where \[r=|X|\], and where both the upper and the lower bound tend to infinity
as \[r\to\infty\].
After a M\"obius automorphism, we may assume that the critical
points are located at \[0\,,\,-1\] and the corresponding critical values at
\[X\,,\,\infty\].
Now the lower bound is obtained by using
the round annulus \[\{z~;~1<|z|<r\}\], while the upper bound follows
from [A, Ch. III].\QED\ss

The space \[\C^2\] is an extremely flabby object, with a very large
group of holomorphic automorphisms. In \S6 we will impose a much more rigid
structure on the moduli space \[\m\cong\C^2\] by partially
compactifying it. The invariant \[X\] will play a key role
in this partial compactification, since it will serve as the projection map
of a canonical fibration, with typical fiber \[L_{X_0}=
\{(f)~;~X(f)=X_0\}\].

\bs\bs

\headline={%
    \ifnum\pageno > 1 {%
        \hss\tenrm\ifodd\pageno 2. FIXED POINTS\hss\folio
        \else\folio\hfill MAPS WITH TWO CRITICAL POINTS\hfill\fi}%
    \else \hss\vbox{}\hss
    \fi}

\cl{\bf \S2. Fixed Points and the curves \[\per_1(\lambda)\]}\ms

Recall that the {\bit multiplier\/} of a rational map \[f\] at a
finite fixed point \[z=f(z)\] is defined to be the first derivative
\[\lambda=f'(z)\]. (In the case of a fixed point at infinity the multiplier
is equal to the limit of \[1/f'(z)\] as \[z\to\infty\].) We first prove
the following.

{\QP{\bf Lemma 2.1.} \it Let \[f\] be a bicritical map of degree \[n\]
with invariants \[X\] and \[Y\]. If \[f\] has a fixed point of multiplier
\[\lambda\], then the product \[\lambda^nY\] can be expressed as a polynomial
function of degree \[2n\] in the variables \[X\] and \[\lambda\].\ss}

{\bf Definition.}
Let \[\per_1(\lambda)\subset\m\] be the set of all conjugacy classes of
bicritical maps which admit a fixed point of multiplier \[\lambda\].

{\QP{\bf Corollary 2.2.} \it For \[\lambda\ne 0\], the curve
\[\per_1(\lambda)\] can be described as the graph of a polynomial function
\[Y={\rm polynomial}_\lambda(X)\]. In particular, for \[\lambda\ne 0\],
each fiber \[L_{X_0}=\{(f)~;~ X(F)=X_0\}\] contains one and only one
conjugacy class \[(f)\] of maps which have a
fixed point of multiplier \[\lambda\].\ss}

{\bf Proof of 2.1 and 2.2.} We will use the normal form (1). First suppose that
\[\lambda\ne 0\]. Then the fixed point of multiplier \[\lambda\] must be
distinct from the two critical points \[0\] and \[\infty\]. After a linear
change of coordinates, we may assume that this fixed point is \[z=1\].
Thus \[1=f(1)=(a+b)/(c+d)\]. Multiplying the coefficients by a common
factor, we may assume that \[a+b=c+d=2\].
If we define parameters \[\mu\] and \[\xi\] by the equations
$$\eqalign{	a\=1+\mu+\xi~,\qquad& b\=1-\mu-\xi~,\cr
		c\=1-\mu+\xi~,\qquad& d\=1+\mu-\xi~,}	\eqno (6) $$
then a straightforward computation shows that
$$	ad-bc\=4\mu~,	\qquad \lambda\=f'(1)\={n(ad-bc)\over(c+d)^2}\=n\mu~,
	\qquad	X\={(1-\mu)^2-\xi^2\over 4\mu}~, $$
and that
$$  Y\={(1+\mu+\xi)^{n+1}(1-\mu-\xi)^{n-1} +(1-\mu+\xi)^{n-1}(1+\mu-\xi)^{n+1}
	\over (4\mu)^n}~. $$
Thus \[\mu^nY\] is equal to a polynomial of degree \[2n\] in the variables
\[\mu\] and \[\xi\].
Note that this expression for \[Y\] is unchanged if we replace \[\xi\] by
\[-\xi\]. (The involution \[\xi\leftrightarrow -\xi\] corresponds to the
conjugation \[f(z)\leftrightarrow 1/f(1/z)\].) Hence \[\mu^nY\] can be
expressed as a polynomial function of \[\mu\] and \[\xi^2\]. Substituting
\[\xi^2=(1-\mu)^2-4\mu X\] and \[\mu=\lambda/n\], we obtain the required
polynomial expression for \[\lambda^nY\], of degree \[n\] in \[\lambda\]
and \[X\]. This proves 2.2.

To prove 2.1, we
must also check that this same polynomial relation remains valid when
\[f\] has a fixed point (necessarily \[0\] or \[\infty\]) of multiplier
\[\lambda=0\]. In that case, the product
\[\lambda^nY\] is certainly zero, and \[X=0\]
so that \[\xi=\pm 1\], and the numerator of the expression for \[Y\]
is identically zero, as required.\QED\ss

We can describe the form of these polynomial relations more precisely as
follows. We will continue to work with the quotient \[\mu=\lambda/n\].

{\QP{\bf Theorem 2.3.} \it For each \[n\] there exist polynomials
$$ P_0(X)\,,~P_1(X)\,,~P_2(X)\,,~\ldots\,,~P_{n+1}(X) $$
with integer coefficients, where each \[P_k(X)\]
has degree \[\le k\], so that
$$\eqalign{	\mu\, Y&\=\cr
	&P_{n+1}(X)-\mu P_n(X)+\mu^2P_{n-1}(X)-+\cdots+(-\mu)^nP_1(X)
	+(-\mu)^{n+1}P_0(X)~.} $$\ms}

(As explicit examples, we have
$${	\mu\,Y\= X\, -\,\mu\,X(2X-1)\, +\,\mu^2(4X+1)\,-\,\mu^3}\qquad 
{\rm for}\quad n=2~, $$
and
$$  \mu\,Y\=X^2-2\mu \,X^2(X-1)+\mu^2(9X^2+4X+1)-2\mu^3(3X+1)+\mu^4\qquad{\rm
for}\quad n=3~.) $$\ss

{\bf Proof of 2.3.} From the proof of 2.1, it follows easily that we
can define polynomials \[P_j(X)\] by the formula
$$	\mu Y=\sum_{i+j=n+1} (-\mu)^iP_j(X)~, \eqno (7) $$
where \[0\le j\le 2n\] or equivalently \[n+1\ge i\ge 1-n\], and where each
\[P_j(X)\] has degree \[\le j\].
What is new in 2.3 is the statement that the \[P_j\]
have integer coefficients, and that \[P_j=0\] for \[j>n+1\].
To prove this, we will derive the
same formula \[(7)\] in a different way. Again we use the normal form
(1) with a fixed point of multiplier \[\lambda\ne 0\] at \[z=1\], but
now we normalize so that \[ad-bc=1\], and set \[u=a+b=c+d\]. Then
computation shows that \[\mu=\lambda/n=1/u^2\]. Furthermore
$$	1\=ad-bc\=(u-b)(u-c)-bc\=u^2-(b+c)u~, $$
or in other words
$$	b+c\=u-1/u~.$$
Since \[bc=X\], it follows that \[b\] and \[c\] are the two roots of the
equation
$$	b^2-(u-1/u)b+X\=0~. $$
Thus
$$	2b\=u-1/u\,\pm\,\sqrt{(u-1/u)^2-4X} ~. $$
For \[|u|\] large, each of the two solutions \[b\] and \[c\] can be expressed
as a Laurent series in \[u\] with coefficients depending on \[X\].
One of these two solutions has the form
$$	b\=   {X\over u}+{X(1+X)\over u^3}+{X(1+X)(1+2X)\over u^5}+\cdots~,
  $$
tending to zero as \[u\to\infty\],
where the successive coefficients are polynomials in \[X\] with integer
coefficients which can be computed by a straightforward induction. The
other solution is then given by \[c=u-1/u-b=u-(1+X)/u-\cdots\], and is
asymptotic to \[u\]. From these we can compute
$$	Y\=(u-b)^{n+1}b^{n-1}\,+\,(u-c)^{n+1}c^{n-1} $$
as a Laurent series in \[u\]. This series begins as
$$\eqalign{	Y&\=u^2X^{n-1} \,-\,X^{n-1}(2X+1-n)\,+\,(n^2X^{n-1}+
\cdots+1)/u^2 \,+\,O(1/u^4)\cr
 &\= X^{n-1}/\mu\,-\,X^{n-1}(2X+1-n)\,+\,(n^2X^{n-1}+\cdots+1)\,\mu\, +\,O(\mu^2)
 ~,} $$
or in other words as
$$	\mu\,Y\= X^{n-1}\,-\,X^{n-1}(2X+1-n)\,\mu\,+\,(n^2X^{n-1}+\cdots+1)
\,\mu^2	\,+\,O(\mu^3) \eqno (8) $$
as \[\mu\to 0\].
Thus there are no terms in \[\mu^{i}\] with \[i<0\]. This proves that
formula (7) reduces to the required form.\QED\ss

Here is some more precise information about the polynomials \[P_k(X)\].
It will be convenient to introduce the abbreviation
$$	\beta(m,k)\=
\Big({m-k\atop k}\Big)\;+\;\Big({m-k-1\atop k-1}\Big) $$
for the sum of two binomial coefficients.

{\QP{\bf Lemma 2.4.} \it Each \[P_k(X)\] with \[0\le k\le n\]
is a polynomial of degree \[k\], however \[P_{n+1}(X)\] is a polynomial of
degree \[n-1\]. We have
$$\eqalign{	&P_0(X)\= 1~,\cr
	&P_1(X) \= 2nX+(n-1)~,\cr
	&\quad\cdots\cr}$$
$$\eqalign{
	&P_{n-1}(X)\=n^2X^{n-1}\,+\,\sum_{j=0}^{n-2}\Big({n+1\atop j}\Big)\,X^j
	~,\cr
	&P_n(X)~~~\= 2X^n\,+\,(1-n)\,X^{n-1}~,\cr
	&P_{n+1}(X)\=X^{n-1}~.} $$
The constant term \[P_k(0)\] in these polynomials is equal to the binomial
coefficient \[\Big({n-1\atop k}\Big)\], while the coefficient of the degree
\[k\] term is equal to the sum \[\beta(2n,k)\].\ss}

{\bf Proof Outline.} The explicit formulas for \[P_k(X)\] with \[k\ge n-1\]
can be derived from the computation (8). For the remaining information,
we again use the normal form (1) with \[ad-bc=1\] and with
\[u=a+b=c+d\]. Let \[s_k=b^k+c^k\], where \[s_1=b+c=u-1/u\].
Starting from the Newton formula
$$	s_{k+1}\=(b+c)s_k\,-\,(bc)\,s_{k-1}\=(u-u^{-1})s_k-X\,s_{k-1}~,$$
it follows inductively that we can express each \[s_k\] as a polynomial
function of \[u-u^{-1}\] and \[X\]. The precise formula is
$$	s_k\=\sum_{0\le j\le k/2}\; \beta(k,j)\,(-X)^j\,(u-u^{-1})^{k-2j}~. $$
(Note that this computation is independent of the degree \[n\].)
Equivalently, recalling that \[\mu=1/u^2\], we can write
$$	{s_k\over u^k}\=\sum_{0\le j\le k/2}\; \beta(k,j)\,(-\mu\,X)^j\,
 (1-\mu)^{k-2j}~, \eqno (9)$$
Now we can compute
$$\eqalign{	Y\=&a^{n+1}b^{n-1}+c^{n-1}d^{n+1}\cr
 \=& (u-b)^{n+1}b^{n-1}\,+\,c^{n-1}(u-c)^{n+1}\cr
	\=& u^{n+1}s_{n-1}\,-\,\Big({n+1\atop 1}\Big)\,u^ns_n\,+\,
  \Big({n+1\atop 2}\Big)u^{n-1}s_{n+1}\,-\,+\,\cdots\,
	\pm s_{2n}~,} $$
or equivalently
$$	\mu^n\,Y\={s_{n-1}\over u^{n-1}}\,-\,\Big({n+1\atop 1}\Big)
{s_n\over u^n}\,+\,\Big({n+1\atop 2}\Big){s_{n+1}\over u^{n+1}}\,
-\,+\,\cdots\,\pm{s_{2n}\over u^{2n}}~.\eqno (10)$$
Substituting \[(9)\] into \[(10)\], we obtain a fairly explicit
formula for \[\mu^n\,Y\]. Further details will be left to the reader.\QED\ss

As an application, we can give a more precise form of 2.2.

{\QP{\bf Corollary 2.5.} \it Each \[\per_1(\lambda)\subset \m\]
with \[\lambda\ne 0\] can be described a smooth curve of the form
$$	Y\= -2X^n\,+\,\big(n(\lambda+\lambda^{-1}+1)-1\big)X^{n-1}\,+\,
	\cdots+\lambda(n-\lambda)^{n-1}/n^n~. $$
If \[\lambda\ne\lambda'\] with \[\lambda\lambda'\ne 0\,,\,1\],
then it follows easily that the two curves \[\per_1(\lambda)\]
and \[\per_1(\lambda')\] have exactly \[n-1\] points of intersection, counted
with multiplicity. \ms}


{\bf Proof.} The first statement is proved by plugging the explicit
values from 2.4 into the equation 2.3. It follows that each intersection
\[\per_1(\lambda)\cap\per_1(\lambda')\] with \[\lambda\ne\lambda'\]
is described by a polynomial equation of the form
$$	n\,\Big(\lambda+{1\over\lambda}-\lambda'-{1\over\lambda'}\Big)\,
X^{n-1}~+~({\rm lower~terms})\= 0~. $$
In the generic case where \[\lambda\lambda'\ne 0,\,1\],
this equation has degree \[n-1\], and the assertion follows.\QED\ss

(On the other hand,
if \[\lambda\lambda'=1\]
then the leading coefficient of this polynomial equation is zero.
In this case, we must count one or more
``intersections at infinity'', in order to get the right number. (See \S6.)
A more significant exception occurs when \[\lambda'=0\]. In fact, to make
the count come out right, we should identify \[\per_1(0)\] with the
curve \[X^{n-1}=0\], or in other words with the locus \[L_0\] counted
\[n-1\] times. (Compare \S8.) In fact, as
\[\lambda'\to 0\] with  \[Y\] bounded,
the locus \[\per_1(\lambda')\] degenerates towards an \[(n-1$-sheeted
covering of the locus \[X=0\].
towards the curve \[X^{n-1}=0\] of multiplicity \[n-1\].)
\ms

{\bf Definition.}
Every rational map \[f\] of degree \[n\] has \[n+1\] fixed points
counted with multiplicity. Let \[\lambda_1\,,\,\ldots\,,\,\lambda_{n+1}\]
be the multipliers
at these fixed points, and let
$$	\sigma_k=\sum_{1\le i_1<\cdots<i_k\le n+1}\;\lambda_{i_1}\cdots
\lambda_{i_k} $$
be the \[k$-th {\bit elementary symmetric function\/}
of these multipliers. It is
convenient to set \[\sigma_0=1\]. Note that the quotient \[\sigma_k/n^k\]
can be described as the \[k$-th elementary symmetric function of the
quotients \[\mu_i=\lambda_i/n\].

{\QP{\bf Theorem 2.6.} \it These elementary symmetric functions can be
computed by the formula
$$\eqalign{	\sigma_k/n^k&\=P_k(X) \qquad{\rm for}\quad 0\le k\le n+1~,~~
 k\ne n~,\quad{\rm but}\cr
	\sigma_n/n^n&\=P_n(X)\,+\,Y~.\cr}\eqno (11) $$\ss}

{\bf Proof.} Each of the \[n+1\] multipliers \[\lambda=\lambda_k\]
trivially satisfies the polynomial equation
$$	\lambda^{n+1}-\sigma_1\lambda^n+\sigma_2\lambda^{n-1}
-+\cdots\pm\sigma_{n+1}\=\prod_{k=1}^{n+1}\big(\lambda-\lambda_k\big)\=
  0~. $$
Hence the quotient \[\mu=\lambda/n\] satisfies
$$	\mu^{n+1}-\sigma_1\mu^n/n+\sigma_2\mu^{n-1}/n^2
-+\cdots\pm\sigma_{n+1}/n^{n+1}\=  0~. $$
On the other hand, from 2.3 we see that
$$	\mu^{n+1}-P_1(X)\mu^n+P_2(X)\mu^{n-1}-+\cdots
	\mp (P_n(X)+Y)\mu\pm P_{n+1}(X)\=0~.$$
The difference of these two polynomial equations is a polynomial of degree
\[n\] in \[\mu\] which vanishes at all \[n+1\] of the \[\mu_i\]. If the
\[\mu_i\] are pairwise distinct (or in other words if the multipliers
\[\lambda_i=n\mu_i\] are pairwise
distinct), then it follows immediately that corresponding coefficients
are equal, which proves the identities (11). These identities follow
in the general case by continuity or by analytic continuation,
since for generic \[(f)\in \m\] the \[n+1\]
multipliers are indeed distinct. To prove this, we need only construct a
single example where the multipliers are distinct. For example if
\[f(z)=z^n+b\] then the multipliers are distinct provided that we exclude
\[n\] very special values of the parameter \[b\]. First we must guarantee
that \[(f)\ne\per_1(1)\], in order to be sure that the \[n+1\] fixed points
are distinct. But if \[(f)\in\per_1(1)\], then the fixed point equation
\[b=z-z^n\] together with the multiplier equation \[nz^{n-1}=1\] imply that
the invariant \[Y=b^{n-1}\] is equal to \[(n-1)^{n-1}/n^n\].
Finally, we must choose
\[b\ne 0\] to guarantee that two distinct fixed points, say \[z\ne\omega z\],
cannot have the same multiplier \[\lambda=nz^{n-1}=n(\omega z)^{n-1}\].
But this would imply that \[\omega^{n-1}=1\], and the fixed point equation
\[z-z^n=b\] would then yield \[(\omega z)-(\omega z)^n=\omega b\ne b\],
provided that \[b\ne 0\]. Also, no finite fixed point has multiplier zero
provided that \[b\ne 0\]. Thus generically the multipliers are distinct,
which completes the proof.\QED\ss

{\bf Remark 2.7.} As an immediate corollary of 2.4 and 2.6:
{\it We could equally well use the two
invariants \[\sigma_1\] and \[\sigma_n\] as coordinates for the moduli
space \[\m\cong\C^2\], in place of the invariants \[X\] and \[Y\] of \S1.}
(In practice, in \S6, it will be convenient to
use \[X\] and \[\sigma_n\] as coordinates.)\ss

{\bf Remark 2.8.} It seems surprising that every one of the
elementary symmetric functions \[\sigma_k\] with \[k\ne n\]
can be expressed as a function of \[X\] alone. Only \[\sigma_n\]
depends also on the coordinate \[Y\]. As an example to illustrate this
statement, consider the family of unicritical polynomials
$$	f(z)\=z^n+b~,$$
with invariants \[X=0\] and \[Y=b^{n-1}\]. For the special case \[b=0\],
there are two fixed points of multiplier zero and \[n-1\] fixed points
of multiplier \[n\], hence
$$	\sigma_k/n^k\=\Big({n-1\atop k}\Big) $$
for every \[k\].
It follows from 2.6 that this same formula holds for any value of the
parameter \[b\], provided that \[k\ne n\]. On the other hand for \[k=n\],
since this binomial coefficient is zero, it follows that
$$	\sigma_n/n^n\= Y\=b^{n-1}~. $$
For example
in the quadratic case \[f(z)=z^2+b\], it follows that the multipliers at the
finite fixed points satisfy \[\lambda_1+\lambda_2=2\] and
\[\lambda_1\lambda_2=4b\].\ss

{\bf Remark 2.9.} 
The {\bit holomorphic fixed point formula\/} asserts that
$$	\sum_1^{n+1}\,{1\over 1-\lambda_j}\= 1 $$ 
provided that \[\lambda_j\ne 1\] for all \[j\].
(See for example [M1].) This
gives rise to a linear relation between the \[\sigma_k\], or equivalently
between the \[P_k\], which takes the form
$$	\sum_0^{n+1}\;(-1)^k(n-k)\,\sigma_k\=\sum_0^{n+1}\;(-n)^k(n-k)\,P_k(X)
\=0~. $$ 
It follows by continuity that this relation still holds also when some of the
\[\lambda_j\] equal \[1\].
Note that the invariant \[Y\] is not involved, since
the coefficient of \[\sigma_n\] in this formula is zero.\ss

{\bf Remark 2.10.} It is sometimes convenient
to consider the moduli space for bicritical maps
with one marked fixed point. In this case, a complete set of invariants is
provided by \[X\] and \[Y\] together with the multiplier \[\lambda=n\,\mu\]
at this marked point. These are subject only to the relation \[\mu Y= X^{n-1}
-\mu P_n(X)+\mu^2 P_{n-1}(X) - + \cdots + (-\mu)^{n+1}\]. We can understand
the topology of the resulting variety better by introducing a new coordinate
\[Y_\mu=Y+P_n(X)-\mu P_{n-1}(X)+-\cdots+(-\mu)^n\] in place of \[Y\].
Then \[X\,,\,Y_\mu\]
and \[\mu\] are subject only to the relation \[\mu\,Y_\mu=X^{n-1}\].
For \[n\ge 3\] this variety has a singular point at \[X=Y=Y_\mu=\mu=0\].

\bs\bs

\headline={%
    \ifnum\pageno > 1 {%
        \hss\tenrm\ifodd\pageno 3. SHIFT LOCUS OR CONNECTEDNESS LOCUS\hss\folio
        \else\folio\hfill MAPS WITH TWO CRITICAL POINTS\hfill\fi}%
    \else \hss\vbox{}\hss
    \fi}

\cl{\bf\S3. Shift Locus or Connectedness Locus}\ms 

By definition, a conjugacy class \[(f)\] of degree \[n\] rational maps
belongs to the {\bit connectedness locus\/} \[\CL\]
if the Julia set \[J_f\] is connected; and belongs to the
{\bit shift locus\/} \[\cal S\] if \[J_f\] is totally disconnected
with \[f|_{J_f}\] topologically
conjugate to the one sided shift on \[n\] symbols. This section will prove
that every conjugacy class of maps with only two critical points
must belong to one or the other:

{\QP{\bf Theorem 3.1.} \it Every \[(f)\in \m\] belongs either to the
connectedness locus or to the shift locus.\ss}

{\bf Note:} 
If \[(f)\] belongs to the shift locus, then evidently both critical points
belong to the Fatou set \[\widehat\C\ssm J_f\], which is connected but far from
simply connected. There are two possibilities.
If \[f\] has an attracting fixed point, and hence is
hyperbolic on its Julia set, then we will say
that \[(f)\] belongs to the {\bit hyperbolic shift locus\/}
\[{\cal S}_{\rm hyp}\].  Otherwise, \[f\] must have a parabolic fixed
point, and we will say that \[(f)\] belongs to the {\bit parabolic shift
locus\/} \[{\cal S}_{\rm par}\]. Thus the
moduli space partitions as a disjoint union
$$	\m\=\CL\cup\SL_{\rm hyp}\cup\SL_{\rm par}~.$$
We will explore this partition of \[\m\] further in Sections 4 and 7.\ss

{\bf Remark 3.2.} In contrast with the polynomial case, we will see that the
connectedness locus is neither closed nor bounded in \[\m\] (although it has
compact closure in the extended moduli space \[\E\]). In
analogy with the polynomial case, one might be tempted to
conjecture that the
interior of the connectedness locus consists only of hyperbolic maps.
In fact this conjecture is true if we restrict attention to the open subset
consisting of \[(f)\] with at least one attracting fixed point. (Compare 0.1.)
However, the connectedness locus
also contains an ``essentially non-polynomial-like'' region \[\CL_{\rm NP}\]
consisting of maps for which
all \[n+1\] fixed points are strictly repelling. This region is certainly
contained in the interior of the connectedness locus, and yet contains many
non-hyperbolic maps. (Compare [R1].) Here is one example. If
$$	f(z)\=\kappa+(1-\kappa)/z^n\qquad{\rm with}\qquad \kappa^n=1\,,
~~\kappa\ne 1~, $$
then the critical points are \[0\] and \[\infty\] with
\[0\mapsto\infty\mapsto\kappa\mapsto 1\], so that
both critical orbits eventually land at the repelling fixed point \[1\].
It follows that the Julia set is the entire Riemann sphere, and that all
periodic orbits are strictly repelling. (This map lies in the locus \[X=-1\],
where one critical point maps directly to the other. Compare [BB].)
{\it It is conjectured that the region \[\CL_{\rm NP}\]
is a topological 4-cell.} When \[n=2\], this can be proved as follows.
Let \[I_j=1/(1-\lambda_j)\] be the holomorphic fixed point index at the
\[j$-th fixed point. (Compare 2.9.)
Then this region in moduli space can be identified with the star
shaped region consisting of unordered triples of complex numbers \[I_j\] with
\[0<{\rm Re}(I_j)<1/2\] and \[I_1+I_2+I_3=+1\]. On the other hand,
for \[n>2\], I don't know even whether \[\CL_{\rm NP}\] is simply
connected. Evidently an understanding of
the topology and dynamics associated with this region \[\CL_{\rm NP}\] would
be fundamental in reaching an understanding of bicritical maps.
\ss

The first ingredient in the proof of 3.1 is the following result which is
due to Shishikura.

{\QP{\bf Theorem A.1.} \it A rational map with two critical points cannot
have any Herman rings.\ss}

\ni In fact a proof of this statement
can be extracted from his paper [Sh] although it is not explicitly
stated there.  See Appendix A for a proof without using quasiconformal
surgery.

The second key ingredient is the following special case of a theorem of
Przytycki and Makienko.

{\QP{\bf Theorem B.5.} \it If a map \[f\] with two critical points has the
property that both critical values lie in a common Fatou component,
then \[(f)\] belongs to the shift locus.\ss}

In fact, more generally, it is shown in [Pr] and in [Ma] that
any rational map with all
critical values in a single Fatou component, when restricted to its Julia set,
is isomorphic to the one sided shift. However, since their argument is rather
complicated,  and since we need only the bicritical case, a proof
of B.5 is included in Appendix B.\ss

{\bf Remark 3.3.} Here is an alternative statement: {\it Suppose that
both critical orbits are eventually absorbed by an invariant Fatou component,
\[\Omega=f(\Omega)\]. Then \[(f)\] belongs to the shift locus.\/} In fact
such a Fatou component
must contain at least one of the two critical points, and hence must be fully
invariant, \[\Omega=f^{-1}(\Omega)\]. Hence it contains both critical
values, and 3.1 applies. (By way of contrast, a map
with three critical points may well have
connected Julia set, even though all
critical orbits are eventually absorbed by an invariant Fatou component.
For example the map \[f(z)=2+2z^3/(27(2-z))\] has critical points
\[0\,,\,3\,,\,\infty\] with orbit \[3\mapsto 0\mapsto 2\mapsto\infty\] ending
on a superattractive fixed point. The immediate
basin of infinity contains no other critical point, since no critical
orbit converges non-trivially to infinity, hence this basin is simply
connected. It follows, as in the proof of 3.1, that every
Fatou component is simply connected.)\ss

The proof of 3.1 will also make use of the following elementary
observation.

{\QP{\bf Lemma 3.4.} \it Let \[P\subset \widehat\C\] be a region bounded by
a simple closed curve which passes through neither
critical value. Then the pre-image of \[P\] under
\[f\] can be described as follows.\ss

{\bf Case 0.} If \[P\] contains no critical value, then \[f^{-1}(P)\] consists
of \[n\] disjoint simply-connected regions bounded by \[n\] disjoint simple
closed curves.

{\bf Case 1.} If \[P\] contains just one critical value, then \[f^{-1}(P)\]
is a single simply connected region bounded by a simple closed curve, and
maps onto \[P\] by a ramified \[n$-fold covering.

{\bf Case 2.} If \[P\] contains both critical values, then \[f^{-1}(P)\]
is a multiply connected region with \[n\] boundary curves, and
maps onto \[P\] by a ramified \[n$-fold covering.\ss}

The proof is straightforward. (Note that Cases 0 and 2
correspond to the ``inside'' and ``outside'' of the same simple closed
curve. In the case of just one critical point in \[P\], the set \[f^{-1}(P)\]
must be connected since \[f\] is locally \[n$-to-one near a critical point,
and the branched covering \[f^{-1}(P)\to P\]
is unique up to isomorphism, since the fundamental group of \[~P\ssm({\rm
critical~ value})~\] is free cyclic, so that there is only one \[n$-fold
covering of this set up to isomorphism.)\QED

{\bf Proof of 3.1.} Suppose that \[(f)\] is not in the shift locus, and
hence that no Fatou component contains more than one critical value.
If \[L\] is a loop in an arbitrary Fatou component, then using Sullivan's
Non-Wandering Theorem we see that some forward
image \[f^{\circ k}(L)\] lies in a simply connected region \[U\]
which is either:

(1) a linearizing neighborhood of some geometrically attracting periodic
point,

(2) a B\"ottcher neighborhood of some superattracting periodic point,

(3) an attracting petal for a parabolic point, or

(4) a Siegel disk.

\ni Here we are using the fact that there are no Herman rings
(Theorem A.1). Using 3.4, it follows by induction
on \[k\] that each component of \[f^{-k}(U)\] is simply connected.
This proves that every Fatou component is simply connected, and hence that
the Julia set is connected.\QED
\bs

\headline={%
    \ifnum\pageno > 1 {%
        \hss\tenrm\ifodd\pageno 4. PARABOLIC SHIFT LOCUS\hss\folio
        \else\folio\hfill MAPS WITH TWO CRITICAL POINTS\hfill\fi}%
    \else \hss\vbox{}\hss
    \fi}

\cl{\bf\S4. The Parabolic Shift Locus \[\SL_{\rm par}\cong\C\ssm\overline\D\]}\ms 

Recall from \S3 that the moduli space \[\m\] splits as a disjoint union
$$	\m\=\CL\,\cup\,\SL_{\rm hyp}\,\cup\,\SL_{\rm par}~.$$
Evidently \[{\cal S}_{\rm hyp}\] is an
open subset of moduli space, disjoint from the curve \[\per_1(1)\cong\C\],
while \[{\cal S}_{\rm par}\] is a relatively open subset of this
curve \[\per_1(1)\].\ss


We will first prove the following.

{\QP{\bf Lemma 4.1.} \it The parabolic shift locus is contained in the common
topological boundary \[\partial {\cal S}_{\rm hyp}=\partial\CL\]. Hence the
closure \[\overline\CL\subset \m\] is equal to \[\CL\cup\per_1(1)\],
with complement \[{\cal S}_{\rm hyp}\].\ss}

{\bf Proof.} For every \[(f)\] in the parabolic shift locus, we must show that
\[f\] can be approximated arbitrarily closely by a map with connected Julia
set, and also by a hyperbolic map with totally disconnected Julia set.
Let \[\{f_t\}\] be a holomorphic one-parameter family of maps with \[f_0=f\].
We will assume that each \[f_t\] has two critical points, and that
this family is not contained in \[\per_1(1)\]. Then for \[|t|\] small but
non-zero, the parabolic fixed point for \[f_0\] splits into two nearby fixed
points, with multipliers say \[\lambda_1\] and \[\lambda_2\]. As \[t\]
traverses a loop around \[t=0\], these two fixed points may be interchanged.
However, if we set \[t=u^2\], then both \[\lambda_1\] and \[\lambda_2\]
can certainly be expressed as single valued holomorphic functions of \[u\],
with \[\lambda_1(0)=\lambda_2(0)=1\]. Since these functions are non-constant,
we can choose \[u\] close to zero so that \[\lambda_1(u)\] takes any required
value close to \[1\].

First let us choose \[u\] so that \[\lambda_1(u)=e^{2\pi i/q}\], with \[q>1\].
Then the corresponding fixed point is parabolic, with at least two attracting
petals. Hence the associated Fatou set is not connected, and \[(f_{u^2})\]
must belong to the connectedness locus \[\CL\].

Now let us choose \[u\] so that \[\lambda_1(u)\] is real, with \[\lambda_1<1\],
so that the corresponding fixed point is strictly attracting.
For \[u\] sufficiently close to zero, we will show that \[(f_{u^2})\]
belongs to the hyperbolic shift locus. Choose a simple arc \[A\]
joining the two critical points with the Fatou set for \[f_0\]. Then for large
\[k\] the image \[f_0^{\circ k}(A)\] lies close to the parabolic point,
and within a sector of small angular size about the attracting direction
for this parabolic point. An easy perturbation argument then shows that
the same description holds for \[f_{(u^2)}\], provided that \[\lambda_1(u)<1\]
with \[u\] close to zero. (See for example [M3, \S4].)
Thus both critical values lie in a common Fatou
component, and it follows that \[(f_{u^2})\in {\cal S}_{\rm hyp}\].\QED


\midinsert
\cl{\psfig{figure=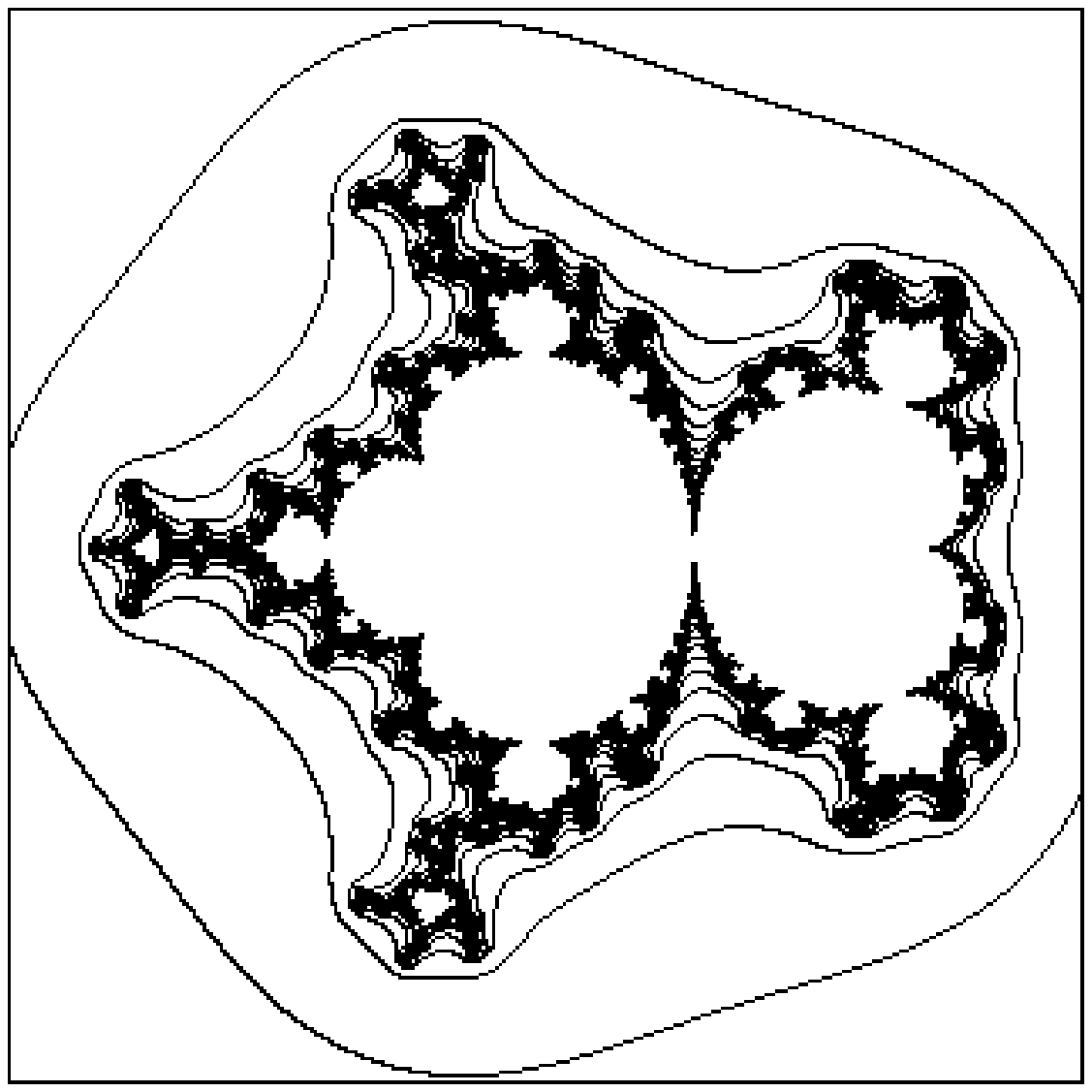,height=1.6in}
\psfig{figure=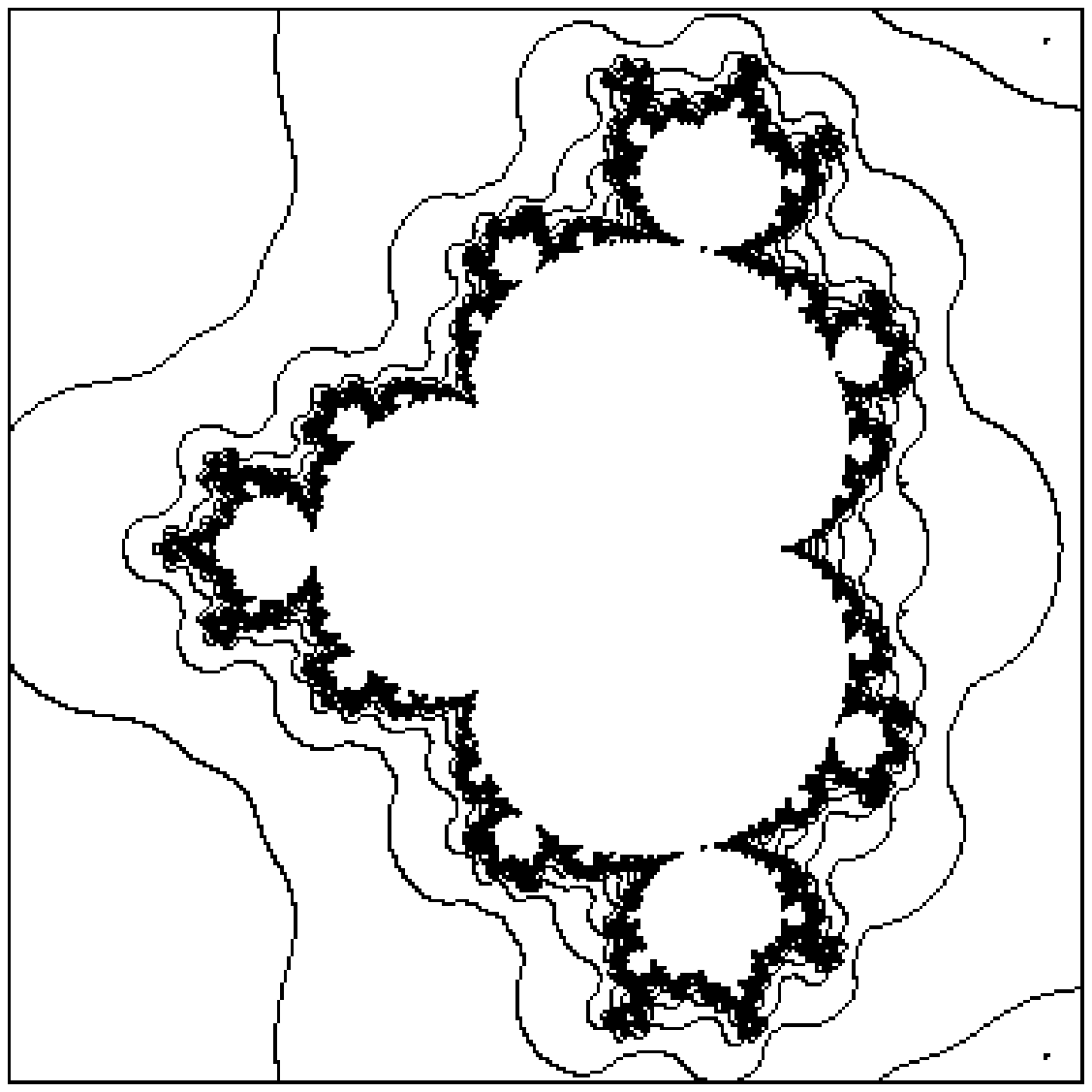,height=1.6in}
\psfig{figure=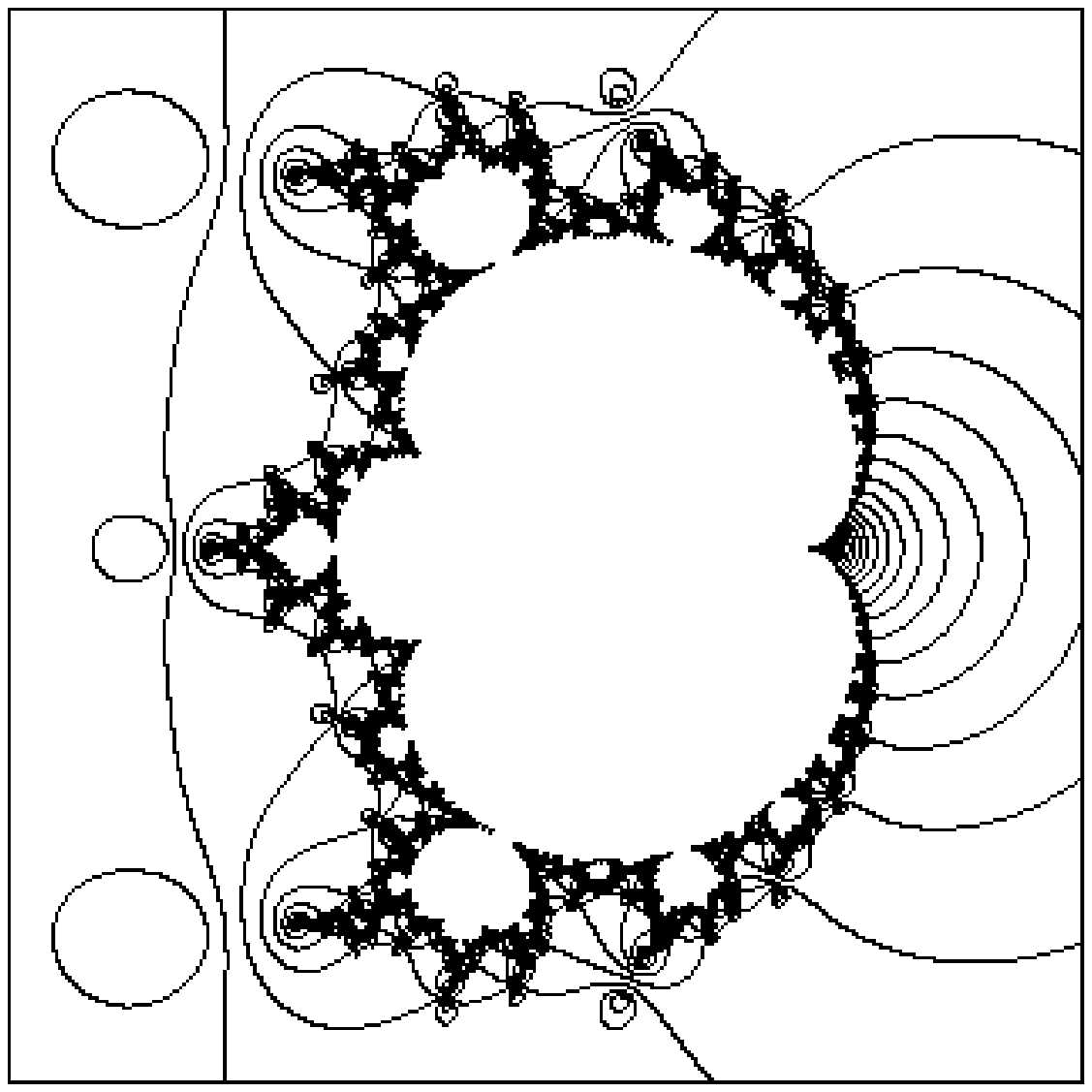,height=1.6in}
\psfig{figure=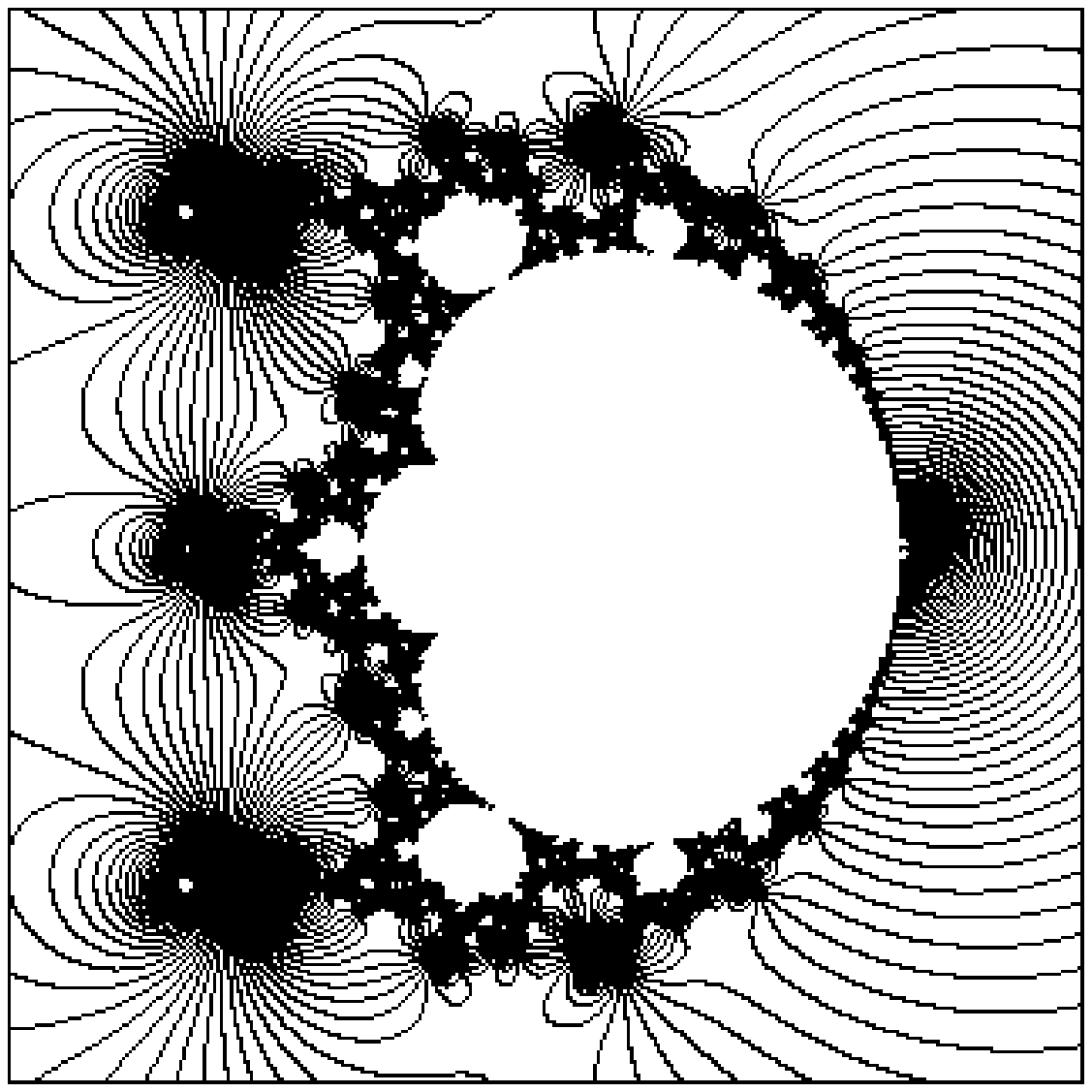,height=1.6in}}

{\ss\leftskip.25in\rightskip.25in\noindent
\bit Figure 1. Pictures of \[\CL\cap\per_1(\lambda)\] for \[\lambda=
0\,,~.01\,,~0.5\], and \[\lambda=1\] respectively, for the degree
\[n=4\]. The configuration deforms continuously for \[0<\lambda< 1\],
and conjecturally as \[\lambda\to 1\] also.
However, there is a qualitative difference between the first and second
pictures due to the fact that as \[\lambda\to 0\] the curve \[\per_1(\lambda)\]
converges not towards \[\per_1(0)\] but rather towards an \[(n-1)$-fold
branched covering of \[\per_1(0)\] 
as shown in Figure 2. Thus Figure 1b is a somewhat squashed version
of Figure 2. The surrounding curves in the 5a, 5b, 5c represent equal rates
of convergence towards the attracting fixed point for the more slowly
converging critical point. (For the corresponding curves
in Figure 1d, see 4.3 below.)\ss}
\endinsert

{\QP{\bf Theorem 4.2.} \it The intersection \[\CL\cap\per_1(1)\] is
a compact, connected full subset of the curve \[\per_1(1)\cong\C\].
Equivalently, the parabolic shift locus \[\SL_{\rm par}=\per_1(1)\ssm
(\CL\cap\per_1(1))\] is always non-vacuous, conformally isomorphic to a
punctured disk.\ms}

(Compare Figure 1d, as well as [M2, Figure 4].)

\midinsert
\cl{\psfig{figure=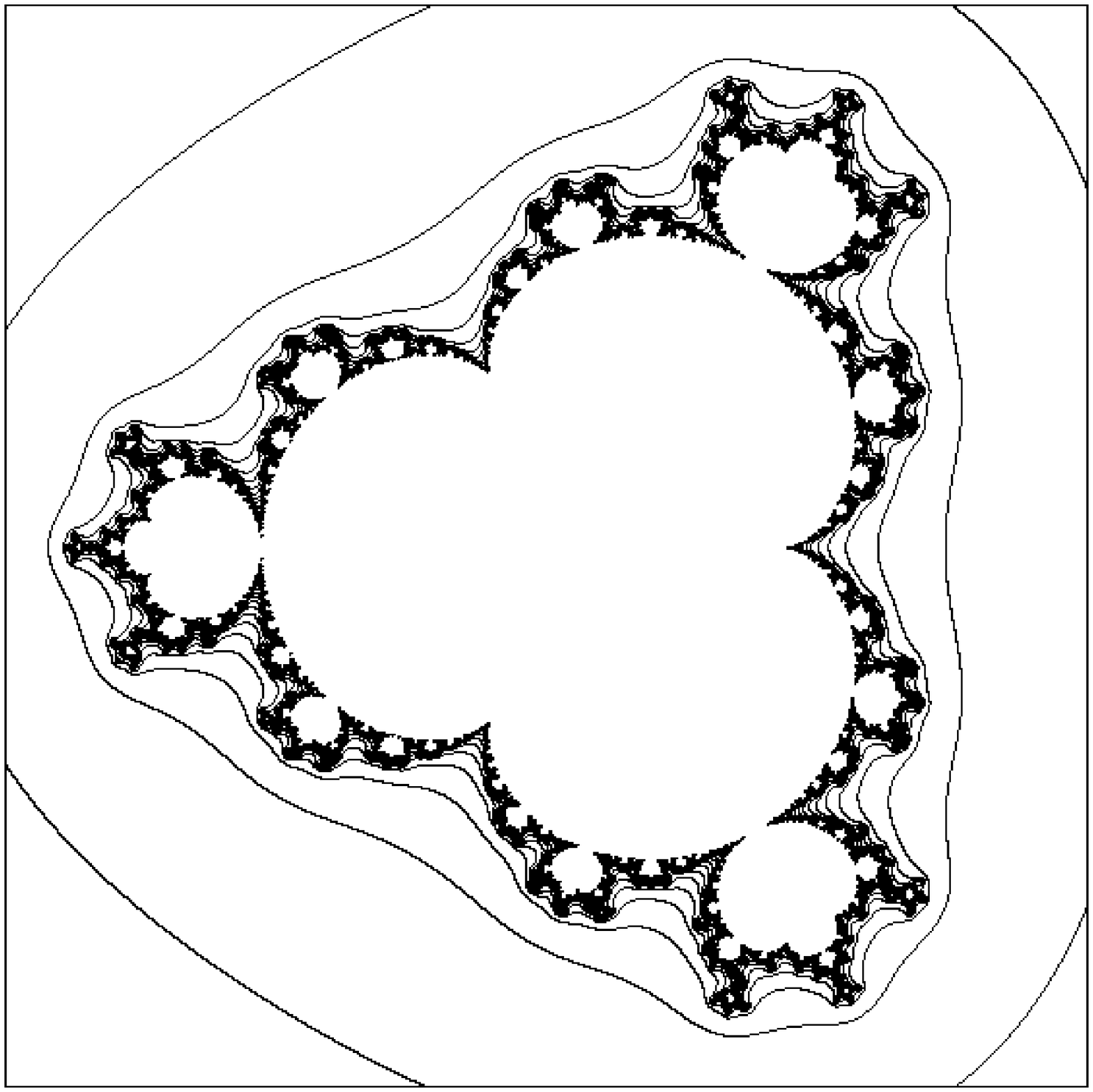,height=2in}}\ss

{\ss\leftskip.25in\rightskip.25in\noindent
\bit Figure 2. The ``Multibrot set'' for degree \[n=4\], that is the
connectedness locus in the \[b$-parameter plane for the family of unicritical
polynomial maps \[z\mapsto z^4+b\]. (See for example [LS2].)
Note the \[(n-1)$-fold rotational symmetry. The corresponding figure in the
\[b^3$-plane is shown in Figure 1a.\ss}
\endinsert


I will outline two proofs of 4.2, one suggested by conversations with
Schleicher, and one suggested by Shishikura. The first begins as follows.\ss

{\bf Proof of compactness.}
It will be convenient to use the normal form (1), choosing
the matrix of coefficients to have the form
$$\left[\matrix{a&b\cr c&d}\right]\=\left[\matrix{n+1+\alpha&n-1-\alpha\cr
	n-1+\alpha&n+1-\alpha}\right]~, \eqno (12)$$
where \[\alpha\in\C\] is a parameter. Then it is easy to
check that the associated rational function \[f\] satisfies
\[f(1)=f'(1)=1\], and that the invariant \[X=X(f)\] of \S1 is given
by \[X=((n-1)^2-\alpha^2)/4n\], so that \[|X|\to\infty\] as \[|\alpha|
\to\infty\].

{\bf Remark.} The special case \[\alpha=0\] corresponds to the unique
conjugacy class \[(f)\] such that \[f\] has a fixed point of multiplier
\[+1\] with two attracting petals. Using this normal form (12),
the corresponding Julia set is the unit circle.\ss

Assuming that \[\alpha\ne 0\], let us set
$$	z\=1+{2\over\alpha w}~,\qquad w\={2\over \alpha(z-1)}~. $$
Then a straightforward computation shows that the
map \[z\mapsto f(z)\] corresponds to
$$	w~\mapsto~F(w)\= 
{2\over\alpha f(1+2/\alpha w)-\alpha}\=
 w+1+O\Big({1\over\alpha w}\Big)~, \eqno (13)$$
where the error estimate holds uniformly provided that both \[|\alpha|\] and
\[|\alpha w|\] are sufficiently large. In particular,
it follows that the region \[\{z~;~{\rm Re}(w)>1/2\}\] maps holomorphically
into itself, and hence is contained in the Fatou set of \[f\], provided that
\[|\alpha|\] is sufficiently large. On the other hand, it is not hard to
check that
$$	f(0)\={b\over d}\=1+{2\over \alpha w_1}~,\qquad f(\infty)\=
{a\over c}\=1+{2\over\alpha w_2} $$
where
$$	w_1\=1-{n+1\over\alpha}\qquad{\rm and}
\qquad w_2\=1+{n-1\over\alpha}~. \eqno (14) $$
Thus, for \[|\alpha|\] sufficiently large, both critical values
belong to the half-plane \[{\rm Re}(w)>1/2\], and hence belong to the same
Fatou component, so that \[(f)\] belongs to the parabolic shift locus.
Therefore the closed
set \[\CL\cap\per_1(1)\] is bounded and hence compact, as asserted.\ss

{\bf Remark 4.3.} We can construct
a holomorphic function \[\Phi:\SL_{\rm par}\to\C\] as follows. For any \[(f)\in
\SL_{\rm par}\], let \[P\subset\widehat\C\ssm J_f\] be an attracting petal for the
parabolic fixed point, and let \[\phi\] be a Fatou coordinate, mapping
\[P\] biholomorphically into \[\C\], and satisfying
$$	\phi(f(z))\=\phi(z)+1~. $$
Then \[\phi\] extends canonically to a holomorphic map which carries the
entire parabolic basin onto \[\C\], satisfying this same
functional equation. In particular, if \[ c_1\] and \[ c_2\]
are the critical points, then the difference \[\phi( c_1)-\phi( c_2)\]
is a well defined complex number, independent of the choice of petal and
Fatou coordinate. In order to make this construction independent of the
numbering of the critical points, we set
$$	\Phi(f)\=\big(\phi( c_1)-\phi( c_2)\big)^2~. $$
then it is not difficult to check that
$$	\Phi\,:\,\SL_{\rm par}~\to~\C $$
is well defined and holomorphic. (Compare the construction of Fatou coordinates
as given in [Ste].) The curves \[|\Phi|={\rm constant}\] are shown
in Figure 1d, and in a much larger region of  \[\per_1(1)\] in Figure 3.
The asymptotic formula
$$	\Phi(f)~\simeq~\Big({2n\over\alpha}\Big)^2~\simeq~{-n\over X(f)}
 \eqno(15) $$
as \[|X(f)|\to\infty\] can be verified as follows. The inequality
$$	{dF(w)\over dw}\= 1~+~O\left({1\over\alpha w^2}\right) $$
follows from (13) together with the Schwarz Lemma. Hence
$$	{F(w_2)-F(w_1)\over w_2-w_1}\= 1~+~O\left({1\over\alpha w_j^2}\right)
  $$
provided that \[w_2/w_1\] is reasonably close to \[1\]. Setting
$$\phi(c_2)-\phi(c_1)\=\lim\,\big( F^{\circ m}(w_2)-F^{\circ m}(w_1)\big) $$
as in [Ste], it follows that
$$	\phi(c_2)-\phi(c_1)\= (w_2-w_1)\,(1+O(1/\alpha) )~. $$
But \[w_2-w_1=2n/\alpha\] by (14); and (15) follows.\ss

\midinsert
\cl{\psfig{figure=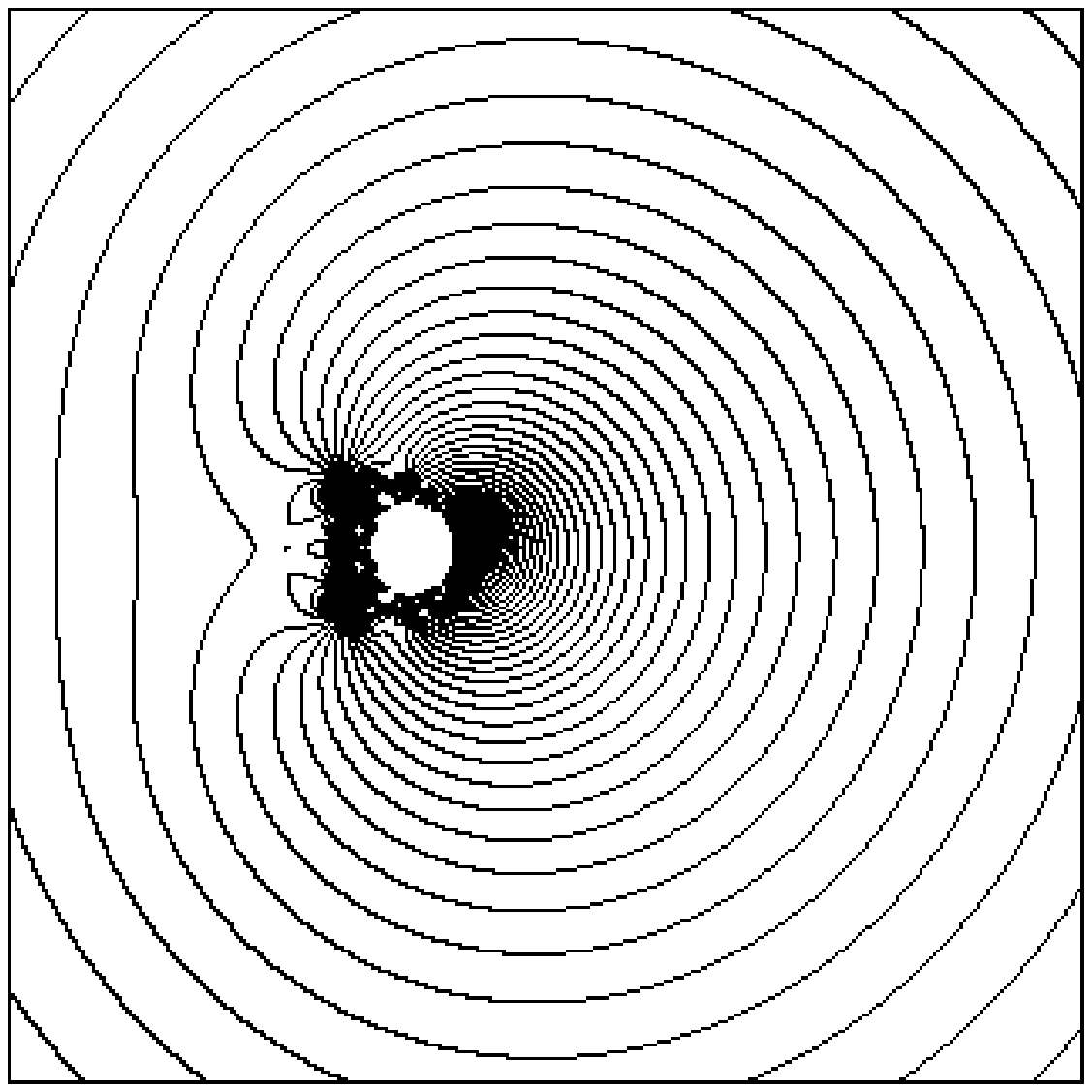,height=2in}}

{\QP\bit Figure 3. Another picture of \[\CL\cap\per_1(1)\] in the degree 4
case, showing a much larger region in order
to illustrate behavior near infinity.\ss}
\endinsert

The proof of 4.2 continues as follows.
In order to show that \[\CL\cap\per_1(1)\] is connected, we
study the limit, in the Hausdorff topology, of the intersection
\[\CL\cap\per_1(\lambda)\] as \[\lambda\] tends to \[1\] through real values
\[\lambda<1\]. First note that
$$	\limsup_{\lambda\nearrow 1}\;\CL\cap\per_1(\lambda)~\subset~
	\CL\cap\per_1(1)~.	\eqno (16) $$
Suppose that \[(f)\] can be approximated arbitrarily closely by elements
of \[\CL\cap\per_1(\lambda)\] with \[\lambda\nearrow 1\]. If \[(f)\]
did not belong to \[\CL\cap\per_1(1)\], then it would have to belong to the
parabolic shift locus. From the proof of 4.1, it would follow that any
approximating map in \[\per_1(\lambda)\] with \[\lambda<1\] must belong to
the hyperbolic shift locus, contradicting our assumption.\ss

On the other hand, we will show that
$$	\liminf_{\lambda\nearrow 1}\;\partial( \CL\cap\per_1(\lambda))~\supset
 ~\partial( \CL\cap\per_1(1))~.\eqno (17)$$
We again use the normal form (12) with marked critical points \[0\]
and \[\infty\], and with parabolic fixed point \[z=1\], writing
\[f=f_\alpha\], where \[\alpha\] is the parameter. Consider a boundary point
\[f_{\alpha_0}\] of \[\CL\cap\per_1(1)\]. We will show that
\[f_{\alpha_0}\] can be approximated arbitrarily closely by maps \[f_\alpha\]
such that one critical orbit of \[f_\alpha\] lands on a repelling
periodic orbit. To fix our
ideas, suppose that the critical point \[\infty\] lies in the parabolic basin
for \[f_{\alpha_0}\], and consider the sequence of maps
$$ \alpha~\mapsto f_\alpha^{\circ k}(0) \eqno (18) $$
for \[k=1\,,\,2\,,\,3\,,\,\cdots\], where \[\alpha\] ranges over some
neighborhood of \[\alpha_0\].

{\bf Case 1.} Suppose that the \[f_\alpha^{\circ k}(0)\] do not form a normal
family throughout any neighborhood of \[\alpha_0\]. Choose a repelling
periodic orbit for \[f_{\alpha_0}\] with period \[\ge 3\]. This orbit varies
holomorphically with the parameter \[\alpha\], throughout some neighborhood
of \[\alpha_0\]. By non-normality, we can choose \[\alpha\] arbitrarily
close to \[\alpha_0\] so that the orbit of \[0\] under \[f_\alpha\]
eventually lands on this periodic orbit. Now perturbing slightly, we can
preserve this critical orbit relation but replace the multiplier at the fixed
point \[z=1\] by some number \[\lambda<1\]. The resulting map must belong to
\[\CL\cap\per_1(\lambda)\], as required.

{\bf Case 2.} Now suppose that the \[f_\alpha^{\circ k}(0)\] do
form a normal family
throughout some neighborhood of \[\alpha_0\].
The hypothesis that every neighborhood of \[f_{\alpha_0}\] intersects the
parabolic shift locus then guarantees that this family of maps must converge
uniformly near \[\alpha_0\] to the constant map \[\alpha\mapsto 1\] as
\[k\to\infty\]. In particular, \[f_{\alpha_0}^{\circ k}(0)\] must converge
to \[1\]. If \[0\] lies in the same parabolic basin as \[\infty\], then it
follows that \[f_{\alpha_0}\] lies in the shift locus, contradicting our
hypothesis. The only other possibilities are that either:

(a) \[\alpha_0=0\], so that \[f_{\alpha_0}\] has two distinct
parabolic basins and \[z=1\] is a fixed point of higher multiplicity, or

(b) some forward image \[f_{\alpha_0}^{\circ k}(0)\] is precisely
equal to the parabolic fixed point \[z=1\].\ss

\ni In case (a), under a slight perturbation within \[\per_1(1)\] this fixed
point splits into one fixed point of multiplier \[+1\]. together with a
second fixed point which can have any multiplier close to
\[+1\]. In particular, if we perturb so that this second fixed point is
attracting, then we must be within the connectedness locus. Therefore
\[f_{\alpha_0}\] is not an isolated point of \[\CL\cap\per_1(1)\]. Hence it
is not an isolated boundary point, and, after
a slight perturbation, we can obtain a contradiction by the argument above.\ss

In case (b), there is only one parabolic basin. Suppose that \[(f_{\alpha_0})\]
were an isolated point of \[\CL\cap\per_1(1)\]. Let \[\alpha\] range over
a small circle centered at \[\alpha_0\], and assume that the corresponding
maps \[f_\alpha\] all belong to the (parabolic) shift locus. Then the
corresponding images \[f_\alpha^{\circ k}(0)\] must loop around the parabolic
fixed point \[z=1\] one or more times, without ever hitting the Julia set
\[J(f_\alpha)\]. By Ma\~n\'e-Sad-Sullivan or Lyubich, this Julia set must
vary continuously as we go around around the loop.
Choose a repelling periodic point of period \[\ge 2\] which is close enough to
\[z=1\] so that it remains inside this loop in the \[z$-plane, for all
parameter values in the circle. A priori, we might worry that this periodic
point comes back to a different periodic point as we go around the circle.
However this cannot happen since we can deform the circle in \[\per_1(1)\]
into a circle in \[\per_1(1-\epsilon)\] which deforming the Julia set
homeomorphically. The corresponding disk in \[\per_1(1-\epsilon)\] bounds
a disk in the hyperbolic shift locus, so the monodromy must be trivial.
Now shrink the parameter loop down to the point \[z=1\].
A winding number
argument shows that at some point during this shrinking, the image
\[f_\alpha^{\circ k}(0)\] must exactly hit the corresponding repelling point,
and hence belong to the Julia set. This contradicts the hypothesis that
\[(f_{\alpha_0})\] was isolated in \[\CL\cap\per_1(1)\]. But if this point
is not isolated, then we see as above that it is indeed possible to approximate
\[f_{\alpha_0}\] by a map with \[0\] eventually mapping to a repelling
periodic point.
	
Now as we vary \[\lambda=1\] to a value
slightly less than 1, the point in moduli space satisfying this critical orbit
relation, say \[f^{\circ k}(0)=f^{\circ \ell}(0)\] deforms continuously,
and necessarily belongs to the connectedness locus. This proves (17).

{\bf Remark.} In fact it is conjectured that \[\CL\cap\per_1(1)\]
is equal to the Hausdorff limit of \[\CL\cap\per_1(\lambda)\] as \[\lambda
\nearrow 1\]. However, our arguments will leave open the possibility
of a ``parabolic queer component'' in \[\CL\cap\per_1(1)\], whose points can be
approximated arbitrarily closely by points in the hyperbolic shift locus.\ss

{\bf Proof that \[\CL\cap\per_1(1)\] is connected.} For \[|\lambda|<1\],
it is known that
\[\CL\cap\per_1(\lambda)\] is compact, connected and full, with connected
boundary. (Compare [GK], [M2].)
If \[\CL\cap\per_1(1)\] were not connected, then we could choose boundary
points in two different components. By (17), these could be approximated by
points in the boundary of say \[\CL\cap\per_1(1-\epsilon)\]. Since this
is true for arbitrarily small \[\epsilon\], it would follow from (16)
that these points must actually belong to the same component.
Thus \[\CL\cap\per_1(1)\] is connected.\ss

In order to prove that \[\CL\cap\per_1(1)\] is full, or equivalently that
the parabolic shift locus is connected, we will need
a sharper form of the construction used to prove B.5.

{\QP{\bf Lemma 4.4.} \it If \[(f)\in\SL_{\rm par}\] has no critical
orbit relation, then there exists an
attracting petal for \[f\] which contains both critical values.\ss}

Here by an {\bit attracting petal\/} we mean a simply-connected open set \[P\]
which eventually captures all orbits in the parabolic basin, and
such that \[f\] maps the closure
\[\overline P\] homeomorphically, with \[f(\overline P)\subset P\cup\{\hat z\}
\] where \[\hat z\]
is the parabolic fixed point. By a {\bit critical orbit relation\/},
we mean some relation of the form \[f^{\circ k}(c_1)=f^{\circ\ell}(c_2)\].
\ss

Let us start with some petal \[P_0\], with smooth
boundary containing no points of the critical orbits. We may assume that \[P_0\] contains no critical value. Inductively
construct \[P_0\subset P_1\subset P_2\subset\cdots\] where \[P_{k+1}\]
is the connected component of \[f^{-1}(P_k)\] which contains \[P_k\].
Define the two integers \[0\le k_1\le k_2\], by setting \[k_j\] equal to
the smallest integer such that \[P_{k_j}\] contains \[j\] distinct
critical values. Then \[P_{k_1}\] is itself a petal, but \[P_{k_1+1}\] is not,
since it contains a critical point (using 3.4). The proof of 4.4 will be by
induction on the difference \[k_2-k_1\]. To start the induction, if
\[k_1=k_2\], then \[P_{k_1}\] is the required petal, and we are done.

Suppose then that \[k_1<k_2\]. Let \[v_1\] be the critical value which is
contained in \[P_{k_1}\]. By 3.4 we know that \[P_{k_1+1}\] is a simply
connected open set which contains the
corresponding critical point \[c_1\] and is a branched \[n$-fold covering
of \[P_{k_1}\]. Let \[x=f^{\circ(k_2-k_1)}(c_2)\] be the unique point in
\[P_{k_1+1}\ssm P_{k_1}\] which belongs to the second critical orbit.
Choose some path \[\gamma\] within \[\overline P_{k_1}\ssm\overline P_{k_1-1}\]
which joins the critical value \[v_1\] to the boundary of \[P_{k_1}\],
and which avoids the point \[f(x)=f^{\circ(k_2-k_1)}(v_2)\]. Then the preimage
of \[\gamma\] under \[f\] is a union \[\gamma_1\cup\cdots\gamma_n\] of
\[n\] paths, each joining the critical point \[c_1\] to the boundary of
\[P_{k_1+1}\]. These \[n\] paths cut the open set \[P_{k_1+1}\] into \[n\]
regions, each of which maps diffeomorphically onto \[P_{k_1}\ssm\gamma\].
Exactly one of these \[n\] regions contains the point \[x\], and exactly
one of these \[n\] regions contains \[P_{k_1}\].

{\QP ASSERTION: {\it It is possible to choose the path \[\gamma\] so that
\[x\] and \[P_{k_1}\] belong to the same connected component of \[P_{k_1+1}
\ssm(\gamma_1\cup\cdots\cup\gamma_n)\].\ss}}

\midinsert
\cl{\psfig{figure=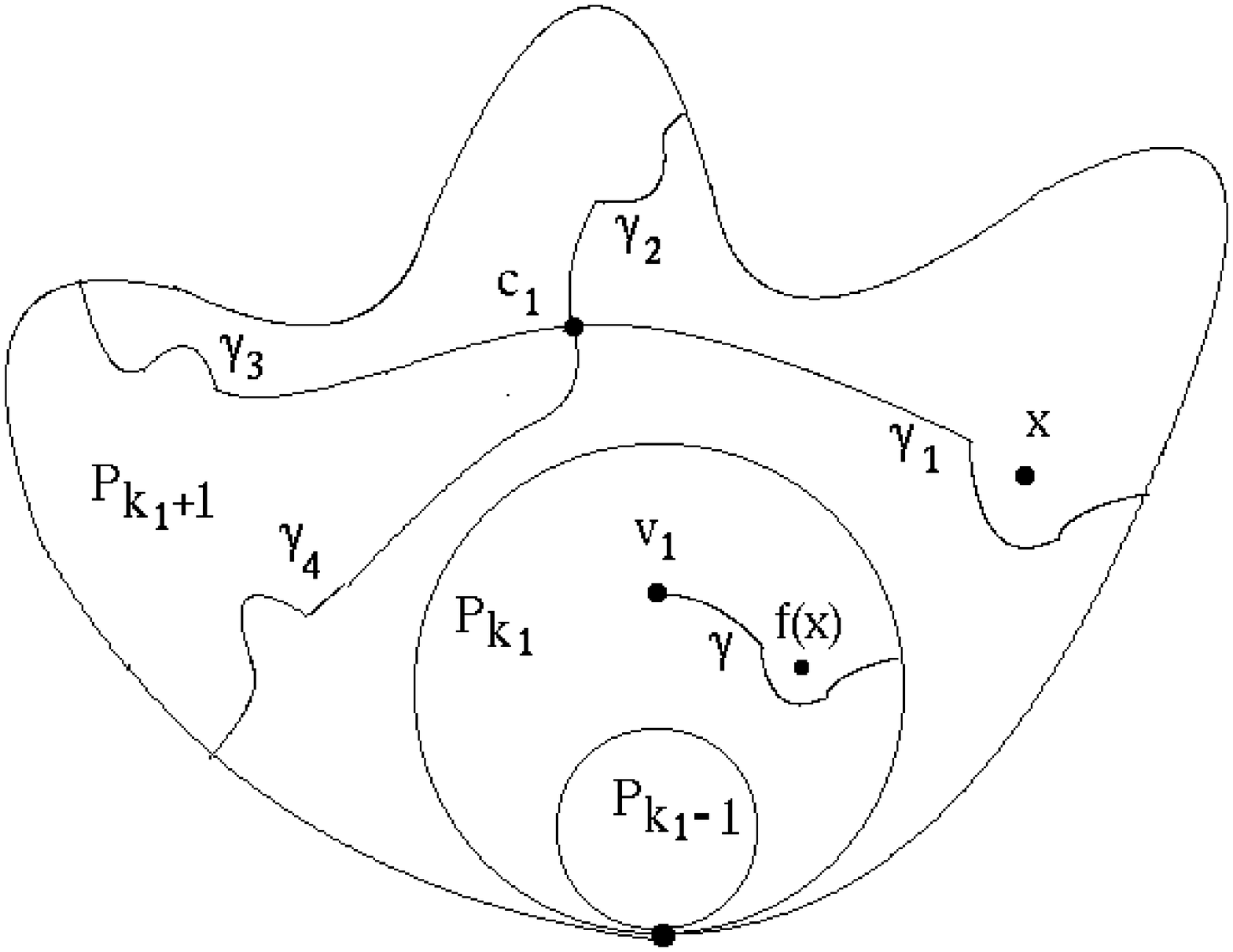,height=2.4in}}
\ss
{\QP\bit Figure 4.  This shows the three petals \[P_{k_1-1}\,,\;P_{k_1}\],
and \[P_{k_1+1}\]. The point \[x=f^{\circ(k_2-k_1)}(c_2)\] belongs to the
orbit of the second critical point. A path \[\gamma\]
from \[v_1\] to \[\partial P_{k_1}\]
within \[P_{k_1}\ssm(\overline P_{k_1-1}\cup\{f(x)\})\] lifts to \[n\]
distinct paths \[\gamma_i\] from \[c_1\] to \[\partial P_{k_1+1}\]. We must
choose this path
\[\gamma\] so that \[x\] and \[P_{k_1}\] belong to the same connected
component of \[P_{k_1+1}\ssm(\gamma_1\cup\cdots\cup \gamma_n)\].\par}
\endinsert

As an example, in Figure 4 this requirement fails as the figure is drawn.
However, if we modify the path \[\gamma\] in a neighborhood of \[f(x)\]
so that it passes above \[f(x)\] rather than below, then the requirement
will be satisfied. More generally, we can choose the path \[\gamma\] from
\[v_1\] so as to loop any number of times around an arc joining \[v_1\] to
\[f(x)\] before terminating on \[\partial P_{k_1}\]. By choosing the number
of loops appropriately, we can easily guarantee that \[x\] lies in the required
component of \[P_{k_1+1}\ssm(\gamma_1\cup\cdots\cup\gamma_n)\]. Details
will be left to the reader.\ss

The proof of 4.4 now proceeds as follows. Construct a new attracting petal
\[P'_0\] which contains no critical value by removing a thin neighborhood of
\[\gamma\] from \[P_{k_1}\]. Then the preferred component \[P'_1\]
of \[f^{-1}(P'_0)\] will consist of the component of
$$  P_{k_1+1}\ssm({\rm neighborhood~of~}\big(\gamma_1\cup\cdots\cup\gamma_n)
\big)$$
which contains \[P_{k_1}\]. By the construction, both \[v_1\] and \[x=
f^{\circ(k_2-k_1)}(c_2)\] belong to this set \[P'_1\]. Therefore,
the new difference \[k'_2-k'_1\] will be equal to \[k_2-k_1-1\]. The conclusion
of 4.4 now follows by induction.\QED
\ms

{\bf Proof of 4.2 (conclusion).}
We must show that the parabolic shift locus is connected. The proof
will make use of a standard quasiconformal surgery argument, as
suggested to me by D. Schleicher. After a small perturbation of \[f\], we
may assume that there are no critical orbit relations.
Choose a petal \[P\] as in 4.4, and choose
an embedded disk \[\overline\Delta\subset P\ssm\overline{f(P)}\] which
contains both critical
values in its interior. Now choose a diffeomorphism from \[\Delta\] to itself
which is the identity near the boundary and which moves the critical
value \[v_2\] arbitrarily close to \[v_1\]. Pulling the standard
conformal structure back under this diffeomorphism, we obtain a new
conformal structure on \[P\ssm\overline{f(P)}\]. Now push this new
structure forward under the various iterates of \[f\], and also pull it
backwards under the various iterates of \[f^{-1}\]. Since the Shishikura
Condition is satisfied: {\it every orbit intersects \[\Delta\] at most once,\/}
it follows that we obtain a well defined measurable conformal structure,
invariant under \[f\]. By the measurable Riemann Mapping Theorem, there is
a quasiconformal mapping \[\Phi\] which conjugates this new structure to the
standard conformal structure. Now \[g=\Phi\circ f\circ\Phi^{-1}\] is a new
rational map, topologically conjugate to \[f\]. By the construction, the
two critical values of \[g\] are close together in an embedded disk
\[\Phi(\Delta)\] which contains no critical point. By 1.8, it follows that
the invariant \[|A(g)|\] is arbitrarily close to infinity.
Thus \[g\] belongs to the unbounded component of the parabolic shift
locus. Since \[f\] and \[g\] clearly belong to the same connected
component of \[\SL_{\rm par}\], this shows that \[\SL_{\rm par}\]
is connected.\QED
\bs

\midinsert
\cl{\psfig{figure=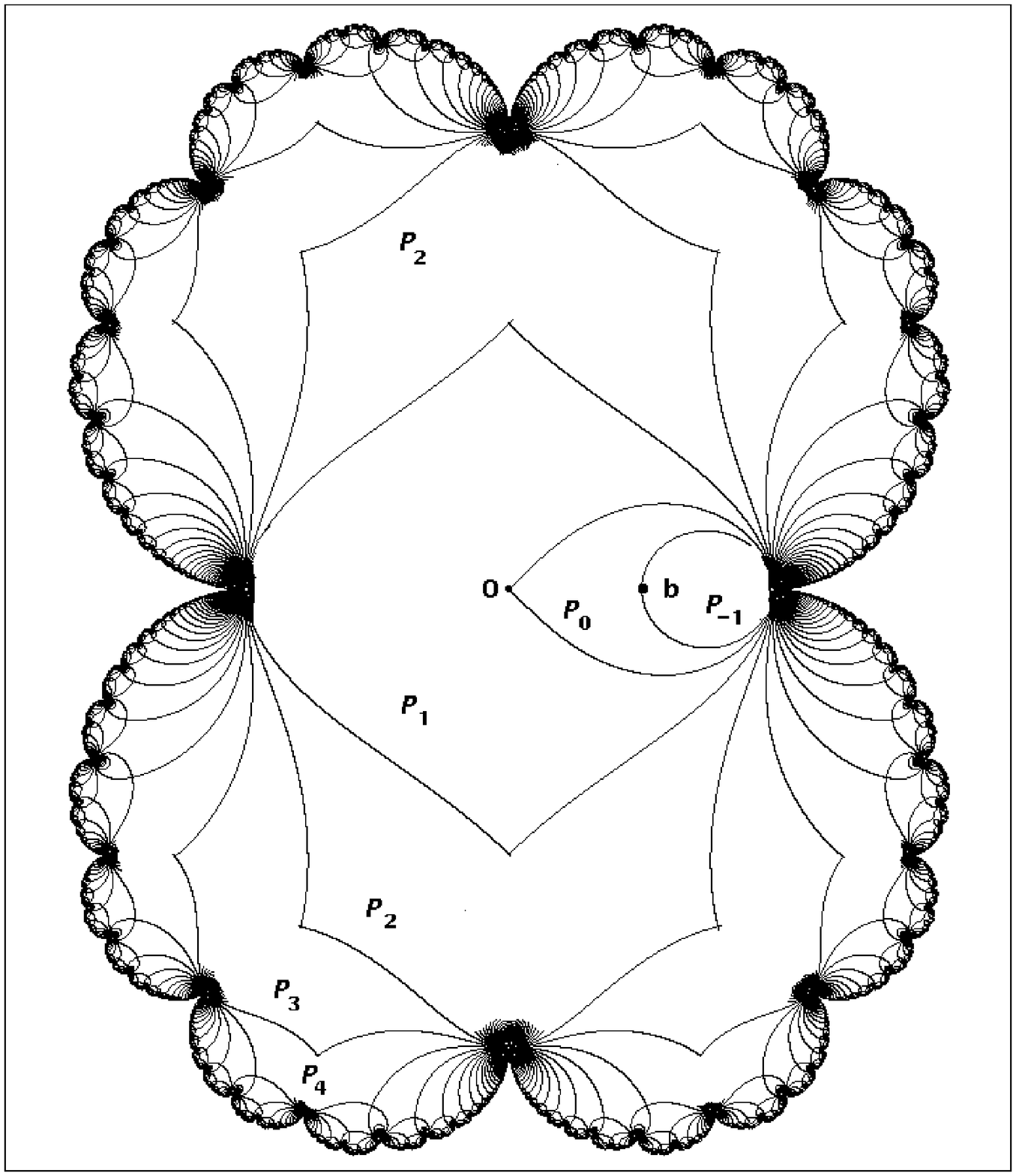,height=4.17in}}\ss
{\QP\bit Figure 5. The regions \[{\cal P}_0\subset{\cal P}_1\subset\cdots\]
for the quadratic case \[n=2\]. (Note: In
the region outside of \[{\cal P}_4\], all of the iterated preimages
of \[\partial{\cal P}_0\] have been drawn in, so that this figure looks
more complicated near the basin boundary.)\ss}
\endinsert

{\bf Alternative Proof of 4.2.} I am grateful to
Shishikura for suggesting a quite different argument,
based on an explicit model for the parabolic shift locus which can
be outlined as follows. Start with the map
$$	f(z)\=z^n+b\qquad{\rm with}\qquad b={(n-1)/ n^{n/(n-1)}}~, $$
which has a parabolic fixed point at \[z=(1/n)^{1/(n-1)}\]. The corresponding
parabolic basin \[\cal B\] is connected and simply-connected, and contains
a single critical point \[0\mapsto b\]. Let \[{\cal P}_0
\subset{\cal B}\] be the largest attracting petal such that the Fatou
coordinate map carries \[{\cal P}_0\] diffeomorphically onto a right
half-plane. Thus the critical point \[0\] belongs to the
boundary of \[{\cal P}_0\], and the critical value \[b\] belongs to the
boundary of the smaller petal \[{\cal P}_{-1}=f({\cal P}_0)\]. Our model space
$$	{\cal S}^*_{\rm par}\=
 ({\cal B}\ssm({\cal P}_{-1}\cup\{b\}))/\alpha~, $$
conformally isomorphic to a punctured disk, is obtained by removing the subset
\[{\cal P}_{-1} \cup\{b\}\] from the basin \[\cal B\], and then
gluing the two halves of \[\partial {\cal P}_{-1}\ssm\{b\}\]
onto each other under the correspondence \[\alpha:z\mapsto \overline z\].
To each point \[v\] in this model space, we construct a multiply connected
parabolic basin \[{\cal B}'_v\] with two critical points as follows. The
original
basin \[\cal B\] can be described as the union of open subsets
$$  {\cal P}_0\,\subset\,{\cal P}_1\,\subset\,{\cal P}_2\,\subset\,\cdots~,$$
where each \[{\cal P}_k\] is the interior of a region bounded by a Jordan
curve, and
where \[f\] maps each \[{\cal P}_{k+1}=f^{-1}({\cal P}_k)\] onto \[{\cal P}_k\]
by an \[n$-fold cyclic covering, branched only over \[b\]. Let \[k_0\] be the
smallest index such that \[v\in{\cal P}_{k_0}\]. Construct a new family
$$	{\cal P}'_0\,\subset\,{\cal P}'_1\,\subset\,{\cal P}'_2\,\subset
\,\cdots $$
as follows. Let \[{\cal P}'_k={\cal P}_k\] for \[k\le k_0\], but let
\[{\cal P}'_{k+1}\] be the \[n$-fold cyclic covering of \[{\cal P}'_k\],
branched over both \[b\] and \[v\] for \[k\ge k_0\]. The covering should
extend over the boundary, so that each boundary curve of \[{\cal P}'_k\]
for \[k\ge k_0\] is covered by \[n\] distinct boundary curves for
\[{\cal P}'_{k+1}\]. Then the
inclusion \[{\cal P}'_{k-1}\hookrightarrow {\cal P}'_k\] lifts inductively
to an inclusion \[{\cal P}'_{k}\hookrightarrow {\cal P}'_{k+1}\]. Let
\[{\cal B}'_v\] be the union of the \[{\cal P}'_k\]. Note that there is a
canonical parabolic map \[f_v\] from this Riemann surface \[{\cal B}'_v\]
to itself, carrying each \[{\cal P}'_{k+1}\] onto \[{\cal P}'_k\] by an
\[n$-fold branched covering.

In the special case where \[v\in\partial{\cal P}_{-1}\ssm\{b\}\], the two
critical values \[b\] and \[v\] are on the boundary of the same petal
\[{\cal P}_{-1}\], and there is a canonical isomorphism \[{\cal B}'_v\to
{\cal B}'_{\overline v}\] which carries \[v\] to \[b\] and \[b\] to \[\overline
v\]. Hence we identify \[{\cal B}'_v\] with \[{\cal B}'_{\overline v}\] in
this special case.

Conversely, suppose that we start with a point \[(g)\] in the parabolic
shift locus. Let \[{\cal Q}_0\] be the largest petal for \[g\] such that
the Fatou coordinate maps \[{\cal Q}_0\] isomorphically onto a right
half-plane. Then \[\partial{\cal Q}_0\]
contains at least one critical point \[c_0\], and there is a unique conformal
isomorphism \[\phi_0:{\cal P}_0\to{\cal Q}_0\] which conjugates \[f\] to \[g\]
and (extended over the boundary) carries \[0\] to \[c_0\]. Now
let \[{\cal Q}_k=g^{-k}{\cal Q}_0\], and let \[k_0\ge 0\] be
the smallest index such that \[{\cal Q}_{k_0}\] contains both critical values
of \[g\]. Then \[\phi_0\] extends uniquely to a conformal isomorphism
\[\phi_{k_0}:{\cal P}_{k_0}\to{\cal Q}_{k_0}\], still conjugating \[f\] to
\[g\]. Let \[v\] be the unique
point of \[{\cal P}_{k_0}\] which maps to the second critical value under
\[\phi_{k_0}\]. In this way we obtain a holomorphic map
\[(g)\mapsto v\] from the parabolic
shift locus \[{\cal S}_{\rm par}\] to the model
space \[{\cal S}^*_{\rm par}\]. (Furthermore, it is not hard to check that
\[\phi_{k_0}\] extends uniquely to a conformal isomorphism \[{\cal B}'_v\to
\widehat\C\ssm J(g)\] which conjugates \[f_v\] to \[g\].)

We must show that this correspondence \[(g)\mapsto v\] is a proper map. That
is, if \[(g_j)\] is any sequence of points in \[{\cal S}_{\rm par}\]
which has no accumulation point within \[{\cal S}_{\rm par}\], then we must
show that the corresponding sequence \[v_j\] has no accumulation point within
the model space. First suppose that the \[(g_j)\] converge, within
\[\per_1(1)\], to a point of the connectedness locus. If the corresponding
sequence of points \[v_j\] converged to a limit within the model space,
then the corresponding sequence \[f_{v_j}\] of model parabolic basins would
possess annuli with moduli bounded away from zero which separate
the two critical points from the Julia set, yielding a contradiction.
Similarly, suppose that the \[(g_j)\] diverge within \[\per_1(1)\] to
the point at infinity. Using 1.8, we see that the distance between the
two critical values of \[f_{v_j}\] must tend to zero. In other words, \[v_j\]
must converge towards the puncture point \[b\] in the model space.
Thus our correspondence \[(g)\mapsto v\] from \[{\cal S}_{\rm par}\] to
\[{\cal S}^*_{\rm par}\] is proper. In order to show that it
is a conformal isomorphism, we need
only check that it has degree one, which follows easily from the
asymptotic formula (15).\QED\ss

{\bf Remark.} The argument shows that the model space obtained from Figure 5
is ``inside out'' with respect to Figures 1d and 3. The puncture point
\[b\] in Figure 5 corresponds to the point at infinity in the earlier
figures.\bs\bs


\def\bic{{\bf bicrit}}
\def\per{{\rm Per}}
\def\D{{\bf D}}

\def\CL{{\cal C}}
\def\SL{{\cal S}}
\def\bic {{\bf Bicrit}}
\def\PER{{\bf Per}}
\def \iso{{\,\buildrel \simeq\over\longrightarrow\,}}

\headline={%
    \ifnum\pageno > 1 {%
        \hss\tenrm\ifodd\pageno 5. REAL MAPS\hss\folio
        \else\folio\hfill MAPS WITH TWO CRITICAL POINTS\hfill\fi}%
    \else \hss\vbox{}\hss
    \fi}

\cl{\bf\S5. Real Maps}\ms

\def\A{\alpha} 

This section will study bicritical maps \[f:\widehat\C\to\widehat\C\] such that
the invariants \[X\] and \[Y\] of \S1 are both real. First a definition:

A homeomorphism \[{\A}:\widehat\C\to\widehat\C\] is {\bit antiholomorphic\/}
if it has the form \[\A(z)=\phi(\overline z)\], where \[\phi\] is some
M\"obius automorphism of \[\widehat\C\] and \[\overline z\] is the complex
conjugate,
and is called an {\bit involution\/} if \[\A\circ\A\] is the identity map.

{\QP{\bf Lemma 5.1.} \it The invariants \[X(f)\] and \[Y(f)\] are both real
if and only if the bicritical map \[f:\widehat\C\to\widehat\C\] commutes with
some antiholomorphic involution \[\alpha\]. This \[\alpha\] is unique
if and only if \[(f)\] lies off the symmetry locus (or off
the real part \[\Sigma_\R\] of the symmetry locus).\ss}

\def\c{{\bf c}}

{\bf Proof.} Let us temporarily introduce the notation \[\c(z)=\overline z\]
for complex conjugation. Then the rational map \[g=\c\circ f\circ\c\]
can be obtained from \[f\] by replacing all of its coefficients
by their complex
conjugates. Evidently \[X(g)=\overline{ X(f)}\] and \[Y(g)=\overline{ Y(f)}\].
Therefore, if \[X\] and \[Y\] are real, then it follows that \[f\] is
holomorphically conjugate to \[g\], say \[f\circ\phi=\phi\circ g\].
Setting \[\alpha=\phi\circ\c\], this means that
$$	f\circ\alpha\=f\circ\phi\circ\c\=\phi\circ g\circ\c\=\phi\circ\c\circ f
\=\alpha\circ f~, $$
so that \[f\] commutes with \[\A\]. It follows that \[f\] commutes with
\[\A\circ\A\], which is a holomorphic map from \[\widehat\C\] to itself.
If \[(f)\] lies off of the symmetry locus of 1.4, this proves that \[\A\circ
\A\] is the identity map, as required. Even for a generic \[f\] in the
symmetry locus which commutes with a unique \[\iota\],
since \[\A\] must either fix or interchange the two critical
points, it follows that \[\alpha\circ\alpha\] must fix both critical points,
and hence be the identity map. Finally, in the two exceptional cases,
corresponding to \[f(z)=z^n\] or \[f(z)=1/z^n\], the assertion is clearly
satisfied. Conversely, if \[f\] commutes with \[\alpha=\phi\circ\c\], then
composing the equation \[f\circ\phi\circ\c=\phi\circ\c\circ f\] on the
right with \[\c\], we see that \[f\circ\phi=\phi\circ g\]. Therefore
\[X(f)=X(g)\] and \[Y(f)=Y(g)\] are real, as asserted.

If there is a second antiholomorphic map \[\beta\] commuting with \[f\],
then \[\alpha\circ\beta\] is a holomorphic map commuting with \[f\], so that
\[f\] belongs to the symmetry locus. Conversely, if
\[(f)\] does belong to the symmetry locus, let \[\iota\]
be a holomorphic involution commuting with \[f\]. Then \[\A\circ\iota\circ\A\]
is also a holomorphic involution commuting with \[f\]. In the generic
case where \[\iota\] is unique, it follows that \[\A\] commutes with \[\iota\],
and that \[\A\circ\iota\] is another antiholomorphic involution commuting
with \[f\]. In other words, \[f\] can be given the structure of a ``real map''
in two essentially distinct ways. The discussion in the two exceptional
cases where \[\iota\] is not unique is straightforward, and will be left
to the reader.\QED\ss

We must now distinguish two different cases. The antiholomorphic involution
\[\A\] may have a circle of fixed points (for example the unit circle if
\[\A(z)=1/\overline z\]). In this case, a M\"obius
change of coordinates will reduce it to the standard form
\[\A(z)=\overline z\] with fixed point set \[\R\cup\{\infty\}\],
and it follows easily that the corresponding map \[f\]
can be expressed as a rational map with real coefficients. On the other hand,
it may happen that \[\A\] has no fixed points at all. In this case, a
M\"obius change of coordinates reduces \[\A\] to the standard {\bit antipodal
map\/} \[\A(z)=-1/\overline z\]. It follows from a theorem
of Borsuk and Hopf that a map \[f\] which commutes with this antipodal
map necessarily has odd degree. (Compare [AH].) Such an \[f\] gives rise to
a conformal self-map with only one critical point on
the non-orientable ``Klein surface''
which is obtained from \[\widehat\C\] by identifying each \[z\] with
\[-1/\overline z\].\ms

\midinsert
\cl{\psfig{figure=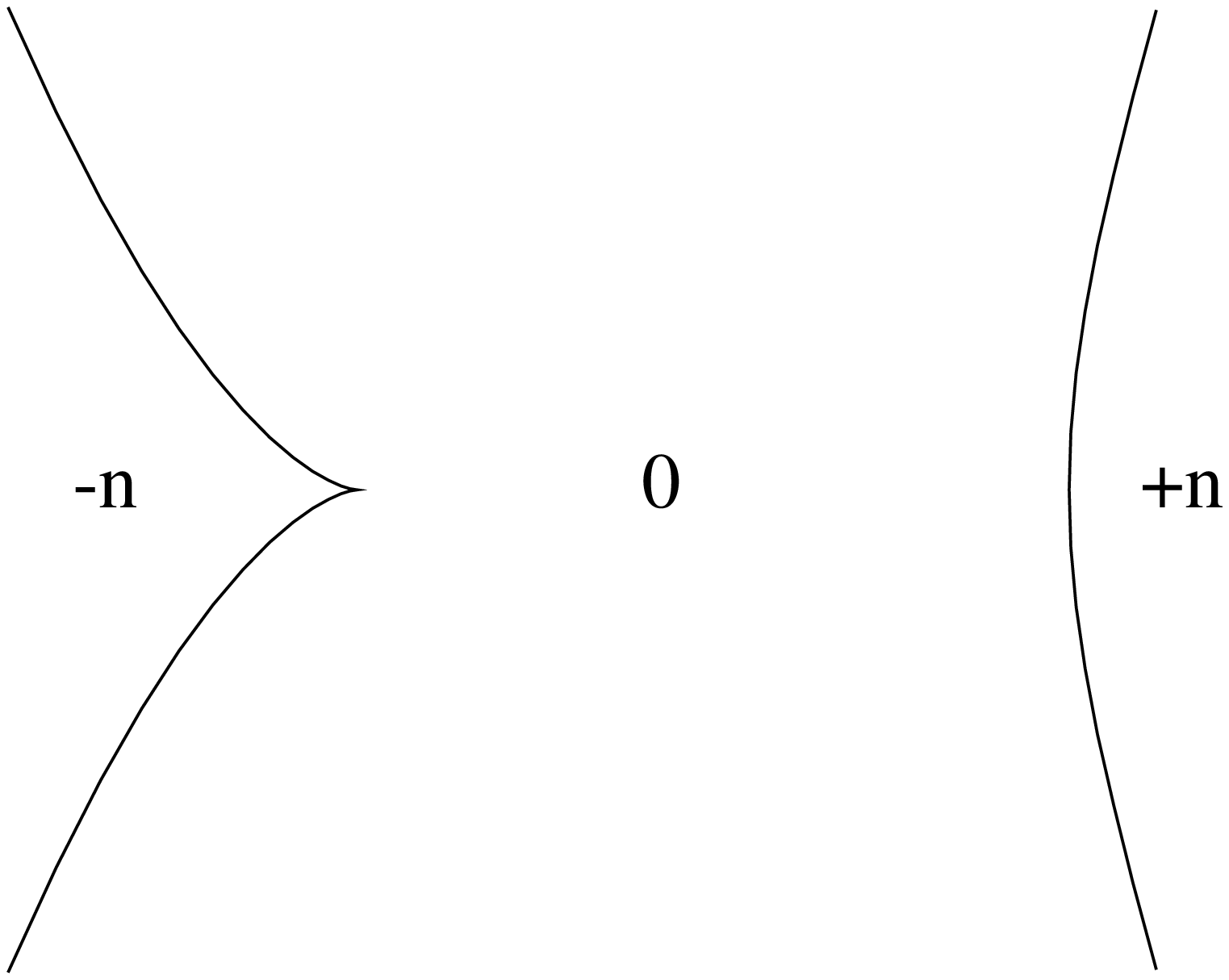,height=1.6in}}
\bs
{\QP\bit Figure 6. Symmetry locus in the real \[(X,Y)$-plane for degree
\[n=2\]. The three complementary domains are labeled according to the degree
of the associated map from the circle of real points to itself. The picture
for higher even values of \[n\] would have another cusp, corresponding
to \[(z\mapsto z^n)\], but would otherwise be similar.\par}
\endinsert

\pageinsert
\cl{\psfig{figure=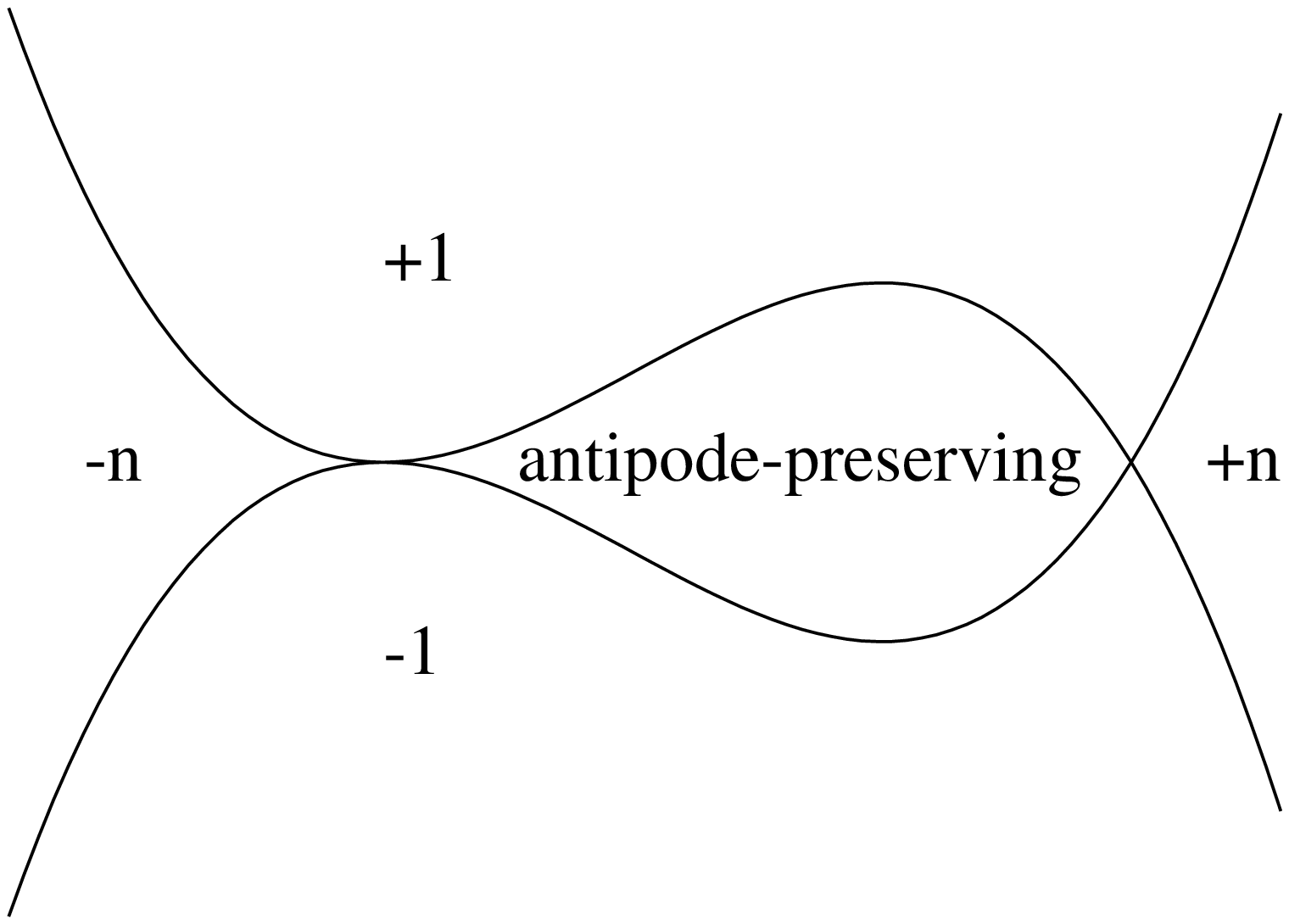,height=2.2in}}
\ss
{\QP\bit Figure 7. Symmetry locus in the real \[(X,Y)$-plane for degree
\[n=3\], with the five complementary domains appropriately labeled.\bs}

\cl{\psfig{figure=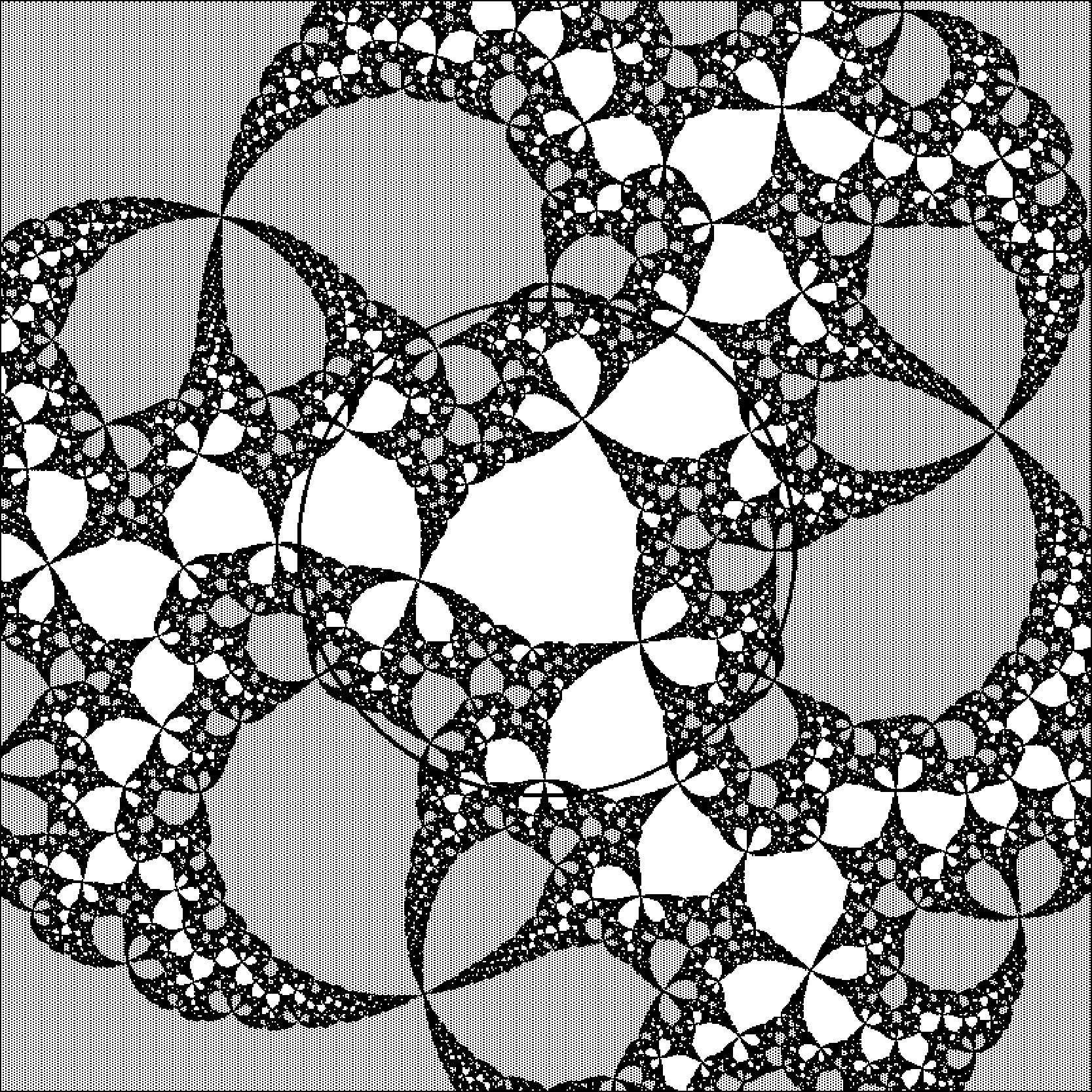,height=3.5in}}

{\QP\bit Figure 8. Dynamic plane for a degree \[n=3\]
bicritical map which commutes
with the antipodal map, with the unit circle drawn in.
(Parameters \[X=-.235\,,~Y=.213\].)
There are two attracting period 4 orbits. The
antipodal map of \[\widehat\C\], reflection in the unit circle
composed with the \[180^\circ\] rotation, interchanges the two
attracting basins, which are colored white and grey
respectively. The Julia set, like that for an
arbitrary degree \[n\] bicritical map, has \[n$-fold rotational symmetry.
It is conjectured that this particular map can be obtained as a mating
of the form \[(z^3+b)\perp\!\!\perp(z^3-\overline b)\], with
\[b\approx.584-.270\,i\]. $($Compare $[T].)$ \ss}
\endinsert

\pageinsert
\cl{\psfig{figure=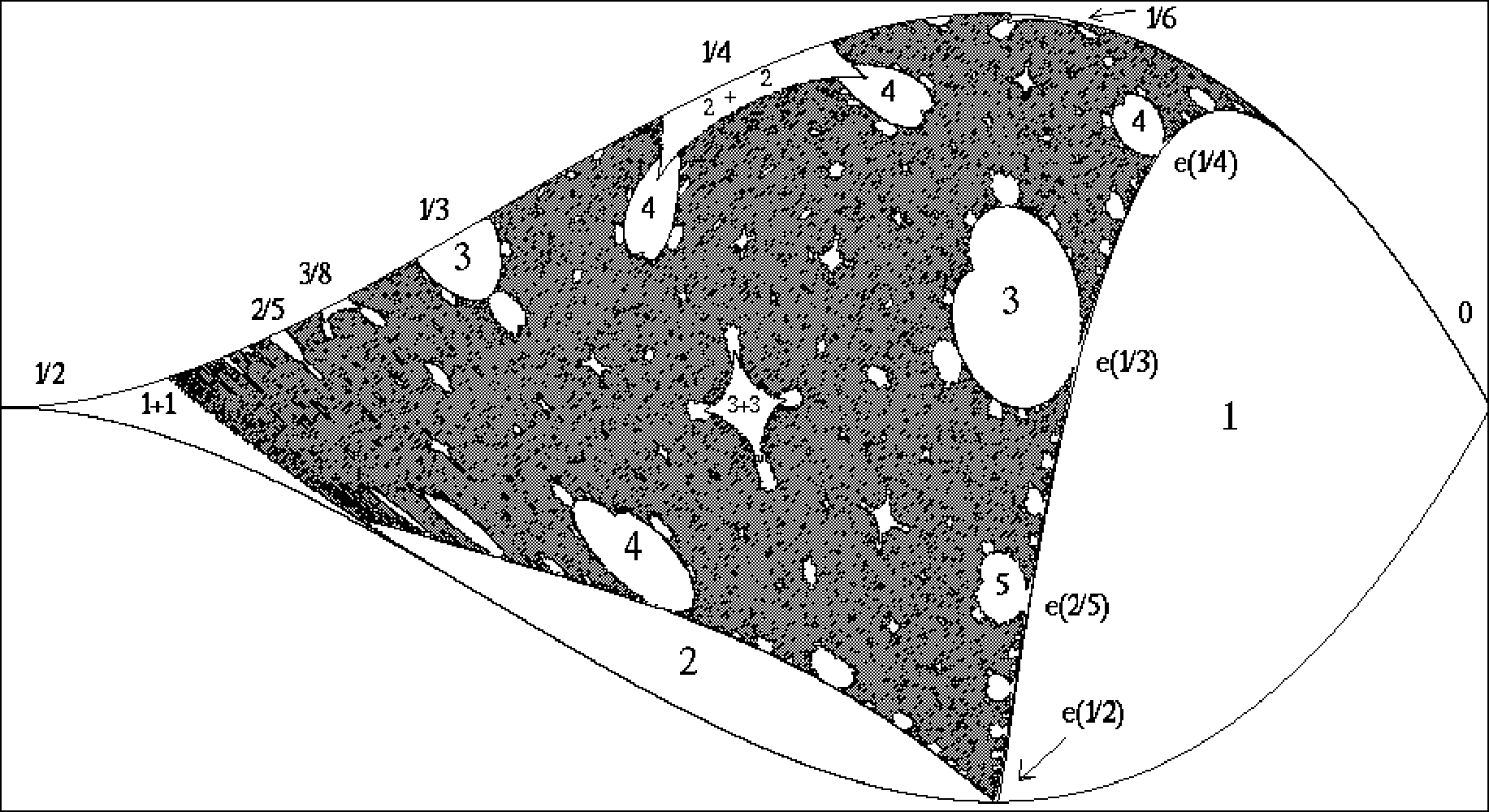,height=3.2in}}
\ms
\cl{\psfig{figure=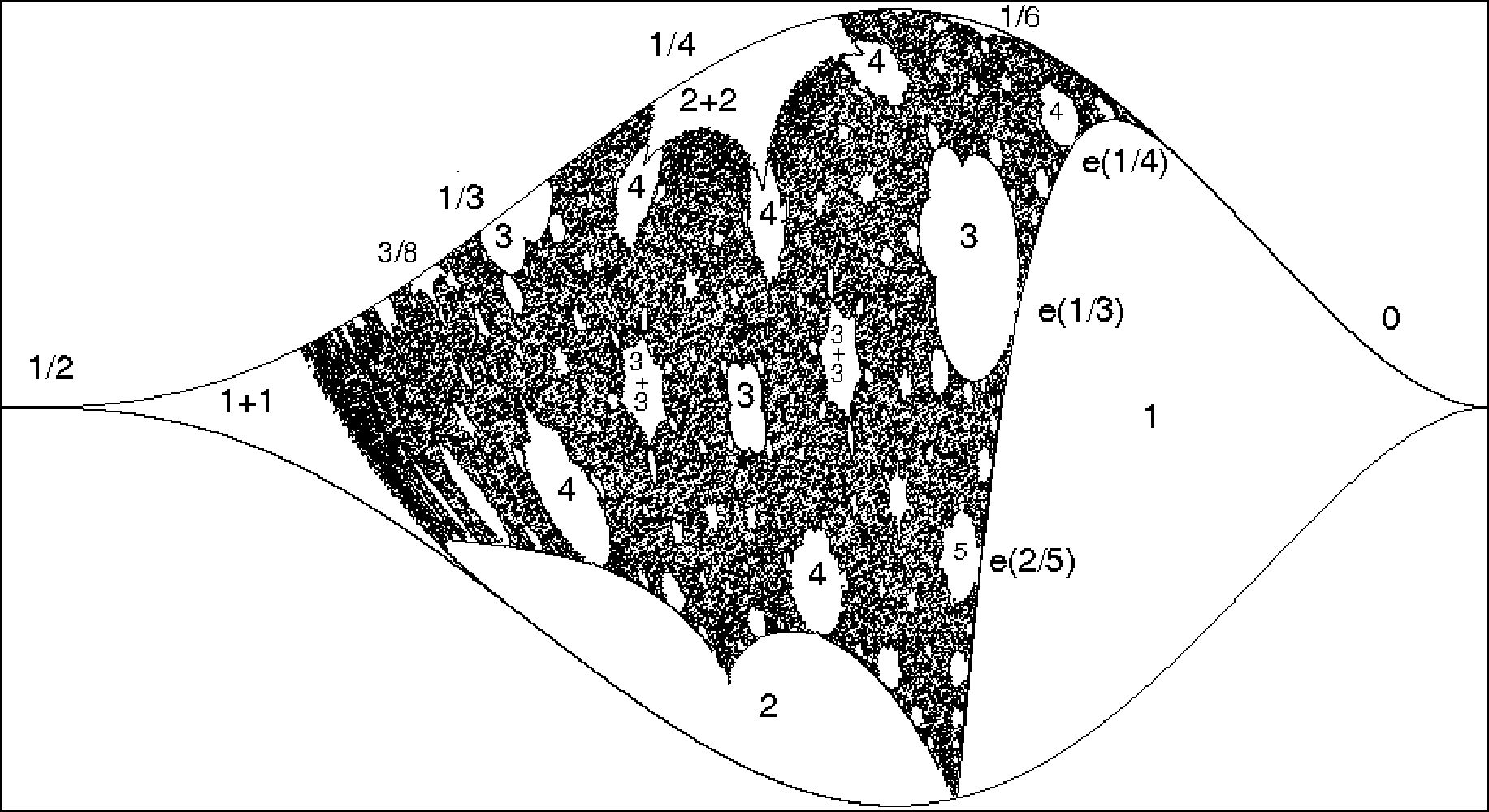,height=3.2in}}

{\QP\bit Figure 9. Antipode-preserving region of the \[(X,Y)$-plane for
degrees 3 and 5, with periods of the larger hyperbolic components labeled.
(Note that the components of type \[p+p\] have a tricorn-like
geometry. Compare [NS]) Along the top edge, these maps give rise to degree one
circle homeomorphisms, with rotation number as indicated. Also, along the
edge of the principal hyperbolic component, the multiplier of one of the
two indifferent fixed points is indicated, using the notation \[e(x)=e^
{2\pi ix}\].\ss}
\endinsert

More precisely, we can distinguish six cases as follows. Suppose
that \[f\] commutes with the antiholomorphic involution \[\A\]. For the first
five cases we suppose that \[\A\] has a circle \[S^1\] of fixed points.
Note that \[f\] induces a map from this circle to itself.\ss

{\bf Case \[+n\].} For any degree \[n\], it may happen that \[f\] induces
an \[n$-fold covering map of degree \[+n\] from the circle \[S^1\] to itself.
In this case the Julia set \[J(f)\]
either coincides with \[S^1\], in which case \[(f)\] has just two Fatou
components, each mapped onto itself by \[f\], or else \[J(f)\] is
a Cantor subset of \[S^1\], so that \[(f)\] belongs
to the shift locus.\ss

{\bf Case \[-n\].} If \[f\] induces an \[n$-fold covering of degree \[-n\],
then again the Julia set may be \[S^1\] or a Cantor subset. Correspondingly,
\[f\] either interchanges the two Fatou components,
or belongs to the shift locus.\ss

{\bf Case \[0\].} If the degree \[n\] is even, then the only other possibility
is that \[f\] maps \[S^1\] onto a proper subset of
itself by a map of degree zero. In this case,
the forward orbits of the two critical points lie in the image \[f(S^1)\],
which is an interval \[I\subset S^1\]. (Compare [M2].) Restricting \[f\]
to this interval, we obtain a map which is either monotone, or unimodal, or
bimodal of a rather restricted type since each point of \[I\] has at most two
preimages in \[I\]. Thus the dynamics is largely controlled by the theory
of smooth interval maps. In particular, any attracting or parabolic cycle
must be contained in \[I\], and if there is no such cycle then the Julia
set must be the entire Riemann sphere.\ss

Note that we can pass between these three cases only by crossing the
symmetry locus. (Compare Figure 6.)\ss

For the remaining three cases, we assume that the degree \[n\] is odd.\ss

{\bf Case \[+1\].} For any odd \[n\], the circle \[S^1\] may map to itself
by a homeomorphism of degree \[+1\] with two critical inflection points.
The dynamics of the critical orbits is then governed by the theory of
monotone degree one circle maps, with a rotation number which is well
defined up to sign. (Compare [BB].)
If the rotation number is \[p/q\], then
both critical orbits converge to (the same or different)
parabolic or attracting cycles of period
\[q\], while if the rotation number is irrational then the Julia set is the
entire Riemann sphere.\ss

{\bf Case \[-1\].} Similarly \[S^1\] may map to itself by a homeomorphism
of degree \[-1\]. In this case, both critical orbits must converge to
parabolic or attracting cycles of period one or two.\ss

{\bf Antipode Preserving Case.} If
$$	f(-1/\overline z)\=-1/\overline{f(z)}~, $$
then
\[(f)\] necessarily belongs to the connectedness locus. For every hyperbolic
component in this region of the \[(X,Y)$-plane, there are either two
antipodal attracting orbits of the same period \[p\], or else a single
attracting orbit of period \[2p\] satisfying the identity
\[\A(x)=f^{\circ p}(x)\]. In the latter case it is convenient to say that
the hyperbolic component is of type \[p+p\].\ss

The division of the real \[(X,Y)$-plane into five regions  for a typical
odd \[n\] is shown in Figure 7.
The antipode-preserving region in this plane is shown in detail in Figure 9
for two values of \[n\], and a typical  associated Julia set is shown in
Figure 8.\bs\bs

\headline={%
    \ifnum\pageno > 1 {%
        \hss\tenrm\ifodd\pageno 6. EXTENDED MODULI SPACE\hss\folio
        \else\folio\hfill MAPS WITH TWO CRITICAL POINTS\hfill\fi}%
    \else \hss\vbox{}\hss
    \fi}

\cl{\bf\S6. The Extended Moduli Space \[\M= \m\cup L_\infty\]}\ms

The next two sections will study the limiting behavior as the conjugacy class
\[(f)\] becomes degenerate. We will say that a sequence of conjugacy
classes \[(f_j)\] {\bit diverges to infinity\/} in the moduli space \[\m\]
if this sequence eventually leaves any compact subset of \[\m\].
Using the results of \S2, it is not hard to check that the
sequence \[\{(f_j)\}\] diverges to infinity within \[\m\] if and only if

$\bullet$ the larger of the two invariants \[|X(f_j)|\] and \[|Y(f_j)|\]
tends to infinity with \[j\],

\ni or if and only if

$\bullet$ the largest, \[\max_i|\lambda_i(f_j)|\], of the fixed point
multipliers tends to infinity with \[j\].

\ni
It will be convenient to number the fixed points of each \[f_j\] so that the
corresponding multipliers \[\lambda_i(f_j)\]
satisfy \[|\lambda_1|\le|\lambda_2|\le\cdots\le |\lambda_{n+1}|\].

{\QP{\bf Lemma 6.1.} \it If \[\{(f_j)\}\] diverges to infinity within \[\m\],
then all but two of the fixed point multipliers must tend to infinity, so that
$$ |\lambda_{n+1}|~\ge~\cdots~\ge~|\lambda_4|~\ge~|\lambda_3|~\to~\infty~. $$
Furthermore, after passing to a suitable subsequence, exactly one of the
following two statements must hold:

$(1)$ either \[\lambda_2\to\infty\] also, but \[\lambda_1\to 0\], or else

$(2)$ \[\lambda_2\] remains bounded and \[\lambda_1\] remains bounded
away from zero, while the product \[\lambda_1\lambda_2\]
converges to \[+1\].\ss}

{\bf Proof.} After passing to a subsequence, we can assume that \[X\] and
\[Y\] and each of the \[\lambda_i\] tends to a well defined limit in
\[\C\cup\{\infty\}\].
First suppose that \[|Y|\to\infty\] while \[|X|\] remains bounded.
Then it follows from 2.6 that the elementary symmetric function \[\sigma_k\]
remains bounded for \[k\ne n\], but that \[|\sigma_n|\to\infty\]. Since
\[\sigma_n\] is the sum of \[n+1\] terms, of which \[\lambda_2\lambda_3\cdots
\lambda_{n+1}\] is the largest in absolute value, it certainly follows
that \[|\lambda_2\lambda_3\cdots\lambda_{n+1}|\to\infty\]. Since the product
\[\sigma_{n+1}\] remains bounded, this implies that \[\lambda_1\to 0\].
We must show that \[|\lambda_2|\to\infty\]. Otherwise, if say \[\lambda_2\,,\,
\ldots\,,\,\lambda_{\ell}\] remained bounded, with \[\ell\ge 2\], then a
straightforward argument would show that the
\[\left({n+1\atop\ell}\right)$-fold sum \[\sigma_{n-\ell+1}\],
with dominant summand \[\lambda_{\ell+1}\cdots\lambda_{n+1}\],
also tended to infinity,
yielding a contradiction. Thus \[|\lambda_2|\to\infty\], and we are in
case 6.1(1).

Now suppose that \[|X|\to\infty\]. In particular, we assume that \[X\ne 0\],
so that
$$	\lambda_1\lambda_2\cdots\lambda_{n+1}\=\sigma_{n+1}\=n^{n+1}X^{n-1}
 \,\ne\, 0~. $$
We will need the following.

{\bf Definition.} Let
$$ \hat\sigma_k=\sigma_{n+1-k}/\sigma_{n+1}$$
be the \[k$-th elementary symmetric function of the reciprocals
\[1/\lambda_1\,,\,\ldots\,,\,1/\lambda_{n+1}\].
This is well defined and finite whenever \[X\ne 0\].
If \[X\to\infty\], we see from 2.4 and 2.6 that
$$\eqalign{ \hat\sigma_k~\to~0&\quad{\rm for}\quad k\ge 3~,\quad{\rm while}\cr
	\hat\sigma_2~\to~1&~.} \eqno (19)$$
(However, the limiting
behavior of \[\hat\sigma_1\] depends not only on \[X\] but also on \[Y\].)
Suppose that exactly \[p\] of the \[\lambda_i\] tend to zero, while
exactly \[q\] of them tend to finite non-zero limits. If \[p\ge 1\], then
it is easy to check that \[\hat\sigma_p\to\infty\] and
also that \[\hat\sigma_{p+q}\to\infty\]. By (19),
it follows that \[p\le 1\], and also that \[q=0\] whenever
\[p=1\]. Thus if \[p=1\], then we are again in case 6.1(1).

Finally, suppose that \[p=0\]. Then evidently \[\hat\sigma_q\] tends to a
finite non-zero limit, while \[\sigma_k\to 0\] for \[k>q\]. Making use of (19),
it follows immediately that \[q=2\], and that we are in case 6.1(2).\QED\ss

A convenient coordinate near the line at infinity
is provided by the sum of reciprocals
$$	\hat\sigma_1\=\sum{1\over\lambda_i}\=
{\sigma_n\over \sigma_{n+1}}\={\sigma_n\over n^{n+1} X^{n-1}}~.$$
If only \[\lambda_1\] and \[\lambda_2\] remain finite, with \[\lambda_1\to
\lambda\] and \[\lambda_2\to \lambda^{-1}\], note
that this coordinate \[\hat\sigma_1\] tends to \[\lambda+\lambda^{-1}\].
We can now define the {\bit extended moduli space\/} \[\E\] 
 to be the disjoint
union \[\m\cup L_\infty\], where \[L_\infty\] is a complex line with
coordinate \[\hat\sigma_1\in\C\]. To make this union into
a complex manifold, we cover it with two coordinate patches, each biholomorphic
to \[\C^2\]. The first coordinate patch is \[\m\] itself, with coordinates \[X\]
and \[\sigma_n\]. The second is \[(\m\ssm L_0)\cup L_\infty\] with coordinates
\[\hat X=1/X\] and \[\hat\sigma_1\]. In the overlap \[\m\ssm L_0\], these two
coordinate systems are related by the biholomorphic map
$$	\hat X\={1\over X}~,\qquad \hat\sigma_1\={\sigma_n\over n^{n+1}
  X^{n-1}}~. \eqno (20) $$

{\QP{\bf Lemma 6.2.} \it The union \[\E\], with complex structure defined
in this way, is a well defined Hausdorff complex manifold.
Furthermore, the coordinate function
\[X\] extends to a locally trivial holomorphic projection
$$	X~:~\E~\longrightarrow~\widehat\C~, $$
where each fiber \[L_X\] is conformally isomorphic to \[\C\].\ss}

The proof is straightforward.\QED\ss

(Similarly, if \[\E_\R\subset\E\] denotes the closure of the real part of
moduli space, as discussed in \S5, then we obtain a real line bundle
$$	\R~\hookrightarrow~\E_\R~\buildrel X\over\longrightarrow~
 \widehat\R=\R\cup\infty~.$$
Topologically, \[\E_\R\] is either a cylinder or a M\"obius band according
as \[n\] is odd or even.)\ss

{\bf Remark 6.3.} Although the fibers of the line bundle
\[\C\hookrightarrow\E\to\widehat\C\] have a sharp
geometric and algebraic interpretation, it is not known whether they have
any dynamic meaning. However,
there are three particular fibers which certainly do
have dynamic interpretations: The fiber \[L_0\] consists of maps
with a superattracting fixed point (compare Figure 1a). The fiber \[L_{-1}\]
consists of maps for which one critical point maps immediately to the other.
(For a study of this family, see Bam\'on and Boberieth [BB].) Finally, the
line at infinity, \[L_\infty\] is certainly quite unique.
Computer pictures suggest that each fiber \[L_X\] may intersect \[\overline
\CL\] in a set which is full and connected, so that \[L_X\cap\SL_{\rm hyp}\]
is conformally isomorphic to \[\C\ssm\overline\D\]. I know of no
reason why this should be true, although it would be compatible with the
description of the fundamental group of \[\SL_{\rm hyp}\]. (Compare 6.6
and 7.8.) I don't
have a good algorithm for making pictures of \[L_X\cap\CL\]. However, if we
pass to the \[(n+1)$-fold branched covering in which one fixed point is marked,
then we can parametrize by its multiplier \[\lambda\]. (Compare 2.10.)
It is relatively easy to make pictures in the disk \[|\lambda|<1\] where
this marked point is attracting. (See Figure 10.) Note that each \[(f)\in L_X\]
with a unique attracting fixed point embeds uniquely in this disk,
provided that \[X\ne 0\]. In particular, \[L_X\cap\SL_{\rm hyp}\] embeds
uniquely. However conjugacy classes with two attracting fixed points are
represented twice.\ss

\def\wideQP{\smallskip\leftskip=.25in\rightskip=.25in\noindent}
\midinsert
\cl{\psfig{figure=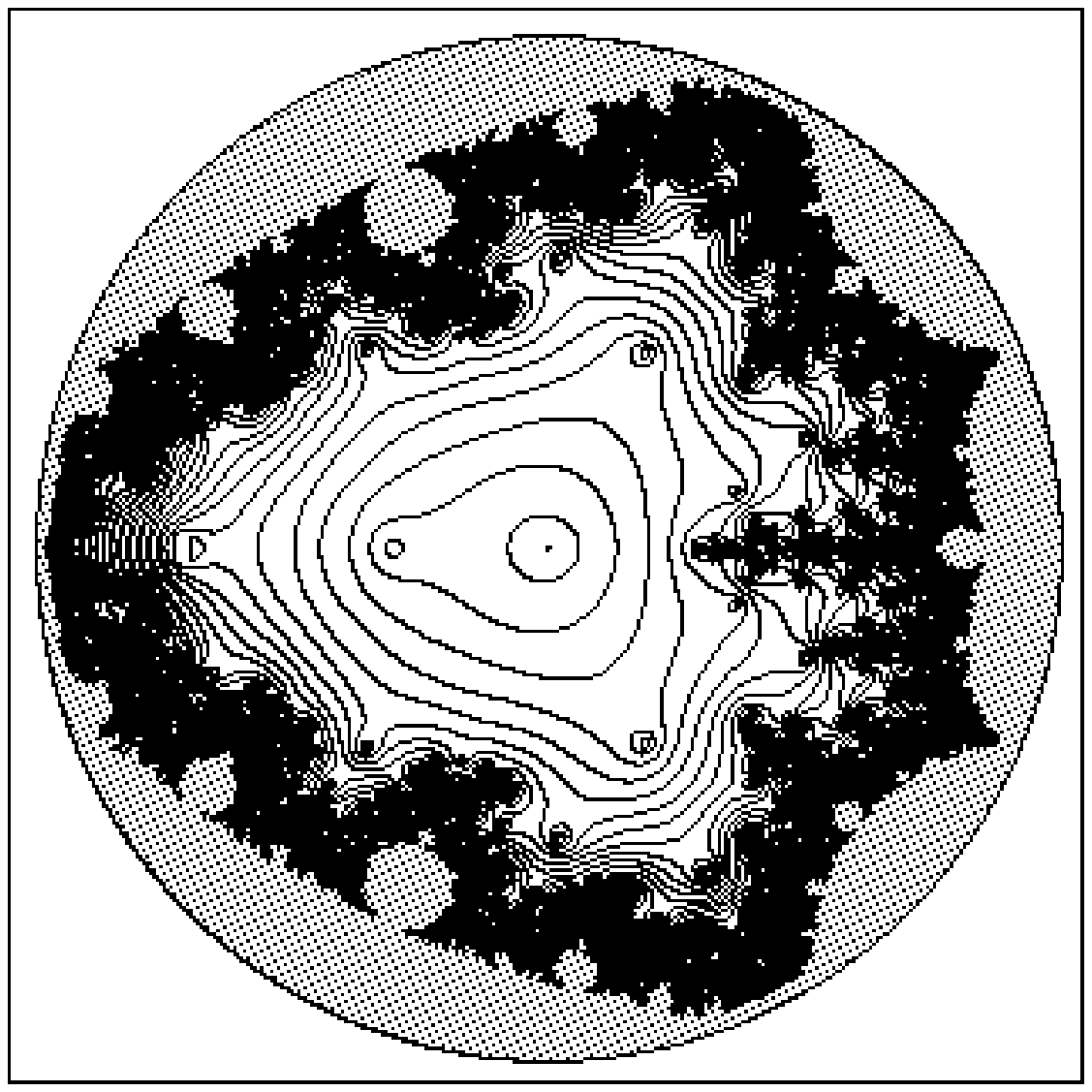,height=2in}\quad
    \psfig{figure=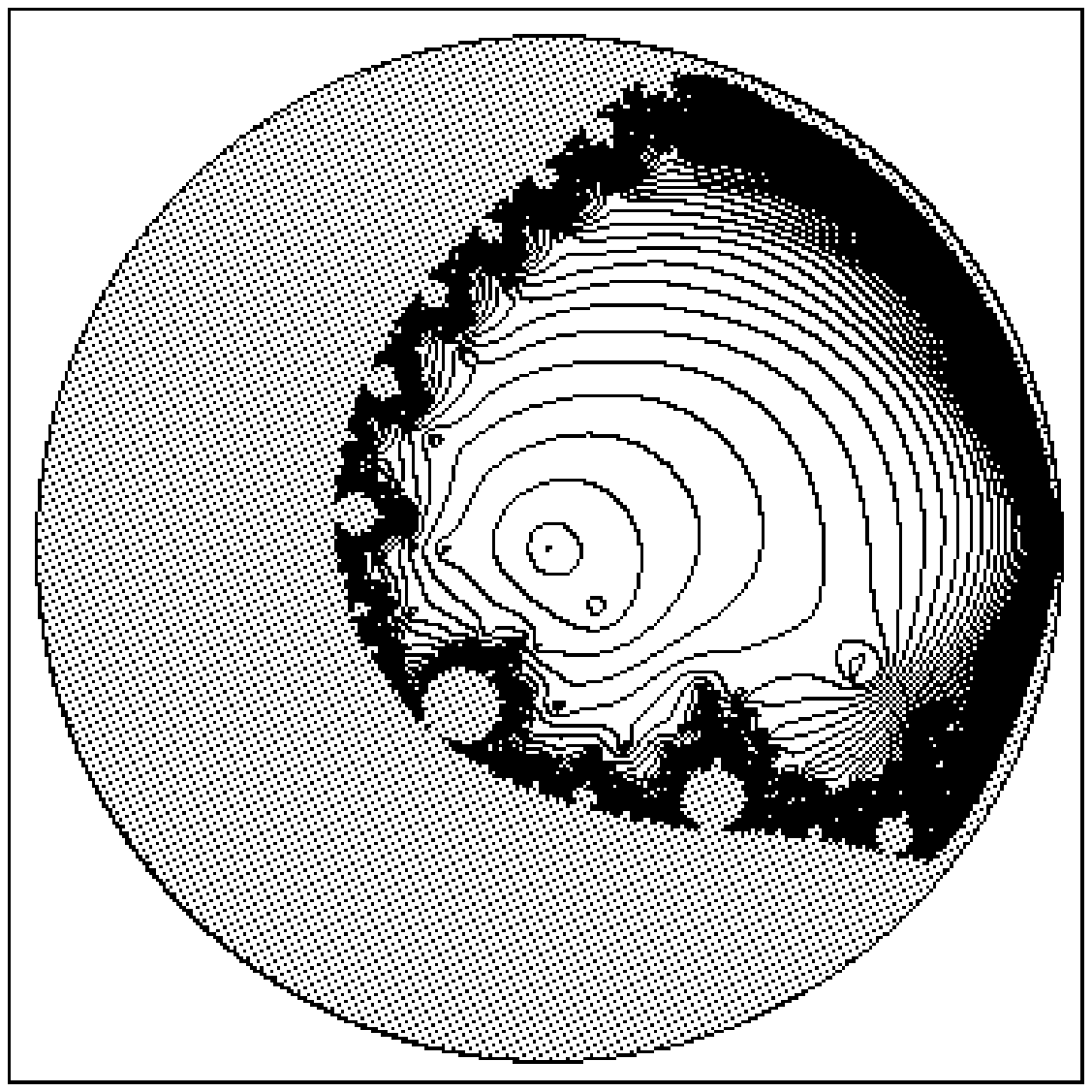,height=2in}\quad
    \psfig{figure=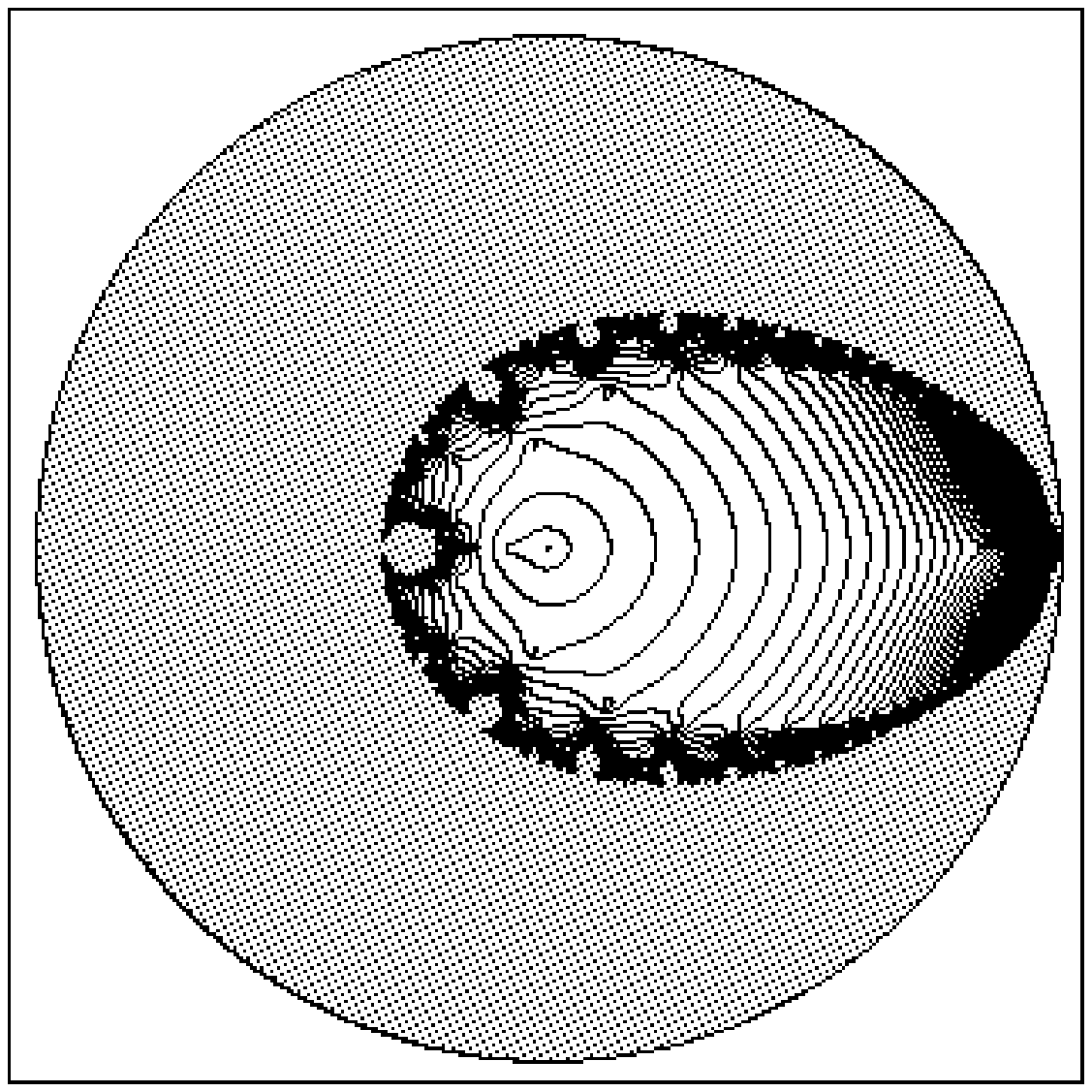,height=2in}}
\ss
{\wideQP\bit Figure 10. Pictures in the disk \[|\lambda|<1\] for \[X\] equal to
\[-.2\,,\;.2i\] and \[.2\] respectively, for degree \[n=2\]. (Compare 6.3.)
The connectedness locus is shaded grey; while
the curves in the shift locus indicate the number of iterations needed
for both critical orbits to reach a small neighborhood of the attracting fixed
point.\ss}
\endinsert

{\QP{\bf Lemma 6.4.} \it The space of all holomorphic sections of the
bundle \[X:\E\to\widehat C\] is an \[n$-dimensional vector space, consisting of
all polynomial functions of the form
$$	\sigma_n\=\sum_{j=0}^{n-1}~c_jX^j~,\qquad
	\hat\sigma_1\=\sum_{j=0}^{n-1}~ c_j\,\hat X^{n-1-j}/n^{n+1}~.$$
It follows that a function \[Y=Y(X)\] gives rise to a holomorphic section
of this bundle if and only if the sum \[Y(X)+2X^n\] is a polynomial of degree
\[\le n-1\] in \[X\].\ss}

The proof is straightforward. (Compare 2.4 and 2.6.)\QED

{\QP{\bf Corollary 6.5.} \it For each \[\lambda\ne 0\], the closure within
\[\E\] of the affine algebraic curve \[\per_1(\lambda)\subset\m\] is a smooth
compact algebraic curve \[\overline\per_1(\lambda)\subset\E\], which can be
described as the image of a smooth section of the complex line bundle\break
\[\M\to\widehat\C\]. The intersection \[\overline\per_1(\lambda)\cap L_\infty\]
of this curve with the line at infinity consists of the point
with coordinates
\[~~	\hat X=0~,\quad \hat\sigma_1\=\lambda+\lambda^{-1}~. \]\ss}

{\bf Proof.} This follows immediately, using 2.3 and 2.4.\QED

Note that the coordinate \[\hat\sigma_1\] on the line
at infinity is real, and belongs to the interval
\[-2\le\hat\sigma_1\le 2\], if and only if \[|\lambda|=1\]. (Compare 7.5.)\ss

(As noted in \S2, the curve \[\per_1(0)\] should be identified with
the fiber \[L_0\] of this fibration counted \[n-1\] times. Similarly,
it may be useful to identify \[\per_1(\infty)\]
with the fiber \[\L_\infty\] counted \[n-1\] times.)

Other examples of smooth sections of this  complex line
bundle are the coordinate
curve \[\sigma_n=0\], and the half symmetry locus \[Y=-2X^{(n-1)/2}
(X+1)^{(n+1)/2}\] of \S1 when \[n\] is odd (but not the locus \[Y=0\]
or the other half symmetry locus).\ss

\midinsert
\cl{\psfig{figure=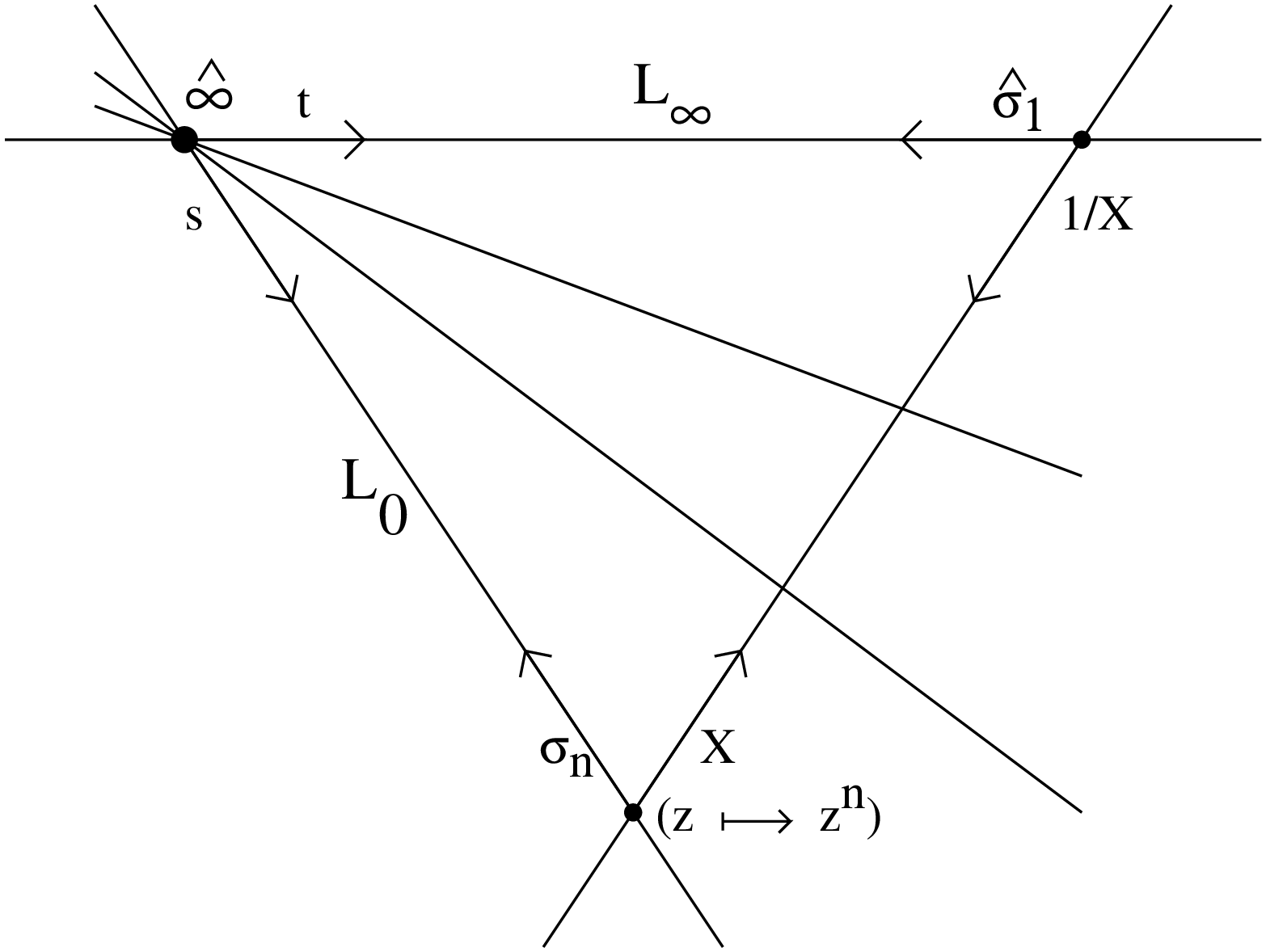,height=2.6in}}

{\QP\bit Figure 11. Schematic picture of the compactified moduli space
\[\E\cup\widehat\infty\], with a singular point at \[\widehat\infty\], showing
the three overlapping systems of coordinates.\par}
\endinsert

{\bf Remark 6.6.}
We can actually compactify moduli space by adding one more point
which we will call \[\widehat
\infty\], contained in a third coordinate neighborhood. However
the result is no longer a manifold when \[n>2\], but rather an orbifold.
In fact a neighborhood of \[\widehat\infty\] is homeomorphic to a cone
over a 3-dimensional
lens space whose fundamental group is cyclic of order \[n-1\].
This third coordinate neighborhood \[W\] can be parametrized by coordinates
\[s,t\] subject to the identifications \[(s,t)=(\alpha s, \alpha t)\],
where \[\alpha\]
ranges over all \[(n-1)$-st roots of unity. They are related to the coordinates
\[(X,\sigma_n)\] on the overlap where \[\sigma_n\ne 0\] and \[s\ne 0\],
and to the coordinates \[(\hat X,\hat\sigma_1)\] on the
overlap where \[\hat\sigma_1\ne 0\] and \[t\ne 0\], by the identities
$$	X\={t\over s}~,\quad \sigma_n\={1\over s^{n-1}}~,\quad{\rm and}\quad
	\hat X\={s\over t}~,\quad \hat\sigma_1\={1\over n^{n+1}t^{n-1}}~. $$
Note that
there is just one point \[\widehat\infty\in W\] which does not belong to \[\M\],
namely the singular point \[\widehat\infty\] which has coordinates \[s=t=0\].
It is not difficult to check that these coordinate transformations
are compatible,
so that \[\M\cup\widehat\infty\] is a well defined compact Hausdorff space.
A small
neighborhood of infinity in \[\M\] can be identified with a neighborhood
of \[\widehat\infty\] in \[W\] with this singular point removed. Evidently,
for reasonable choice of neighborhood, the fundamental group will be cyclic
of order \[n-1\].

Intuitively, we can think of the fibers of our holomorphic line bundle
as a pencil of lines though this exceptional point \[\widehat\infty\]. This
pencil of lines sweeps out the compactified moduli space \[\M\cup\widehat
\infty\]. Using this construction, we can reinterpret Lemma 6.1 as follows.
If a sequence of points \[(f_j)\in\m\] converges
within \[\E\cup\widehat\infty\] to the point \[\widehat\infty\], then
we are in
Case (1) of 6.1, with all but one of the fixed point
multipliers tending to infinity.
On the other hand, if a sequence converges to some point of \[L_\infty
\subset\E\], then we are in Case (2) of 6.1, so that two of the
multipliers converge to finite non-zero limits with product \[+1\] while the
others tend to infinity.
\bs\bs

\headline={%
    \ifnum\pageno > 1 {%
        \hss\tenrm\ifodd\pageno 7. EXTENDED SHIFT LOCUS\hss\folio
        \else\folio\hfill MAPS WITH TWO CRITICAL POINTS\hfill\fi}%
    \else \hss\vbox{}\hss
    \fi}

\cl{\bf\S7. The Extended Hyperbolic Shift Locus \[\widehat\SL_{\rm hyp}\subset
\M\]}\ms

We first prove the following.

{\QP{\bf Theorem 7.1.} \it  The connectedness locus \[\CL\] is contained in
a compact subset of the extended moduli space \[\M\]. Hence the closure
\[\overline\CL\subset\E\] is compact.\ss}

It follows that the complement \[\E\ssm\overline\CL\]
is a neighborhood of infinity in \[\E\].
By definition, this complementary open neighborhood of infinity,
consisting of \[{\cal S}_{\rm hyp}\]
together with the points of \[L_\infty\ssm(\overline\CL\cap L_\infty)\],
will be called the {\bit extended hyperbolic shift locus\/}
\[\widehat{\cal S}_{\rm hyp}\]. (See 7.4 for a precise description
of \[\overline\CL\cap L_\infty\].)
The following statement is completely equivalent to 7.1.

{\QP{\bf Corollary 7.2.} \it There exists a constant \[k\] depending only on
the degree \[n\] with the following property. If all but one of the fixed
point multipliers of \[f\] are greater than \[k\] in absolute value, then
\[(f)\] belongs to the hyperbolic shift locus.\ss}

(Compare [M2, \S8.8], which shows that the best value for \[k\] in the
degree 2 case lies between 3 and 6.)\ss

{\bf Proof of 7.1 and 7.2.} We will first find a rough criterion for
belonging to the hyperbolic shift locus depending on
the invariants \[X\,,\,Y_1\,,\,Y_2\,,\,Y\] of \S1, and then relate this
to the topology of \[\E\]. 
Consider a sequence of maps \[f_j\] with marked critical points so that the
conjugacy classes  \[(f_j)\] tend to infinity in \[\m\]. Thus, passing to
a subsequence if necessary,
at least one of the invariants \[X(f_j)\] and \[Y(f_j)\]
must tend to infinity. Since \[Y_1+Y_2=Y\] and
\[Y_1Y_2=(X+1)^{n+1}X^{n-1}\] by (3), it follows that \[Y_1\]
or \[Y_2\] must tend to infinity. Interchanging the two critical points if
necessary, we may assume that \[Y_2(f_j)\to\infty\] as \[j\to\infty\].

In particular, putting the \[f_j\] into the normal form (1), we may assume
that the coefficients \[c\] and \[d\] are non-zero and hence, after a linear
change of coordinate, we may assume that \[c=d=1\]. In this way, taking
the critical values to be \[f(0)=v\] and \[f(\infty)=v+\delta\], we obtain
the normal form
$$	f(z)\=f_j(z)\={(v+\delta)z^n+v\over z^n+1}\=v+{\delta z^n\over z^n+1}~.
	\eqno (21) $$
(We will write \[v=v_j\] and \[\delta=\delta_j\] whenever it is necessary for
clarity.) Note that \[f^{-1}(\infty)\] is the set of \[n$-th roots of \[-1\].
The invariants
$$	X\=v/\delta\,,\quad Y_2\=1/\delta^n $$
can be computed from $(2)$. It follows that \[\delta\to 0\], and that
the critical value \[v\] satisfies the equation \[v^n=X^n/Y_2\]. It will
be convenient to set \[\epsilon=\sqrt{|\delta|}>0\]. Thus
\[\epsilon=\epsilon_j\] also tends to zero as \[j\to\infty\].\ss

\midinsert
\cl{\psfig{figure=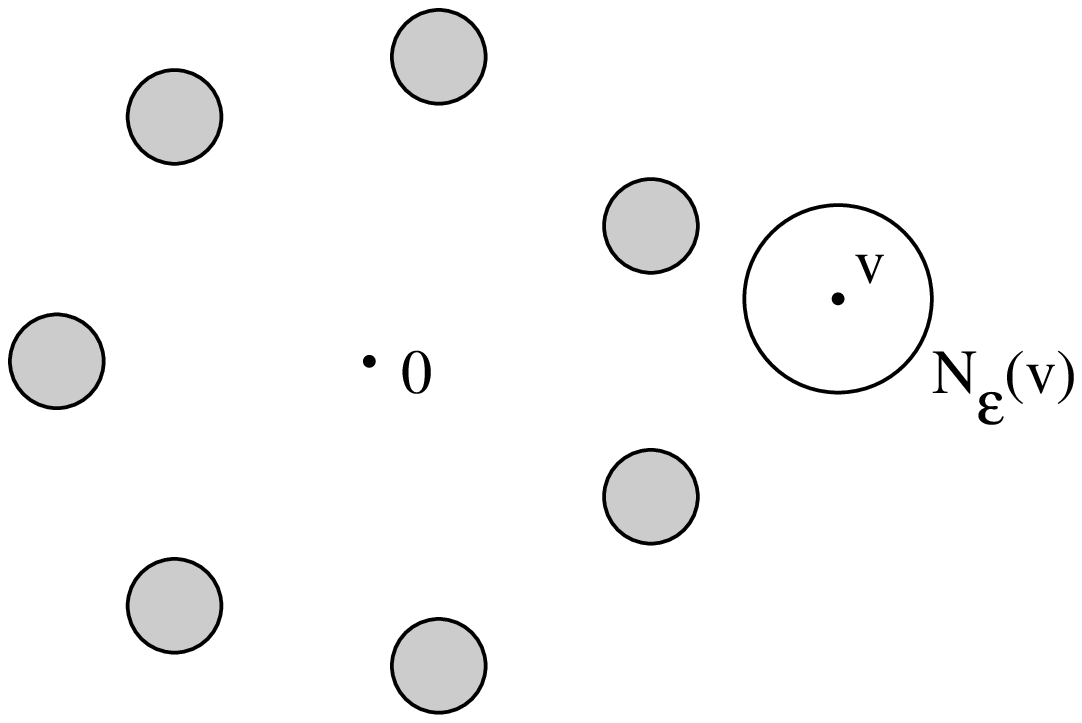,height=1.5in}}\ss
{\QP\bit Figure 12. Dynamics for \[(f)\] near \[\widehat\infty\]. Here the
critical points are at \[0\] and \[\infty\], while the critical values
as well as the attracting fixed point are in the disk \[N_\epsilon(v)\].
The Julia set is contained in the \[n\] shaded disks \[\overline U_i\],
each of which maps
diffeomorphically onto the complement \[\widehat\C\ssm\overline
 N_\epsilon(v)\]. As \[(f)\to\widehat\infty\], each of these \[n+1\] disks
shrinks to a point.\ss}
\endinsert

{\bf Case 1.} Suppose that \[v^n=X^n/Y_2\] is bounded away from \[-1\], or
equivalently that \[v\] is bounded away from the \[n$-th roots of \[-1\],
as \[j\to\infty\]. Then we will prove that \[(f_j)\] belongs to the
hyperbolic shift locus for large \[j\]. Let
\[N_\epsilon(v)\] be the open
neighborhood of radius \[\epsilon=\sqrt{|\delta|}\] about \[v\].
Then the pre-image \[W=f^{-1}N_\epsilon(v)\] consists of all \[z\]
with
$$	|z^n+1|~>~ |\delta|/\epsilon\=\epsilon~. $$
If \[\epsilon\] is small,
this set can be described as the complement \[W=
\widehat\C\ssm(\overline U_1\cup\cdots\cup \overline U_n)\]
of the union of closed neighborhoods \[\overline U_1\,,\,\ldots
\,,\,\overline U_n\] of the various \[n$-th roots of \[-1\], where each
\[\overline U_j\] has radius roughly \[\epsilon/n\].
Thus \[W\] is a connected open set which contains both
critical points. (Compare the figure, where \[W\] is the complement
of the shaded region.) Since \[v\] is bounded away from the roots of \[-1\],
it follows that the image \[f(W)=N_\epsilon(v)\] is compactly
contained in \[W\] provided that \[j\] is sufficiently large. Using the
Poincar\'e metric for \[W\], it follows that all orbits
in \[W\] converge to a common attracting fixed point. Using Theorem B.5 of
Appendix B, it follows
that \[(f_j)\] is contained in the hyperbolic shift locus.\ss

{\bf Case 2.} Suppose in fact that \[v^n\] converges to \[-1\]
(and hence that \[|v|\to 1\]), but that
$$	|1+v^n|/\epsilon~\to~\infty\qquad{\rm as}\qquad j\to\infty~,$$
so that \[\epsilon\ll |1+v^n|\] for large \[j\]. Note that the distance
between \[v\] and the closest \[n$-th root of \[-1\] is approximately
\[|1+v^n|/n\] for large \[j\]. It now follows, just as in Case 1, that
\[f(W)\] is compactly contained in \[W\], and hence that \[(f_j)\] belongs
to the hyperbolic shift locus, provided that \[j\] is large.\ss


{\bf Remark 7.3.} More explicitly, in both Cases 1 and 2, it follows
for large \[j\] that the Julia set \[J\] is the disjoint
union of compact subsets \[J_k=J\cap \overline U_k\], where \[f\] maps each \[\overline U_k\]
diffeomorphically onto the strictly
larger set \[\widehat\C\ssm N_\epsilon(v)\].
Hence \[f\] is hyperbolic, that every point of \[J\] is
uniquely determined by its itinerary with respect to the partition
\[J=J_1\cup\ldots\cup J_n\]. In fact, as \[j\to\infty\], it follows that
the Julia set becomes more and more hyperbolic, in the sense that
the multiplier of any periodic orbit in \[J\] tends to infinity.

{\bf Case 3.} We now return to the proof of 7.1 and 7.2.
Suppose finally that the ratio \[|1+v^n|/\epsilon\] remains
bounded as \[j\to\infty\], although \[\epsilon\to 0\] hence
\[1+v^n\to 0\].
(By Cases 1 and 2, this must be the case if each \[(f_j)\] belongs to the
connectedness locus.) Then we will show that the \[(f_j)\] remain within
some compact subset of \[\E\]. Setting
$$	1+v^n\=1+X^n/Y_2\=\eta~,$$
we see that \[|\eta|\] is less than some constant time \[\epsilon\], so that
\[|\eta^2|\] is less than a constant times \[\epsilon^2=|\delta|\sim 1/|X|\].
We can now solve for
$$\eqalign{	Y_2\=&-X^n(1+\eta+\eta^2+\cdots)\=-X^n-X^n\eta +O(X^{n-1})~,\cr
	Y_1\=& (X+1)^{n+1}X^{n-1}/Y_2\= -(X+1)^{n+1}(1-\eta)/X
	\=-X^n+X^n\eta+O(X^{n-1})~,\cr} $$
and hence
$$	Y\=Y_1+Y_2\=-2X^n+O(X^{n-1})~.$$
Since \[X\] tends to infinity, it is appropriate to pass to the
coordinates \[\hat X=1/X\] and \[\hat\sigma_1=\sum 1/\lambda_k\] of \S6.
Recall from 2.6 that
$$	\hat\sigma_1\={\sigma_n\over\sigma_{n+1}}\= {n^n(Y+P_n(X))\over
 n^{n+1}P_{n+1}(X)}\={Y+2X^n+O(X^{n-1})\over nX^{n-1}}	~. $$
It follows that \[\hat\sigma_1\] remains bounded as \[j\to\infty\] and
\[X\to\infty\], so that the \[(f_j)\] remain within a compact subset of the
extended moduli space \[\E\].\ss

Theorem 7.1 now follows immediately, and 7.2 is proved as follows.
Recall that \[\lambda_2\] is defined to be the fixed point multiplier
with next-to-smallest absolute value. Suppose, for every positive integer
\[j\], that we could find a conjugacy class \[(f_j)\in\m\cap\overline\CL\] with
\[|\lambda_2(f_j)|>j\]. The resulting sequence
cannot be in Case 1 or Case 2, and therefore
must be in Case 3, with \[X(f_j)\to\infty\] and \[\hat\sigma_1(f_j)\] bounded.
But this contradicts the hypothesis that \[|\lambda_2(f_j)|\to\infty\],
and completes the proof of 7.2.\QED\ss

We can describe the set \[\overline\CL\subset\E\] and its complement
\[\widehat{\cal S}_{\rm hyp}\] more precisely as follows.

{\QP{\bf Theorem 7.4.} \it The closure of the connected locus \[\CL\] within
the extended moduli space \[\E=\m\cup L_\infty\] is a compact set which
consists of \[\CL\cup\per_1(1)\subset\m\], together with a closed line segment
consisting of all points in the line \[L_\infty\cong\C\] such that the
coordinate \[\hat\sigma_1\] is real with \[-2\le\hat\sigma_1\le 2\].\ss}

Equivalently, we can say that \[\overline\CL\cap L_\infty\] consists of
all points of \[L_\infty\] with \[\hat\sigma_1\] of the form
$$	\hat\sigma_1\=2\cos\theta\=\lambda+\lambda^{-1}\qquad{\rm where}
	\qquad \lambda=e^{i\theta}~, $$
so that \[\lambda\] ranges over the unit circle. It follows that the complement
$$	\widehat{\cal S}_{\rm hyp}\=\E\ssm\overline\CL $$
consists of \[{\cal S}_{\rm hyp}~\subset\m\], together with all points
of \[L_\infty\] with \[\hat\sigma_1\in\C\ssm[-2,2]\] or in other words
all points for which
$$	\hat\sigma_1\=\lambda+\lambda^{-1}\qquad{\rm with}\qquad
	0~<~|\lambda|~<~1~. $$

The proof of 7.4 begins as follows. To see that \[\overline\CL\] contains
all points of \[L_\infty\] with \[\hat\sigma_1=e^{i\theta}+e^{-i\theta}\],
note that each such point can be described as the intersection
$$	L_\infty\cap\overline\per_1(e^{i\theta}) ~. $$
But clearly \[\per_1(e^{i\theta})\] is disjoint from the hyperbolic
shift locus, and hence contained in \[\overline\CL\].

The proof in the other direction
will depend on a study of dynamic behavior for conjugacy
classes \[(f)\in\m\subset\M\] which are sufficiently close to the line
\[L_\infty\subset\E\].
We must first find a normal form for such classes.
Suppose the multipliers at the \[n+1\] fixed points are approximately
\[\lambda\,,\,\lambda^{-1}\,,\,\infty\,,\,\ldots\,,\,\infty\]. Here we must
exclude the case \[\lambda=0\], and also the case
\[\lambda=1\] which turns out to be particularly difficult. We will use
the normal form (1) with \[ad-bc=1\] and with \[a+b=c+d=u\], so that there is
a fixed point of multiplier \[\lambda=n/u^2\] at \[z=1\]. Suppose that there is
another fixed point at \[z=\kappa\]. (Here we must assume that \[\kappa^n\ne
 1\], since \[\kappa^n=1\] would imply that \[f(\kappa)=f(1)=1\]. In
practice, we will only be interested in the limiting case as \[\kappa\to 1\].)
Solving the fixed point equation \[a\kappa^n+b=(c\kappa^n+d)\kappa\] (together
with the equations \[ad-bc=1\] and \[a+b=c+d=u\]), we find the unique
solution  
$$\eqalign{	a\=-P+\kappa Q~,\quad& b\=\kappa^nP-\kappa Q~,\cr
	c\=~-P+Q~,~\quad& d\=\kappa^nP-Q~,} $$
where we have temporarily introduced the abbreviations
$$	P\={u\over \kappa^n-1}~,\qquad Q\={1\over u(\kappa-1)}~.$$
It follows easily that
$$	X\=bc\= -(\kappa^nP-\kappa Q)(P-Q) $$
is given by the asymptotic formula
$$	X~\sim~ -{(\lambda-1)^2\over n\lambda}\,{1\over (\kappa-1)^2}\qquad
	{\rm as}\qquad \kappa\to 1~. $$
{\it Thus, given \[\lambda\ne 0\,,\,1\], we can realize any large value of the
invariant \[X\] by suitable choice of the parameter \[\kappa\approx 1\].\/}
(Note by 2.1 that \[X\]
and \[\lambda\] together determine the remaining invariant \[Y\].)
The multiplier at \[\kappa\], call it \[\lambda'\], can be computed by the
formula
$$	\lambda'\={n\kappa^{n-1}\over (c\kappa^n+d)^2}\=n\kappa^{n-1}
  \left({u(\kappa-1)\over \kappa^n-1}\right)^2~\sim~u^2/n\=1/\lambda $$
as \[\kappa\to 1\]. For \[\kappa\approx 1\],
the two critical values \[f(\infty)=a/c\] and \[f(0)=b/d\]
are close to \[1\] and \[\kappa\], and extremely close to each other.
In fact the difference is given by
$$	{a\over c}-{b\over d}\= {1\over cd}~\sim~{1\over X} \=O(\kappa-1)^2~,$$
and we can compute
$$	{a\over c}\=t\cdot 1+(1-t)\cdot \kappa\qquad{\rm with}\qquad
	t\={1\over 1-Q/P}~\sim~{1\over 1-\lambda}~, $$
and similarly \[b/d=t'\cdot 1+(1-t')\cdot \kappa\] with
\[t'\sim 1/(1-\lambda)\].\ms

It will be convenient to conjugate \[f\] by the M\"obius involution
$$	\phi(z)\=\phi^{-1}(z)\=\kappa\,{z-1\over z-\kappa}~. $$
which interchanges the critical point \[0\] with the fixed point \[1\],
and interchanges the critical point \[\infty\] with the fixed point \[\kappa\].
Thus the conjugate map
$$	g\=\phi\circ f\circ\phi $$
has critical points at \[\phi(0)=1\] and at \[\phi(\infty)=\kappa\], and has a fixed
point of multiplier \[\lambda\] at \[\phi(1)=0\] and a fixed point
of multiplier \[\approx\lambda^{-1}\] at \[\phi(\kappa)=\infty\]. The critical
values \[g(1)\] and \[g(\kappa)\]
are both close to \[\lambda\]. In fact computation shows that
$$	g(\kappa)\=\kappa Q/P\={\kappa+\kappa^2+\cdots+\kappa^n\over u^2}\=
	\lambda+O(\kappa-1)\qquad{\rm as}\qquad \kappa\to 1~, $$
and similarly \[g(1)=g(\kappa)/\kappa^n=\lambda+O(\kappa-1)\].

{\QP{\bf Theorem 7.5.} \it Suppose that we fix \[\lambda\ne 0\,,\,1\].
Then as \[\kappa\to 1\] the rational function \[g=g_\kappa\] converges uniformly
to the linear function \[w\mapsto \lambda\,w\] throughout any compact subset
of \[\widehat\C\ssm\{1\}\].\ss}

(Compare [M2, \S4].
Here `uniform convergence' refers to convergence with respect to
the spherical metric on the target space \[\widehat\C\].) Thus, as \[\kappa\to
1\] the two critical points of \[g\] crash together, and \[g\] converges to a
linear map, except in a small neighborhood of \[1\].

The proof will be based on the following. Let \[\D_r\] be the disk
\[\{z\in\C~;~|z|<r\}\]. If \[U\] is an open neighborhood
of \[0\] in \[\C\], it will be convenient to define the {\bit inradius\/}
$$	s(U_i)\=\max\{s~;~\D_s\subset U\}\={\rm dist}(0\,,\,\C\ssm U) $$
as the maximum radius of a disk centered at the origin and contained in \[U\].

{\QP{\bf Lemma 7.6.} \it Given simply connected open sets \[U_i\subset\C\] with
inradius \[s(U_i)\] tending to infinity, there exist radii \[r_i\to\infty\]
and conformal isomorphisms\break
\[h_i:\D_{r_i}{\buildrel\cong\over\longrightarrow}U_i\] so that \[\{h_i\}\]
converges uniformly to the identity map on any compact subset of \[\C\].\ss}

{\bf Proof.} Let \[u\] range over univalent maps \[u:\D\to\C\]
on the unit disk with \[u(0)=0\]
and \[u'(0)=1\]. Note that the space of all such \[u\] is compact, and that
the inradius of the image always satisfies
$$ s\Big(u(\D)\Big)~\le~ 1~. \eqno(22) $$
(See for example [CG, \S I.1].) It follows that the second derivative of \[u\]
is uniformly bounded on any compact subset of \[\D\], say
$$	|u''(w)|~\le~ 2k\qquad{\rm for}\qquad |w|~\le~ 1/2~. $$
Hence
$$	|u(w)-w|~\le~k|w^2|\qquad{\rm for}\qquad  |w|~\le~ 1/2 \eqno (23) $$
and for all such maps \[u\], where \[k\] is a uniform constant.
Let \[\rho_i:\D\buildrel\cong\over\longrightarrow U_i\]
be the Riemann map, satisfying
\[\rho_i(0)=0\] and \[\rho_i'(0)=r_i>0\], and define \[h_i:\D_{r_i}\to U_i\]
by \[h_i(z)=\rho_i(z/r_i)\]. Applying (22) and (23) to the map \[u(w)=
\rho_i(w)/r_i =h_i(r_iw)/r_i\], we find easily that
\[s(U_i)\le r_i\] so that \[r_i\to\infty\], and that
$$	|h_i(z)-z|~\le~ k|z^2|/r_i\qquad{\rm for}\qquad |z|~\le ~r_i/2~. $$
Thus \[\{h_i\}\] converges uniformly to the identity for \[z\] restricted to
any compact subset of the plane.\QED\ss

{\bf Proof of 7.5.}
Let \[\Gamma=\Gamma_\kappa\] be a circle arc joining the two critical
values \[g(1)\] and \[g(\kappa)\]. Then the preimage \[g^{-1}(\Gamma)\] is
a union of \[n\] circle arcs joining the two critical points
\[1\] and \[\kappa\], where two
adjacent arcs span an angle of \[2\pi/n\] at either common endpoint.
Furthermore, each of the \[n\]
connected components of \[\widehat\C\ssm g^{-1}(\Gamma)\] maps
biholomorphically onto \[\widehat\C\ssm \Gamma\]. In fact
the corresponding statement is true for any bicritical map, since any such map
can be expressed as a composition of the \[n$-th power map for which it is
clearly true, together with M\"obius transformations which carry circles to
circles. (Compare the normal form (1).)

Both the diameter of \[\Gamma\] and the diameter of \[g^{-1}(\Gamma)\] depend on the
precise choice of circle arc \[\Gamma\]. However, with a little care we can choose
\[\Gamma=\Gamma_\kappa\] so that both of these diameters
tend to zero as the distance \[|\kappa-1|\] between the two critical points
tends to zero. It then follows that there is one
largest component \[U=U_\kappa\] of \[\widehat\C\ssm g^{-1}(\Gamma)\], and that
the diameter of the complement \[\widehat\C\ssm U\] tends to zero
as \[\kappa\to 1\]. After conjugating by a rotation
of the Riemann sphere which interchanges \[1\] and \[\infty\], we can
use the Lemma to construct conformal isomorphisms
\[h_\kappa:\Delta_\kappa\to U_\kappa\], where \[\Delta_\kappa\subset
\widehat\C\ssm\{1\}\] is the complement  of a small round disk about \[1\],
so that \[h_\kappa\] converges
uniformly to the identity map on any compact subset of \[\widehat\C\ssm
\{1\}\]. Further, without loss of generality, we may assume that \[h_\kappa\]
fixes the points \[0\] and \[\infty\].

Similarly we can choose conformal isomorphisms \[h'_\kappa:\widehat
\C\ssm \Gamma_\kappa\to\Delta'_\kappa\] fixing \[0\] and \[\infty\],
where \[\Delta'\] is the complement of a small round disk about \[g_\kappa(1)
\approx \lambda\], so that \[h'_\kappa\] converges
uniformly to the identity on any compact subset of \[\widehat\C\ssm
\{1/\lambda\}\]. Now the composition \[h'_\kappa\circ g_\kappa\circ
h_\kappa\] maps the round disk \[\Delta_\kappa\]
conformally onto the round disk \[\Delta'_\kappa\], and hence
extends to a M\"obius automorphism of the Riemann sphere. Since it fixes
\[0\] and \[\infty\], it must have the form
$$	h'_\kappa\circ g_\kappa\circ
h_\kappa(w)\=\lambda_\kappa w~, $$
where \[\lambda_\kappa\to\lambda\] as \[\kappa\to 1\]. Now since \[h_\kappa\]
and \[h'_\kappa\] converge uniformly to the identity except near \[1\] and
\[g_\kappa(1)\] respectively, it follows that \[g_\kappa(w)\] converges to
\[\lambda w\] except near \[1\].\QED\ss

{\QP{\bf Corollary 7.7.} \it For \[0<|\lambda|<1\] and for \[|X|\] greater
than some constant depending continuously on \[\lambda\],
the conjugacy class \[(f)\] belongs to the hyperbolic shift locus.\ss}

{\bf Proof.} Suppose that \[|\lambda|<1-\epsilon<1\]. Then \[|g(w)|\approx
|\lambda w|\le |w|(1-\epsilon)\] for \[\kappa\]
close to \[1\], or equivalently for \[|X|\]
large, provided that \[w\] is bounded away from \[1\]. Thus every \[w\]
in the disk \[D_{1-\epsilon}\] 
belongs to the attracting basin of \[0\]
for \[|X|\] large. Since the two critical
values of \[g\] are close to \[\lambda\], with \[|\lambda|<1-\epsilon\],
they both belong to this basin. By B.5 in the Appendix, it follows
that the conjugacy class \[(f)=(g)\] belongs to \[{\cal S}_{\rm hyp}\].\QED\ss

Evidently Theorem 7.4 follows immediately.\QED\ss

{\QP{\bf Lemma 7.8.} \it 
This open set \[{\cal S}_{\rm hyp}\subset\M\] is connected. Its fundamental
group is cyclic of order \[n-1\].\ss}

The proof can be outlined as follows. Note first that there is a
holomorphic map
$$	\pi\,:\,{\cal S}_{\rm hyp}~\longrightarrow~\D $$
to the open unit disk, which maps each \[(f)\] to the multiplier at its unique
attracting fixed point. (The extension of \[\pi\] to the line at infinity is
straightforward: If \[\hat\sigma_1=\lambda+\lambda^{-1}\] with \[|\lambda|<1\],
then \[\pi\] takes the value \[\lambda\].) For any \[\lambda\in\D\ssm\{0\}\],
the fiber
\[\pi^{-1}(\lambda)\] is isomorphic to the open disk, with a point on
\[L_\infty\] as center. (Compare the discussion in [GK] and [M2].)
 In this way we obtain a topological
fibration
$$	\D~~\to~~\pi^{-1}(\D\ssm\{0\})~~\to~~\D\ssm\{0\}~. $$
However, the fiber \[\pi^{-1}(0)\] is different in two ways. First it is
missing a point at infinity, and hence is a punctured disk. Second, there is
no local cross-section near zero. Instead we have something like a Seifert
fibration. In fact, since \[\sigma_{n+1}/n^{n+1}=X^{n-1}\] by Corollary 2,
it follows that
$$	\lambda~\cong~ \hat\sigma_1^{-1}\=n^{n+1}\,X^{n-1}/\sigma_n $$
for \[\lambda\] close to zero and \[A\] close to \[+1\] in the shift locus.
Thus a small loop around the line \[A=1\] in the shift locus corresponds
to a loop in the \[\lambda$-plane which winds \[n-1\] times around the origin.
Thus \[n-1\] times a generator for the fundamental group of
\[\pi^{-1}(\D\ssm\{0\}\], which is free cyclic, maps to zero in the
fundamental group of \[{\cal S}_{\rm hyp}\].
\bs\bs

\headline={%
    \ifnum\pageno > 1 {%
        \hss\tenrm\ifodd\pageno 8. THE CURVES \[\per_p(\lambda)\]\hss\folio
        \else\folio\hfill MAPS WITH TWO CRITICAL POINTS\hfill\fi}%
    \else \hss\vbox{}\hss
    \fi}

\cl{\bf \S8. The Curves \[\overline\per_p(\lambda)\subset \E\]}\ms

Roughly speaking, the curve \[\per_p(\lambda)\] is the set of all
conjugacy classes with an orbit of period \[p\ge 1\] and multiplier \[\lambda
\in\C\].
However, we must be somewhat careful to give a definition which yields a
well behaved algebraic curve which depends continuously on \[\lambda\],
even in exceptional cases. (One
difficulty arises when the multiplier is a root of
unity, since a sequence of orbits of
period \[pq\] with multiplier converging to \[+1\] may converge to an orbit
of period \[p\] with a \[q$-th root of unity as multiplier. Another difficulty
arises for \[\lambda=0\], since the curves \[\per_p(\lambda)\] converge towards
an \[(n-1)$-fold branched covering of \[\per_p(0)\] as \[\lambda\to 0\].)\ss

Let us first count periodic orbits. As in [M2], [M3],
define positive integers \[\nu_n(p)\] by the formula
$$	n^p\=\sum_{q|p}\nu_n(q)\qquad\Longleftrightarrow\qquad \nu_n(p)\=
	\sum_{q|p} \mu(p/q)\,n^q~, $$
to be summed over all divisors \[1\le q\le p\], where the M\"obius
function \[\mu(p/q)\] is defined to be
\[(-1)^m\] if \[p/q=\ell_1\cdots \ell_m\] is a product of \[m\] distinct primes
\[\ell_j\], and \[\mu(p/q)=0\] otherwise. 
The first few values are as follows.\ms

\halign{\quad\hfil#\hfil&\qquad\hfil#\hfil&\quad\hfil#\hfil&\quad\hfil#\hfil&\quad\hfil#\hfil&
\quad\hfil#\hfil&\quad\hfil#\hfil\cr
period \[p\] &1& 2&3&4&5&6\cr
number \[\nu_n(p)\]&n&  \[n^2-n\]& \[n^3-n\]& \[n^4-n^2\]&
\[n^5-n\] & \[n^6-n^3-n^2+n\]\cr}
\ms

\ni Since \[f^{\circ p}\] is a rational map of degree \[n^p\], with \[n^p+1\]
fixed points in the generic case, it follows easily that a generic \[f\]
has \[\nu_n(p)\] points of period \[p\] for \[p\ge 2\]. (However, for \[p=1\]
it has \[\nu_n(1)+1=n+1\] fixed points.)

We next construct a commutative diagram
$$\matrix{	\PER_p & \to & \bic\cr
	\downarrow & & \downarrow\cr
	\per_p & \to & \m } $$
where both horizontal arrows represent branched coverings having degree
\[\nu_n(p)\] for \[p\ge 2\], but having degree \[\nu_n(1)+1\] for \[p=1\].
(To simplify the notation, the subscript \[n\] on \[\bic\] and \[\m\] has
been suppressed.) Start with the variety
$$	{\bf V}_p~\subset~\bic\times\widehat\C $$ 
consisting of all pairs \[(f,z)\] where \[f\]
is a bicritical map of degree \[n\] and \[z\] is a point of \[\widehat\C\]
satisfying \[f^{\circ p}(z)=z\]. Here \[p\] can be any positive integer.
Note that \[{\bf V}_q\subset {\bf V}_p\] for every divisor \[q\] of \[p\], Let
\[\PER_p\subset\bic_n\times\widehat\C\]
be the Zariski closure of the set
$$	{\bf V}_p\ssm\bigcup_{q|p,\,q<p}{\bf V}_q $$
consisting of points in \[{\bf V}_p\] which do not belong to
\[{\bf V}_q\] for any \[q<p\]. Clearly \[{\bf V}_p\] can be expressed as
the union of \[\PER_q\] as \[q\] ranges over all divisors of \[p\]
including \[p\] itself.

By definition, two points \[(f,z)\] and \[(g,w)\] of \[\PER_p\]
are {\bit conjugate\/} if there is a M\"obius automorphism \[\phi\] of
\[\widehat\C\] so that \[g=\phi\circ f\circ\phi^{-1}\] and \[w=\phi(z)\].
The orbifold consisting of all conjugacy classes \[((f,z))\] of points of
\[\PER_p\] will be denoted by \[\per_p\]. Evidently the
correspondence \[((f,z))\mapsto(f)\] yields the required branched covering
\[\per_p\to\m\].\ss

{\bf Remark.} These varieties \[\PER_p\] and \[\per_p\] are
actually irreducible.
The analogous statement for unicritical polynomials was proved by Bousch [Bo]
in the quadratic case and by Lau and Schleicher [LS1] for all degrees.
Irreducibility of these varieties for bicritical rational maps follows,
since any irreducible component of \[\per_p\]
must intersect the locus of polynomial maps.\ss

{\bf Definition.} The {\bit multiplier map}
$$	\lambda:\per_p~\to~\C $$
carries each \[((f,z))\] to the derivative
of \[f^{\circ p}\] at \[z\]. Let \[F_{\lambda_0}\subset\per_p\]
be the fiber, consisting of all \[((f,z))\] with \[\lambda((f,z))=\lambda_0\],
and let \[\per_p(\lambda_0)\] be its image under the projection to \[\m\].
(Intuitively, a point of \[\per_p(\lambda)\] is a conjugacy class of maps
\[(f)\] which posses a period \[p\] orbit of multiplier \[\lambda\],
while a point of \[F_{\lambda}\] consists of such an \[(f)\] together
with a specific choice of period \[p\] point for a representative map \[f\].)
\ss

{\bf Remark.}
Each \[k\] modulo \[p\] in the cyclic group \[\Z_p=\Z/p\Z\] acts on the
varieties \[\PER_p\to\per_p\] 
by the correspondence \[(f,z)\mapsto(f\,,\,f^{\circ k}(z))\]. Thus we obtain
quotient orbifolds and a larger commutative diagram
$$ \matrix{	\PER_p & \to &\PER_p/\Z_p & \to & \bic\cr
	\downarrow && \downarrow && \downarrow\cr
	\per_p & \to &\per_p/\Z_p &\to & \m \cr
	&&\downarrow\lambda\cr
	&&\C~~ } $$
where now the horizontal arrows represent branched coverings of degree
\[p\] on the left and \[\nu_n(p)/p\] on the right. Starting with any
\[\lambda_0\in\C\] we can form the preimage in \[\per_p/\Z_p\] and then
project to the subset \[\per_p(\lambda_0)\subset\m\].\ss

{\QP{\bf Theorem 8.1.} \it Let \[\overline\per_p(\lambda)\] be the
closure of \[\per_p(\lambda)\] within the extended moduli space \[\E\].
If \[\lambda\ne 0\], then \[\overline\per_p(\lambda)\]
is a compact, not necessarily irreducible, algebraic
curve in the extended moduli space \[\E\]. The projection map \[((f,z))
\mapsto X(f)\] carries \[\overline\per_p(\lambda)\] onto \[\widehat\C\]
with degree \[\nu_n(p)/n\]. Equivalently, the number of intersections
of \[\overline\per_p(\lambda)\] with any fiber \[L_X\] of the fibration
\[X:\E\to\widehat\C\], counted with multiplicity, is independent of \[\lambda\]
and \[X\], being equal to \[\nu_n(p)/n\]. The same statements are true for
\[\lambda=0\] and \[p\ge 2\] provided that the point set \[\per_p(0)\] is
counted with multiplicity \[n-1\].\ss}

In other words, for \[p\ge 2\],
it is asserted that the set \[\per_p(0)\] intersects a generic fiber \[L_X\]
in \[\nu_n(p)/(n(n-1))\] distinct points, each of which must be counted with
multiplicity \[n-1\] in order to get the correct count of \[\nu_n(p)/n\].
(Here is one intuitive explanation for this multiplicity: For
\[(f)\in\per_p(0)\], suppose that we perturb \[f\] within
the much larger space consisting of all rational maps of degree \[n\].
Generically, the periodic critical point will split up into \[n-1\] nearby
critical points, any one of which can be periodic. Thus the locus
\[\per_p(0)\subset\m\] splits up locally, in this larger context, into \[n-1\]
nearby sheets.)\ss

{\bf Outline Proof of 8.1.} We need only consider the case \[p\ge 2\], since
the period one case has been discussed in \S2 and \S6.
Since \[\per_p(\lambda)\] is an algebraic
curve in \[\m\cong\C^2\], it can be defined by a single polynomial equation
in the coordinates \[X\] and \[\sigma_n\]. Substituting \[\hat X=1/X\] and
\[\hat\sigma_1=\sigma_n/(n^{n+1}X^{n-1})\], and multiplying through by an
appropriate power of \[\hat X\], we obtain a corresponding polynomial
equation relating \[\hat X\] and \[\hat\sigma_1\]. This shows that
\[\overline\per_p(\lambda)\] is an algebraic curve in the variety \[\E\].
If \[|\lambda|>1\] so that the associated periodic orbit is contained
in the Julia set, then
it follows from Remark 7.3 that this curve is contained in a compact
subset of \[\E\], and hence is itself compact. On the other hand, if
\[|\lambda|\le 1\] with period \[p\ge 2\], then \[\per_p(\lambda)\] is
contained in the connectedness locus, which has compact closure within \[\E\].

To compute the degree of the projection map \[\overline
 \per_p(\lambda)\to\widehat\C\],
we first look at the exceptional special case \[\lambda=0\] and count the
number of intersections of \[\per_p(0)\]
with the fiber \[L_0\]. If \[f_b(z)=z^n+b\] with
invariants \[X=0\] and \[Y=b^{n-1}\], then the equation
\[f_b^{\circ p}(0)=0\]
has degree \[n^{p-1}\] in the unknown \[b\], so there are \[n^{p-1}\]
solutions \[b\], counted with multiplicity. Subtracting off the numbers
of solutions for proper divisors of \[p\], we see that there are
\[\nu_n(p)/n\] choices for \[b\]. Since we assume that \[p\ge 2\],
it follows that \[b\ne 0\].
Hence there are \[n-1\] choices of \[b\] for every choice of \[Y=b^{n-1}\].
Thus the projection \[X:\overline\per_p(0)\to\widehat\C\] has degree
\[\nu_n(p)/(n(n-1))\] for \[p\ge 2\]. The first few values are as follows:\ms

\halign{\quad\hfil#\hfil&\qquad\hfil#\hfil&\quad\hfil#\hfil&\quad\hfil#\hfil&\quad\hfil#\hfil&
\quad\hfil#\hfil&\quad\hfil#\hfil\cr
period \[p\] & 2&3&4&5&6\cr
\[\nu_n(p)/(n(n-1))\] &1&  \[n+1\]& \[n(n+1)\]& \[(n+1)(n^2+1)\]&
 \[(n+1)(n^3+n-1)\]\cr}
\ms

\ni
To check what happens as we perturb the multiplier \[\lambda\] away from zero,
we can simply apply the following result. See [LS1, \S2.2].

{\QP{\bf Proposition (Lau and Schleicher).} \it If \[W\] is any hyperbolic
component in the \[b$-parameter plane for the family of maps \[z\mapsto
 z^n+b\], then the multiplier map \[W\to\D\] is an \[(n-1)$-fold branched
covering, ramified only at the unique inverse image of zero.\ss}

\ni Thus, as we perturb \[\lambda\] away from zero, each intersection of
\[\per_p(0)\] with \[L_0\] splits\break into \[n-1\] distinct intersection,
and we obtain the required count of \[\nu_n(p)/n\] points in \[\per_p(\lambda)
\cap L_0\] for \[\lambda\ne 0\].\QED\ss

{\bf Example 8.2.} For the special case \[p=2\],
it follows that \[\overline\per_2(0)\] is the image of
a smooth section of the fibration \[X:\E\to\widehat\C\]. We can compute
this section explicitly as follows. Using the normal form (1) with \[ad-bc=1\],
recall that
$$	X\=bc~,\quad X+1\=ad~,\quad Y_1=a^{n+1}b^{n-1}~, $$
and that \[f(0)=b/d\]. Thus the critical point \[0\] has period exactly two
if and only if
$$	f(b/d)\={ab^n+bd^n\over cb^n+d^{n+1}} $$
is zero, but \[b\ne 0\]. This yields the equation \[ab^{n-1}+d^n=0\].
Multiplying by \[a^n\], we obtain \[Y_1+(X+1)^n=0\]. If \[Y_1\ne 0\],
it follows that
$$	Y_2\=(X+1)^{n+1}X^{n-1}/Y_1\= -(X+1)x^{n-1}~, $$
and therefore that
$$	Y\=Y_1+Y_2\= -(X+1)^n-(X+1)X^{n-1} ~.$$
On the other hand, if \[Y_1=0\] then it is easy to check that \[X+1=Y_2=Y=0\]
also, so that the equation is still satisfied. Similarly, if the other
critical point \[\infty\] is periodic, then interchanging the roles of
\[Y_1\] and \[Y_2\] we obtain the same equation.\ss

{\bf Remark 8.3.}
It seems likely that the curve \[\overline\per_p(\lambda)\] is usually
irreducible, however there are certainly exceptions. As an example,
if \[n\] is odd then it was noted in \S1 that the curve
\[\overline\per_2(n^2)\] is reducible, since it contains the half symmetry
locus \[\Sigma_-\] as an irreducible component. (This period two curve
is strictly larger that \[\Sigma_-\] since it
has degree \[n-1\ge 2\] over \[\widehat\C\] while \[\Sigma_-\] has
degree one.)

Here is another example. If \[1\le q<p\] is a proper divisor of \[p\],
and if \[\xi\] is a primitive \[(p/q)$-th root of unity, then it is not
hard to show that
$$	\per_q(\xi)~\subset~\per_p(1) $$
since, under a small perturbation of \[f\], any period \[q\] orbit of
multiplier \[\xi\] will split off a period \[p\] orbit with multiplier
close to \[1\]. In general,
it follows that \[\per_p(1)\] is reducible for \[p\ge 2\]. The only exception
occurs in the degree two case, where \[\per_2(1)\] actually coincides
with the curve \[\per_1(-1)\]. (Thus a quadratic rational map
cannot have any period two orbit with multiplier \[1\].)

We can describe the intersection of \[\overline\per_p(\eta)\] with the line
at infinity rather precisely as follows.

{\QP{\bf Theorem 8.4.} \it For \[p\ge 2\] the intersection \[\overline\per_p
(\eta)\cap L_\infty\] depends only on the period \[p\] and not on the
multiplier \[\eta\]. This intersection
is a finite set consisting of points with
coordinate \[\hat\sigma_1\] of the form \[\lambda+\lambda^{-1}\] where
\[\lambda\] is a \[q$-th root of unity for some \[q\le p\].\ss}

{\bf Proof.} First consider the case \[\eta=0\], and suppose that
\[\lambda^q\ne 1\] for \[1\le q\le p\]. In other word, we suppose that the complex
numbers
$$	1\,,~\lambda\,,~\lambda^2\,,~\ldots\,,~\lambda^p $$
are all distinct. Let \[(f)\in\m\] be a conjugacy class
with \[\hat X\] close to zero and with \[\hat\sigma_1\]
close to \[\lambda+\lambda^{-1}\]. Using 7.5 we can find
a normal form for \[(f)\] so that both critical points are close to \[1\],
with the first \[p+1\] points of both critical orbits close to
$$	1~\mapsto~\lambda~\mapsto~\lambda^2~\mapsto~\cdots~\mapsto~\lambda^p~\ne~1~. $$
Thus, if the pair \[(\hat X\,,\,\hat\sigma_1)\] is sufficiently close
to \[(0\,,\,\lambda+\lambda^{-1})\], then it follows that neither critical
orbit can have period \[p\], hence \[(f)\not\in\per_p(0)\], as required.

(This argument does not apply if \[\lambda^q=1\] for some \[q<p\].
For then the \[q$-th forward
image of the critical point is close to \[1\], lying in a small region whose
image under \[f\] covers the entire Riemann sphere. Hence the further orbit
of the critical points cannot be predicted without more information.)

Now consider \[\per_p(\eta)\] with variable \[\eta\]. Consider the
coordinate \[\hat\sigma_1\] for the various
points of \[\per_p(\eta)\cap L_\infty\]. The elementary symmetric functions
of these \[\nu_n(p)/n\] (not usually distinct)
coordinate values can be expressed as holomorphic functions of
the parameter \[\eta\]. But if \[|\eta|<1\]and \[p\ge 2\]
then \[\per_p(\eta)\] is contained in the connectedness locus, so that
$$	\overline\per_p(\eta)\cap L_\infty~\subset~\overline\CL\cap L_\infty
	~\cong~[-2,2]~\subset~\R~. $$
In other words, every one of the \[\nu_n(p)/n\] values of \[\hat\sigma_1\]
is real, hence the elementary symmetric functions of these coordinate values
are also real. But a
holomorphic function from \[\eta\in\C\] to \[\C\] which takes only
real values throughout the open disk \[|\eta|<1\] must be constant. This
proves the the intersection is independent of \[\eta\] for all
\[\eta\in\C\], and completes the proof of 8.4.\QED
\ss

{\bf Remark 8.5.} More precisely, for \[p\ge 2\]
it can be conjectured that the point
of \[L_\infty\] with coordinate
\[\hat\sigma_1=\lambda+\lambda^{-1}\] belongs to \[\overline\per_p(\eta)\]
if and only if \[\lambda\] is a primitive \[q$-th root of unity for
some \[2\le q\le p\], or in other words if and only if
$$	\hat\sigma_1\=2\cos(2\pi r/q)\qquad{\rm for~some~integers}\qquad
	0<r<q\le p~. $$
Thus it is believed that the point of \[L_\infty\] with coordinate
\[\hat\sigma_1=2\] does not lie on any \[\overline\per_p(\eta)\] with
\[p\ge 2\]. (Compare [Sti] and [R3], which give a more precise description of
the curves \[\overline\per_p(0)\] near \[L_\infty\] in the degree two case,
and see also [M2].)\ss

\def\lr{\longrightarrow}
\def\lrm{~\buildrel n-1\over\lr}
\def\lrn{~\buildrel n\over\lr}

\midinsert\ss
\cl{\hbox{\vrule\vbox{\tenpoint
\halign{\quad#\hfil\quad &\hfil #\hfil\quad&\hfil #\hfil & \hfil #\hfil &
 \hfil #\hfil\cr
\noalign{\hrule\ss}
\[p\] & \[r/q\] & limiting orbit of \[g(1)\approx\lambda\]
 & number & total\cr
\noalign{\ss\ss\hrule}\cr
2 & 1/2 & \[\lambda\lr \;1\] & 1 \cr
&&&& 1 \cr
\noalign{\ss\hrule}\cr
3 & \[1/3,\, 2/3\] & \[\lambda\lr\lambda^2\lr \;1\] & \[2\times 1\]\cr
& 1/2 & \[\lambda\lr \;1 \lrm \;1\] & \[n-1\]\cr
&&&& \[n+1\]\cr
\noalign{\ss\hrule}\cr
4 & \[1/4,\,3/4\] & \[ \lambda\lr\lambda^2\lr\lambda^3\lr\;1\] &
 \[2\times 1\]\cr
& \[1/3,\,2/3\] & \[ \lambda\lr\lambda^2\lr\;1\lrm\;1\] & \[2\times(n-1)\]\cr
& \[1/2\] & \[\lambda\lr\;1 \lrn \lambda\lr\;1\] & \[n\]\cr
& \[1/2\] & \[\lambda\lr\;1 \lrm\;1\lrm\;1\] & \[(n-1)^2\]\cr
&&&& \[n^2+n+1\]\cr
\noalign{\ss\hrule}\cr
5 &  \[1/5,\,2/5,\,3/5,\,4/5\] &\[ \lambda\lr\lambda^2\lr\lambda^3\lr\lambda^4
 \lr\;1\] & \[4\times 1\]\cr
& \[1/4,\,3/4\] & \[ \lambda\lr\lambda^2\lr\lambda^3\lr\;1\lrm\;1\] &
 \[2\times(n-1)\]\cr
& \[1/3,\,2/3\] & \[ \lambda\lr\lambda^2\lr\;1\lrm\lambda^2\lr\;1\] &
 \[2\times(n-1)\]\cr
& \[1/3,\,2/3\] & \[ \lambda\lr\lambda^2\lr\;1\lrm\;1\lrm\;1\] &
 \[~2\times(n-1)^2\]\cr
& 1/2 & \[\lambda\lr\;1\lrn\lambda\lr\;1\lrm\;1\] & \[n(n-1)\]\cr
& 1/2 & \[\lambda\lr\;1 \lrm\;1\lrn\lambda\lr\;1\] & \[n(n-1)\]\cr
& 1/2 & \[\lambda\lr\;1\lrm\;1\lrm\;1\lrm\;1\] & \[(n-1)^3\]\cr
&&&& \[n^3+n^2+n+1\]\cr
\noalign{\ss\hrule}\cr}}\vrule}}\ms

{\QP\bit Table listing the numbers of points \[(g)\in L_X\] such that
\[g\] admits a periodic critical orbit of period dividing \[p\],
with \[|X|\] very large,
grouped by the conjectured limiting behavior of this orbit as \[X\to\infty\].
By definition, \[\lambda=e^{2\pi ir/q}\] in each row.
The sum on the right must be equal to \[1+n+n^2+\cdots+n^{p-2}\]
for all values of \[p\]. \par}
\endinsert

Following is an intuitive argument which attempts to justify this statement,
and also to compute the multiplicities of the various points of \[\per_p(0)
\cap L_\infty\]. I suspect that it could be made into a rigorous
proof, but certainly have not done so.
Let us count solutions to the equation \[g(z)=k\], where \[k\in\C\]
is some given constant and where \[g\]
is the bicritical map of 7.5, with \[|X|\] very large so that \[g(z)\approx
\lambda z\] outside of a small neighborhood of \[1\]. If \[k\] is bounded away
from \[\lambda\], then there is a unique solution \[z\approx k/\lambda\]
which is bounded away from \[1\]. Hence the remaining \[n-1\] solutions
must all lie
in some small neighborhood of \[1\]. Let us express these facts symbolically,
in the limit as \[|X|\to\infty\], by writing
$$	k/\lambda\lr k\qquad{\rm but}\qquad 1\lrm k~, \qquad{\rm
     for}\qquad k\ne\lambda~. $$
In the special case \[k=\lambda\], all \[n\]
of the solutions must be very close to \[1\], so we write
$$	1\lrn\lambda~. $$
The numbers on the right in the table are computed by muliplying these
factors of \[n-1\] or \[n\] by the number of values of \[r/q\]
which are listed on the left. If \[N(p,q)\] denotes the number of solutions
of period dividing \[p\] which are computed in this way, corresponding to the
case where \[\lambda\] is a primitive \[q$-th root of unity, \[q\ge 2\], then
it is not hard to check that these numbers can be computed recursively as
follows:
$$	N(p,q)\=\cases{ 0 & if \[p<q\]\cr
	\phi(q) & if \[p=q\]\cr
	nN(p-q,q)+\sum_{j=1}^{q-1} (n-1)N(p-j,q) & if \[p>q\],} $$
where \[\phi(q)\], the Euler \[\phi$-function, is the number of primitive
\[q$-th roots of unity. If \[N(p,q)\] is defined by this recursion, then
it seems empirically that the sum over \[q\] is given by
$$	\sum_{q\ge 2} N(p,q)\= 1+n+n^2+\cdots+n^{p-2}~.$$
For any given \[p\], this is easy to
check by computer; but I don't know a general proof.

On the other hand, for every \[p\ge 2\], the actual number of
solutions of period exactly \[p\] is known to be \[\nu_n(p)/(n(n-1))\].
There are no solutions of period 1, hence
the number of solutions with period \[d\] dividing \[p\] must be
equal to
$$	\sum_{d|p\,,\;d\ne 1}{\nu_n(d)\over n(n-1)}\={n^p-n\over n(n-1)}\=
	1+n+n^2+\cdots+n^{p-2}~.$$
This agrees precisely with the intuitive computation, which supports the
conjecture that the computed
values \[N(p,q)\] are indeed correct, and that there are no solutions
at all corresponding to the case \[\lambda=1\].

It is interesting to note that, for each \[p\], the contribution from
\[q=2\] is by far the largest. Indeed, for large \[n\],
the contribution of \[(n-1)^{p-2}\],
corresponding to a critical orbit for \[g\] which tends asymptotically to
$$	1~\mapsto~-1~\mapsto~1~\mapsto~1~\mapsto~\cdots~\mapsto~1~, $$
contributes much more than all of the other asymptotic behaviors combined.
The associated point of \[\per_p(0)\cap L_\infty\] has
coordinate \[\hat\sigma_1=\lambda+\lambda^{-1}=-2\].
\bs\bs

\def\bic{{\bf bicrit}}
\def\per{{\rm Per}}
\def\D{{\bf D}}

\def\CL{{\cal C}}
\def\SL{{\cal S}}
\def\M{ E} 
\def \iso{{\,\buildrel \simeq\over\longrightarrow\,}}


\headline={%
    \ifnum\pageno > 1 {%
        \hss\tenrm\ifodd\pageno A. NO HERMAN RINGS\hss\folio
        \else\folio\hfill MAPS WITH TWO CRITICAL POINTS\hfill\fi}%
    \else \hss\vbox{}\hss
    \fi}

\cl{\bf Appendix A. No Herman Rings}\ss

The following statement is a straightforward consequence
of results due to Shishikura. 

{\QP{\bf Theorem A.1.} \it A rational map $f$ with only two critical
points cannot have any Herman rings.\ms}

The following proof is closely modeled on Shishikura's argument (compare
[Sh, \S8], together with the remark on page 4 of [Sh]),
although it avoids the use of surgery constructions.\ss

Let $f:\widehat\C\to\widehat\C$ be a rational map of degree $d\ge 2$ which
possesses a cycle of Herman rings $H=H_1\cup\cdots\cup H_p$ of period
$p\ge 1$. That is, assume that the $H_i$ are disjoint annuli
with $\partial H_i\subset J(f)$, where $i$ ranges over the group
of integers modulo $p$, 
and assume that $f$ maps each $H_i$ diffeomorphically onto $H_{i+1}$.
Choosing some base point $b\in H$, the closure
of the forward orbit of $b$ is a union $\Gamma=\Gamma_1\cup\cdots\cup\Gamma_p$
of smooth simple closed curves, where $\Gamma_i\subset H_i$ and
$f(\Gamma_i)=\Gamma_{i+1}$.
Since $f(\Gamma)=\Gamma$, we have
$$	\Gamma~\subset~ f^{-1}(\Gamma)~\subset~f^{-2}(\Gamma)~\subset~\cdots~.
  $$
(If $\Gamma$ is chosen to be disjoint from the postcritical set, then each
$f^{-n}(\Gamma)$ will be a union of at most $pd^n$ disjoint simple closed
curves.)
The connected components of $\widehat\C\ssm f^{-n}(\Gamma)$ are open sets which
will be called {\bit Shishikura puzzle pieces\/} of {\bit level\/} $n$.
Evidently $f$ maps each puzzle piece of level $n>0$ onto a puzzle piece of
level $n-1$ by a possibly branched covering map.

Each  $\Gamma_i$  separates the Riemann sphere into two disks. In order to
label these disks, choose orientations for the
loops  $\Gamma_i$  compatible with the mapping. Then for each  $\Gamma_i$
 we can speak of the  disk $D^{\rm L}_i$ to the {\bit left\/} of  $\Gamma_i$
with boundary  $\Gamma_i$, and the disk $D^{\rm R}_i$ to the
{\bit right\/} of  $\Gamma_i$  with boundary  $\Gamma_i$ .
Fixing some level $n$, note that
each  $\Gamma_i$  lies on the boundary of exactly two puzzle pieces. Let  $L_i
=L^{(n)}_i\subset D^{\rm L}_i$ be the adjacent piece to the left of  $\Gamma_i$
and  $R_i=R^{(n)}_i\subset D^{\rm R}_i$  the adjacent piece
to the right. Define the {\bit level $n$ left neighborhood\/} $L^{(n)}$
of $\Gamma$
to be the union $L^{(n)}_1\cup\cdots\cup L^{(n)}_p$, and define the {\bit right
neighborhood\/} $R^{(n)}=R^{(n)}_1\cup\cdots\cup R^{(n)}_p$ similarly. Then
we have the following. (See [Sh,\S7].)

{\QP{\bf Theorem A.2.} \it The left neighborhood $L^{(n)}$ of $\Gamma$
contains at least one critical point of $f$, and similarly the right
neighborhood $R^{(n)}$ contains at least one critical point.\ms}

{\bf Proof.} We may assume that $n\ge 1$.
If $L_1\cup\cdots\cup L_p$ contains no critical point, then we will show
that $f$ must map each $L_i$ onto $L_{i+1}$, and hence that
$f^{\circ p}$ must map $L_i$ onto itself. This would
imply that each $L_i$ is contained in the Fatou set $\C\ssm J(f)$.
But that is impossible, since $L_i$ intersects the boundary of the
ring $H_i$, which is contained in the Julia set. This contradiction will
prove the Theorem.

For each left hand puzzle piece $L_i^{(n)}$ we have the following diagram.
$$ L_i^{(n)}~{\buildrel f\over\longrightarrow}~
 L_{i+1}^{(n-1)}~\supset~ L_{i+1}^{(n)}~. $$
If there is no critical point in $L_i^{(n)}$, then $f$ maps $L_i^{(n)}$ onto
$L_{i+1}^{(n-1)}$  by a covering map (which may be one-to-one).
Hence, if we pass to universal covering spaces, this diagram takes the form
$$	\tilde L_i^{(n)} ~\cong
~\tilde L_{i+1}^{(n-1)}~\longleftarrow~\tilde L_{i+1}^{(n)}~. $$
In other words, we can choose a single valued holomorphic branch of $f^{-1}$
which maps $\tilde L_{i+1}^{(n)}$ into $\tilde L_i^{(n)}$.
Each puzzle piece $L_i^{(n)}$ is a hyperbolic Riemann surface, and hence
has a Poincar\'e metric. Fixing $n$, for any smooth curve segment
$\alpha:[0,1]\to L_i^{(n)}$,
let $\ell_i(\alpha)$ be the Poincar\'e arclength.
Then one of the following two possibilities must hold:

{\bf Case 1.} If $f(L_i^{(n)})$ is precisely equal to $L_{i+1}^{(n)}$, 
then this branch of $f^{-1}$ is a Poincar\'e isometry from
$\tilde L_{i+1}^{(n)}$ onto $\tilde L_i^{(n)}$. In this case $f$
preserves Poincar\'e arclength, so that $$\ell_{i+1}(f\circ\alpha)\=
\ell_i(\alpha)$$ for every smooth curve in $L_i^{(n)}$.

{\bf Case 2.} If $f(L_i^{(n)})$ is strictly larger than $L_{i+1}^{(n)}$, then
this branch of $f^{-1}$ is strictly distance reducing. Hence $f$
must strictly increase Poincar\'e arclength, in the sense that
$$	\ell_{i+1}(f\circ\alpha)~>~\ell_i(\alpha) $$
for any non-constant curve in $L_i^{(n)}$ which maps into $L_{i+1}^{(n)}$.

Now choose some orbit closure $\Gamma'\subset H\cap L^{(n)}$, and let
$\alpha:\R/\Z\to \Gamma'\cap H_0$ be a smooth parametrization of one
connected component. Then $f^{\circ i}\circ\alpha$ is a smooth parametrization
of $\Gamma'\cap H_i$. Comparing the discussion above, we have
$$	\ell_0(\alpha)~\le~\ell_1(f\circ\alpha)~\le~\ell_2(f^{\circ 2}
 \circ\alpha)~\le~\cdots~\le~\ell_p(f^{\circ p}\circ\alpha)\=\ell_0(\alpha)~,$$
since $f^{\circ p}\circ\alpha$ and $\alpha$ parametrize the same loop.
Thus equality must hold throughout, and Case 2 can never occur.
This proves that $f$ must map every $L_i^{(n)}$ onto $L_{i+1}^{(n)}$.
As noted above, this leads to a contradiction. Therefore
$L^{(n)}$ must contain at least one critical point, which completes the
proof of A.2.\QED\ms

{\bf Remark A.3.} {\it If $n\ge p-1$, then $L^{(n)}$ is disjoint from
$R^{(n)}$.\/} This statement is clear when $p=1$, since $L_i^{(n)}\cap
R_i^{(n)}=\emptyset$. If
$L^{(p-1)}\cap R^{(p-1)}~\ne~\emptyset$ with $p>1$, then
$L_{i_0}^{(p-1)}$ would be precisely equal to $R_{i_0+\delta}^{(p-1)}$
for some $i_0$ and some $\delta\not\equiv 0$.
This would imply that $L_{i_0+1}^{(p-2)}=R_{i_0+1+\delta}^{(p-2)}$. Continuing
inductively it would follow that $L_i^{(0)}=R_{i+\delta}^{(0)}$ for every
$i$. But this would imply that
the entire disk $D^{\rm R}_i$ to the right of $\Gamma_i$ is contained in
$D^{\rm R}_{i+\delta}$. Therefore
$$	D^{\rm R}_0~\subset~D^{\rm R}_\delta~\subset~ D^{\rm R}_{2\delta}~\subset~\cdots
	~\subset~D^{\rm R}_{p\delta}\=D^{\rm R}_0~,$$
which is impossible.\QED\ms

Now let us specialize to the case of a rational map of degree $d$ with only two
critical points. Note the following basic observation. If a simple closed
curve $\Gamma_i\subset\widehat\C$ bounds
a closed disk $\overline D$ which is disjoint
from the two critical values, then $f^{-1}(\overline D)$ is the union of $d$
disjoint topological disks, each of which contains no critical point
and maps homeomorphically onto $\overline D$. On the other hand, if $\Gamma_i$
separates the two critical values, then $f^{-1}(\Gamma_i)$ is a single simple
closed curve, which maps onto $\Gamma_i$ by a $d$-fold covering map.
Evidently this second case can never occur for the loops $\Gamma_i$
associated with a cycle of Herman rings.

It will be convenient to put
one of the two critical points at infinity. Then one of the two components
of the complement of any  $\Gamma_i$  is bounded, and maps diffeomorphically onto
its image under  $f$ . I will call this component the {\bit inside\/}
$D_i^{\rm in}$ of $\Gamma_i$ . The other 
component $D_i^{\rm out}$ is unbounded, and contains both critical points.

Call  $\Gamma_k$  a {\bit separating\/} loop if
it separates the critical points from the critical values. {\it Then
the inside $D_i^{\rm in}$ of each  $\Gamma_i$  maps diffeomorphically onto:

\hskip 1in the outside   $D^{\rm out}_{i+1}$  if  $\Gamma_{i+1}$  is separating,

\hskip 1in
the inside  $D^{\rm in}_{i+1}$  if  $\Gamma_{i+1}$  is non-separating.}

\ni Call  $\Gamma_k$  {\bit minimal\/}  if the union $\Gamma$ is disjoint
from the open disk $D_k^{\rm in}$.

{\QP{\bf Lemma A.4.} \it  There exists a (necessarily unique)
$\Gamma_k$  which is both separating and minimal.\ms}

{\bf Proof} [Sh, \S8]. There certainly exists at least one minimal  $\Gamma_i$ .
If every  $\Gamma_i$  were minimal and non-separating, then the inside of every
$\Gamma_i$  would map diffeomorphically onto the inside of  $\Gamma_{i+1}$  and we would
have a cycle of Siegel disks rather than Herman rings. Therefore, either
there is some  $\Gamma_i$  which is minimal and separating, as required, or else
there is some  $\Gamma_i$  which is not minimal. In the latter case, we can
choose some non-minimal  $\Gamma_i$  so that  $\Gamma_{i+1}$  is minimal. Then
$D_i^{\rm in}$, which contains other $\Gamma_j$,  cannot map diffeomorphically 
onto $D^{\rm in}_{i+1}$ which does not.
Hence in this case  $\Gamma_{i+1}$  must be a minimal separating loop, as
required.  \QED\ss

{\bf Proof of A.1.} Let  $\Gamma_k$  be the loop of A.4.
Since $\Gamma_k$ is minimal,
the closure $\overline D_k^{\rm out}$ contains the union
$\Gamma=\Gamma_1\cup\cdots\cup\Gamma_p$. Since $\Gamma_k$ is separating, the
preimage $f^{-1}(\overline D_k^{\rm out})$
is a union of $d$ disjoint bounded closed disks. Evidently this union
contains $f^{-1}(\Gamma)$ and has boundary $f^{-1}(\Gamma_k)$.
It follows that the complementary region
$f^{-1}(D_k^{\rm in})$ is the unique unbounded puzzle piece of level one.
This puzzle piece contains both critical points, and
has no  $\Gamma_i$  other than  $\Gamma_{k-1}$  on its boundary.
Therefore the level one
puzzle piece to the inner side of  $\Gamma_{k-1}$  is bounded and contains no
critical point, while for any  $i\ne k-1$ both of the pieces $L^{(1)}_i$
and  $R^{(1)}_i$  are bounded and contain no critical point. Thus $L^{(1)}$
and $R^{(1)}$ cannot both contain critical points, which contradicts
A.2 and completes the proof that a bicritical map has no Herman ring. \QED
\bs\bs

\headline={%
    \ifnum\pageno > 1 {%
        \hss\tenrm\ifodd\pageno B. TOTALLY DISCONNECTED JULIA SETS\hss\folio	        \else\folio\hfill MAPS WITH TWO CRITICAL POINTS\hfill\fi}%
    \else \hss\vbox{}\hss
    \fi}

\cl{\bf Appendix B. Totally Disconnected Julia Sets.}\ms

We will first prove the following result, with no restriction on the number
of critical points. Let \[f\] be a rational function of degree \[n\ge 2\].

{\QP{\bf Theorem B.1.} \it The Julia set \[J\] of \[f\] is totally
disconnected and contains no critical point  if and only if
all of the critical values of \[f\] lie in a single Fatou component.\ms}

In the hyperbolic case, this was proved by Rees [R]. I attempted to extend
the argument to the parabolic case in [M2], however the details of the
argument were not quite right. The following proof is rather different and
perhaps easier. It will be based on the following ideas.\ss

{\bf Definition.}
A map \[g:K\to K'\] between metric spaces is called {\bit locally
distance increasing\/} if every point of \[K\] has a neighborhood \[V\] so that
$$	{\rm dist}(g(x),g(y))~>~{\rm dist}(x,y) \eqno (24)$$
for all \[x\ne y\] in \[V\]. If \[K\] is compact, then an equivalent
condition is that for some \[\epsilon>0\] we have \[(24)\] whenever
\[0<{\rm dist}(x,y)<\epsilon\].\ss

{\QP{\bf Lemma B.2.} \it Let \[g:K\to K\] be a locally distance increasing map
from a compact metric space to itself. If \[g\] is injective on each connected
component of \[K\], then the space \[K\] must be totally disconnected.\ms}

{\bf Proof.} First note the following:
{\it There exists
a constant \[\epsilon'>0\] so that the inequality \[(24)\] holds whenever
\[x\] and \[y\] are distinct points which belong
to the same connected component of \[K\] and satisfy \[{\rm dist}
\big(g(x),g(y)\big)<\epsilon'\].} For otherwise, with \[\epsilon\]
as above, there would exist pairs \[(x_i\,,\,y_i)\] belonging to the
same component \[K_i\] of \[K\]
so that \[{\rm dist}\big(g(x_i)\,,\,g(y_i)\big)\] converges to zero as
\[i\to\infty\], but with \[{\rm dist}(x_i\,,\,y_i)\ge\epsilon\]. After
passing to a subsequence, we may assume that the points \[x_i\] converge
to some point \[x\] and that the \[y_i\] converge to some \[y\], where
\[{\rm dist}(x,y)\ge\epsilon\] and where \[g(x)=g(y)\].
But this is impossible since \[x\] and \[y\] belong to the
same connected component of \[K\]. In fact, if \[L\] denotes the set of all
accumulation points of \[K_i\] as \[i\to\infty\], then \[L\] is a connected
subset of \[K\] containing both \[x\] and \[y\].\ss

The {\bit diameter\/} \[{\rm diam}(K)\ge 0\] of a non-vacuous compact
metric space is defined to be the maximum distance between two of its points.
More generally, for any integer \[m\ge 1\] define \[{\rm diam}_m(K)\] to be the
largest number \[\delta\] so that there exist \[m+1\] points \[x_0\,,\,x_1\,,
\,\ldots\,,\,x_m\] in \[K\] which are {\bit \[\delta$-separated\/} in the sense
that \[{\rm dist}(x_i\,,\,x_j)\ge\delta\] for \[i\ne j\]. (Think of
\[m+1\] strongly repelling points, which try to get as far as possible
from each other.) Thus
$$	{\rm diam}(K)={\rm diam}_1(K)~\ge~{\rm diam}_2(K)~\ge~{\rm diam}_3(K)
 ~\ge~\cdots\ge 0~. $$
Note that \[{\rm diam }_m(K)=0\] if and only if
\[K\] is a finite set with at most \[m\] elements.\ss

The proof of B.2 now proceeds as follows.
Choose \[m\] large enough so that \[K\] can be covered by subsets \[X_1\,,\,
\ldots\,,\,X_m\] of
diameter \[<\epsilon'\]. It then follows that \[{\rm diam}_m(K)<\epsilon'\],
since a collection of \[\epsilon'$-separated points can have at most one point
in each \[X_i\]. Let \[\delta_{\rm max}\]
be the supremum of the numbers \[{\rm diam}_m(K_\alpha)\] as \[K_\alpha\]
ranges over all connected components of \[K\]. We want to prove that \[\delta_{\rm max}
=0\]. Otherwise, if \[\delta_{\rm max}>0\], we will obtain a contradiction by
constructing a ``largest'' component \[K_\mu\], with \[{\rm diam}_m(K_\mu)\]
equal to \[\delta_{\rm max}\],
and then showing that \[{\rm diam}_m\big(f(K_\mu)\big)>\delta_{\rm max}\].

Let \[K_i\] be a sequence of components of \[K\] such
that the numbers \[\delta_i={\rm diam}_m
(K_i)\] converge to the supremum \[\delta_{\rm max}\], and
choose points \[x_0(i)\,,\,x_1(i)\,,
\,\ldots\,,\,x_m(i)\in K_i\] which are \[\delta_i$-separated.
After passing to a subsequence, we may assume that each sequence
\[x_j(1)\,,\,x_j(2)\,,\,x_j(3)\,,\,\ldots\] converges to a limit \[x_j\in K\].
Let
\[L\] be the set of all accumulation points of \[\{K_i\}\] as \[i\to\infty\].
Then \[L\] is connected, and hence is contained in some connected component
\[K_\mu\]. Furthermore \[{\rm diam}_m(K_\mu)=\delta_{\rm max}\] since
\[L\] contains \[m+1\] points \[x_j\] which are \[\delta_{\rm max}$-separated.

Now, assuming that \[\delta_{\rm max}>0\], we will obtain a contradiction
by showing that\break \[{\rm diam}_m\big(f(K_\mu)\big)\] must be strictly
larger than \[\delta_{\rm max}\].
In fact, for each \[0\le i<j\le m\] we have either
$$	{\rm dist}\big(f(x_i)\,,\,f(x_j)\big)~>~{\rm dist}(x_i\,,\,x_j)
~\ge~\delta_{\rm max}~, $$
or else
$$	{\rm dist}\big(f(x_i)\,,\,f(x_j)\big)~\ge~\epsilon'~>
	~{\rm diam}_m(K)~\ge~\delta_{\rm max}~. $$
This contradiction proves that \[\delta_{\rm max}\] must be zero.
Hence every connected component \[K_\alpha\] is finite, and hence consists
of a single point. This proves B.2.\QED\ss

{\bf Proof of B.1.} In one direction this is clear, since a totally
disconnected set cannot separate the Riemann sphere. For the proof in the
other direction, suppose that all critical values lie in a single Fatou
component \[U\].
Choose a smoothly embedded closed disk \[\Delta^*\subset
U\] which contains all of the critical values in its interior. Let
\[\Delta=\widehat\C\ssm{\rm interior\,}(\Delta^*)\] be the complementary closed
disk which contains the Julia set.
Since \[\Delta\] contains no critical values, its pre-image \[f^{-1}(\Delta)\]
splits as a disjoint union
$$	f^{-1}(\Delta)\=\Delta_1\cup\cdots\cup\Delta_n $$
of smoothly
embedded closed disks, each of which maps homeomorphically onto \[\Delta\]
under \[f\]. Now every connected component \[J_\alpha\] of the Julia set
must be contained in one of the \[\Delta_j\], and hence must map
homeomorphically under \[f\]. Thus, to prove B.1, we need only show that
\[f\] restricted to the Julia set is locally distance increasing with respect
to a suitably chosen metric. In fact we will prove the following.
(For related statements, see [DH1] and [TY].)\ss

{\QP{\bf Lemma B.3.} \it A rational map \[f\], restricted to its Julia set
\[J\], is
locally distance increasing with respect to a suitably chosen metric if and
only if \[f\] has no critical points in \[J\].\ms}

{\bf Proof of B.3.} If there is a critical point \[c\] in the Julia set
\[J\], then \[f|_J\] is not even one-to-one near \[c\], so it certainly cannot
be locally distance increasing. Suppose then that there are no critical
points in \[J\].
Let \[P\] be the closure of the {\bit postcritical set\/} of \[f\] (the
union of all forward orbits of critical values).
If we exclude trivial cases where \[P\] has only two elements,
then each connected component of the complement \[U=\widehat\C\ssm P\] has
a well defined Poincar\'e metric. The associated Poincar\'e
distance function will be denoted by \[\dist_U(x,y)\].
Let \[U_0=f^{-1}(U)\subset U\]. Then \[f\] maps
\[U_0\] onto \[U\] by a covering map. Since \[U_0\] is a proper subset
of \[U\], it follows that
$$ \dist_U(x,y)~<~\dist_{U_0}(x,y)\=\dist_U(f(x)\,,\,f(y)) $$
for any two points \[x\ne y\] in \[U_0\] which are sufficiently close
to each other.
Thus \[f\] is locally distance increasing on \[U_0\] with respect to the
metric \[\dist_U\].

Since there are no critical points in \[J\], it follows from the Sullivan
classification of Fatou components that every critical orbit must converge
to an attracting or parabolic cycle. Therefore the intersection \[\Pi=P\cap J\]
can be described as the set of all parabolic periodic points. If there are
no parabolic points, then
\[J\subset U_0\], and it follows that \[f\] is locally distance increasing
on \[J\] with respect to the metric \[\dist_U\]. However, if \[f\]
has parabolic points, then the Poincar\'e metric becomes
infinite at such points. We will modify this metric near these points
so as to obtain a better behaved metric.

It will be convenient to choose coordinates so
that the Julia set \[J\] is contained in the finite plane \[\C\].
The Poincar\'e metric on \[U=\widehat\C\ssm P\] has the form \[\rho(z)|dz|\]
throughout \[U\cap\C\], where \[\rho\] extends to a function
$$	\rho:\C~\to~(0\,,\,\infty] $$
which is continuous everywhere, but takes the value \[+\infty\] on the points
of \[P\cap\C\]. In particular, \[\rho(z)\] tends to infinity as \[z\] tends
to any point in the set \[\Pi=P\cap J\] of parabolic
points. First suppose that these parabolic points are all fixed
under \[f\], with multiplier \[+1\]. Let \[M>0\] be a constant
which is large enough so that
$$	M~>~\rho(z)/|f'(z)|\eqno (25) $$
for every \[z\]
in the finite set \[f^{-1}(\Pi)\ssm\Pi\]. (This
ratio \[\rho(z)/|f'(z)|\] is always finite at the points of \[f^{-1}(\Pi)
\ssm\Pi\], since such points are non-critical and
outside \[P\].) Let \[N_\epsilon (\Pi)\] be the open neighborhood
of Euclidean radius \[\epsilon\] about the parabolic set \[\Pi\]. We
define a new Riemannian metric \[\eta(z)|dz|\] on
the open set \[U'=N_\epsilon(\Pi)\cup
 (\C\ssm P)\] by setting
$$	\eta(z)\=\cases{M & for \[z\in N_\epsilon(\Pi)\],\cr
 \rho(z) & for \[z\not\in N_\epsilon(\Pi)\].\cr}\eqno (26) $$   
This function \[\eta\] has a jump discontinuity
on the circle of radius \[\epsilon\] about each parabolic point; but
this will not cause any difficulty.
As usual, we use this Riemannian metric to define a distance
\[\dist_\eta(x,y)\] between points of \[U'\] as the infimum of lengths
$$	{\rm length}_\eta(\Gamma) \=\int_\Gamma\;\eta(z)|dz|~,  $$
where \[\Gamma\] ranges over all piecewise smooth paths joining the two points
within \[U'\].

{\QP\it
If \[\epsilon\] is sufficiently small, then we will prove that \[f|_J\]
is locally distance increasing with respect to this metric \[\dist_\eta\].\ss}

\ni Since \[\eta(\Pi)=\infty\], we can certainly choose \[\epsilon\] small
enough so that \[\rho(z)>M\] throughout
the neighborhood \[N_{2\epsilon}(\Pi)\]. Evidently this will guarantee that any
Euclidean straight line segment lying within \[N_\epsilon(\Pi)\] is the
unique \[\eta$-geodesic of minimal \[\eta$-length joining its endpoints.
To study behavior of \[f\]
near a parabolic point \[z_0\], we proceed as follows. Assuming for
convenience that \[z_0=0\], we can write the local power series expansion as
$$	f(z)\=z\Big(1\,+\,az^m\,+\,{\rm(higher~terms)}\Big)~, $$
where \[m\ge 1\] and \[a\ne 0\]. Using the Leau-Fatou Flower Theorem, we
obtain the estimate
$$	{az^m\over |az^m|}~\to~ 1 $$
as \[z \to 0\] within the Julia set \[J\]. Hence
$$	f(z)\= z\Big(1+|az^m|+o(z^m)\Big)~,\qquad
	f'(z)\= 1+(m+1)|az^m|\,+\,o(z^m)  \eqno (27) $$
as \[z\to 0\] within \[J\]. It follows easily that \[f|_J\] is distance
increasing throughout some neighborhood of each parabolic point.
Furthermore, if \[\epsilon'\] is sufficiently small, it follows that
$$	|f(z)|~>~|z|\qquad{\rm and}\qquad
	|f'(z)|~>~1\qquad{\rm for}\qquad 0~<~|z|~<~\epsilon'~.\eqno (28) $$
To prove that \[f\] is locally distance increasing at points of \[J\ssm\Pi\]
we will prove the infinitesimal form of the required inequality. That is
we will prove that
$$	\eta(f(z))\,|f'(z)|~>~ \eta(z) 	\eqno(29)$$
throughout \[J\ssm\Pi\], and hence that
$$	\int_\Gamma \eta(f(z))\,|df(z)|\=\int_\Gamma \eta(f(z))\,|f'(z)dz|~>~
	\int_\Gamma \eta(z)\,|dz| $$
for any path \[\Gamma\] in some neighborhood of \[J\ssm\Pi\].
There are four cases to consider: If both \[z\]
and \[f(z)\] belong to \[N_\epsilon(\Pi)\], then (29) follows from (28).
If neither \[z\] nor
\[f(z)\] belongs to \[N_\epsilon(\Pi)\], then it follows from the corresponding
property of the Poincar\'e metric. If \[z\] belongs to \[N_\epsilon(\Pi)\]
but \[f(z)\] does not, then it follows since
$$	\eta(f(z))\,|f'(z)|\=\rho(f(z))\,|f'(z)|~>~\rho(z)~>~M\= \eta(z)~.$$
Finally, if \[z\in f^{-1}N_\epsilon(\Pi)\ssm N_\epsilon(\Pi)\], then we proceed
as follows. A compactness argument shows that \[f^{-1}N_\epsilon(\Pi)\]
shrinks down to \[f^{-1}(\Pi)\] as \[\epsilon\;{\scriptstyle\searrow}\; 0\].
Thus, given \[\epsilon'\] we can find \[\epsilon\] so
that every point \[z\] of \[f^{-1}N_\epsilon(\Pi)\] has Euclidean distance
at most \[\epsilon'\]
from some point \[\hat z\in f^{-1}(\Pi)\]. If \[\hat z\in\Pi\], then since
\[f(z)\in\N_\epsilon(\Pi)\] it
follows from (28) that \[z\in N_\epsilon(\Pi)\] also, contradicting the
hypothesis that \[z\not\in N_\epsilon(\Pi)\].
On the other hand, if \[\hat z
 \ne\Pi\] then it follows from (25) that 
\[M\,|f'(\hat z)|>\eta(\hat z)\]. If \[\epsilon\] is sufficiently small,
it evidently follows that \[M\,|f'(z)|>\eta(z)\] also.\ss

This proves B.3 in the special case where every parabolic point is a fixed
point
of multiplier \[+1\]. To handle the general case, choose a positive integer
\[k\] so that every parabolic point of \[f^{\circ k}\] is a fixed point
of multiplier \[+1\], and let \[\dist_\eta(x,y)\] be a metric with the
required property for \[f^{\circ k}\]. Then \[f\] itself will be locally
distance increasing on \[J\] with respect to the metric
$$	\dist'(x,y)\=\sum_{j=0}^{k-1}\dist_\eta\Big(f^{\circ j}(x)\,,\,
	f^{\circ j}(y)\Big)~. $$
This proves B.3, and completes the proof of B.1.\QED\ss

Next we will prove the following.
Let \[f\] be rational of degree \[n\ge 2\] with Julia set \[J\].

{\QP{\bf Lemma B.4.} \it
Suppose that there exists a closed disk \[\Delta^*\subset\widehat\C\] which:

{\rm (a)} contains all of the critical values of \[f\] in its interior,

{\rm (b)} satisfies \[f(\Delta^*)\subset\Delta^*\], and

{\rm (c)} eventually absorbs every orbit in the Fatou set.

\noindent Then \[f|_{J}\] is topologically conjugate to the one-sided
shift on \[n\] symbols.\ss}

{\bf Proof.} It follows from (b) that the interior of such a disk \[\Delta^*\]
is contained in the Fatou set. Therefore, by (a) and B.1, the Julia set \[J\]
is totally disconnected. As in
the proof of B.2, the complementary closed disk \[\Delta\] contains \[J\],
and \[f^{-1}(\Delta)\] splits as a disjoint union \[\Delta_1\cup\cdots\cup
\Delta_n\], where each \[\Delta_i\] maps homeomorphically onto \[\Delta\].
But now we have the additional information that
each \[\Delta_i\] is a subset of \[\Delta\]. It follows inductively that each
finite intersection
$$	\Delta_{i_0\,i_1\,\cdots\,i_k}\=	\Delta_{i_0}
  \cap f^{-1}(\Delta_{i_1})\cap\cdots\cap f^{-k}(\Delta_{i_k})$$
is a closed topological disk. In fact
$$	f:\Delta_{i_0\,i_1\,\ldots\,i_k}\iso
  \Delta_{i_1\,i_2\,\ldots\,i_k}\iso\cdots\iso \Delta_{i_{k-1}\,i_k}\iso
  \Delta_{i_k}\iso \Delta~, ~ $$
Now given an infinite sequence \[I=(i_0\,,\,i_1\,,
\,i_2\,,\,\ldots)\], it follows that the intersection \[J_I\] of the nested
sequence
$$	\Delta_{i_0}\,\supset\,\Delta_{i_0\,i_1}\,\supset
\,\Delta_{i_0\,i_1\,i_2}\,\supset\,\cdots $$
is compact, connected, and non-vacuous. This intersection is contained in
the Julia set, since by condition (c) every orbit outside of the Julia set
eventually leaves the disk \[\Delta\]. Hence by B.1 each \[J_I\] is a single
point.

Equivalently, to each point \[z\] of \[J\] we can assign a sequence \[I(z)=
(i_0\,,\,i_1\,,\, i_2
\,,\,\ldots)\] of integers between \[1\] and \[n\]
by the condition that \[z_j\in D_{i_j}\] where
$$      f\;:\;z=z_0\mapsto z_1\mapsto z_2\mapsto\cdots $$
is the orbit of \[z\]. Evidently \[J_I\] consists of the unique point which
has symbol sequence equal to
\[I\]. Since \[f\] maps \[J_I\] to \[J_{\sigma(I)}\], where
$$ \sigma(i_0\,,\,i_1\,,\,\ldots)\=(i_1\,,\, i_2\,,\,\ldots) $$
is the shift operator, this proves that the correspondence
\[z\mapsto I(z)\] is a topological conjugacy from \[f|_J\] onto the one-sided
shift on \[n\] symbols.\QED

We will only try to apply this lemma in the following very special case.

{\QP{\bf Theorem B.5.} \it If \[f\] has only two critical values, both
in the same Fatou component, then \[f|_J\] is topologically conjugate
to a one-sided shift.\ss}

Note: The hypothesis that \[f\] has only two critical values implies that
\[f\] has only two critical points (making use of the fact
that the twice punctured sphere has
free cyclic fundamental group). By way of contrast, a Chebyshev polynomial
of degree \[n\ge 3\] has only
three distinct critical values on the Riemann sphere but has \[n\]
distinct critical points.\ss

{\bf Proof of B.5.} Let \[U_0\] be a simply connected open set bounded by
a smooth
Jordan curve, which satisfies \[f(U_0)\subset U_0\], and which eventually
absorbs every orbit in the Fatou set. Thus in the hyperbolic case \[U_0\]
will be a neighborhood of the attracting point, while in the parabolic
case \[U_0\] will be a carefully chosen attracting petal. We
can assume that the boundary of \[U_0\] is disjoint from the two critical
orbits. Now define a sequence of connected open sets with smooth boundary
$$	U_0~\subset~U_1~\subset~U_2~\subset~\cdots $$
by setting \[U_{k+1}\] equal to the connected component of \[f^{-1}(U_k)\]
which contains \[U_k\]. If \[U_k\] is simply connected and contains
at most one critical value, then it is easy to show that \[U_{k+1}\] is also
simply connected. Thus we can continue this construction
until we obtain some \[U_m\] which
contains both critical values. The closure \[\overline U_m\] will then be
the required disk \[\Delta^*\], so that we can apply B.4.\QED
\bs\bs

\headline={%
    \ifnum\pageno > 1 {%
        \hss\tenrm\ifodd\pageno C. CROSS-RATIO FORMULAS\hss\folio
        \else\folio\hfill MAPS WITH TWO CRITICAL POINTS\hfill\fi}%
    \else \hss\vbox{}\hss
    \fi}

\def\(({{ \Big(\!\!\Big(}}
\def\)){{ \Big)\!\!\Big)}}
\def\vr{\vrule height .2in width 1pt}
\def\ab{\above 1pt}
\def\xx#1#2#3#4{{#1~\ab#3~}\!\vr\!{~#2\ab~#4}}

\cl{\bf Appendix C. Cross-Ratio Formulas.}\ms

In 1.6 we saw that the invariant \[X\] of 1.1 can be defined as a cross-ratio.
This appendix will give corresponding formulas for the invariants
\[Y_1\] and \[Y_2\]. I will use the nonstandard notation
$$\xx{p}{q}{r}{s}\={(p-q)(r-s)\over(p-r)(q-s)}~\in~\widehat\C $$ 
(product of row differences, divided by product of column differences).
This cross-ratio symbol can be characterized as follows. It
is well defined unless three of the four points \[{p~q\atop r~s}\] in
\[\widehat\C\] coincide, and it takes the value:
$$\cases{ 0 & if two entries in the same row are equal,\cr
	\infty & if two entries in the same column are equal,\cr
	1 & if two diagonally opposite entries are equal.} $$
Furthermore, if three of the four variables \[p,q,r,s\]
are fixed and distinct, then
it represents a M\"obius transformation of the remaining variable.

As an example, since a M\"obius transformation which fixes three points
is the identity, it follows that
$$	\xx{\,x\,}{\,0\,}{\infty}{1}\=x~.$$
Furthermore, it follows easily that
$$	\xx{\phi(p)}{\phi(q)}{\phi(r)}{\phi(s)}\=\xx{p}{q}{r}{s} $$
for any M\"obius transformation \[\phi\].
This symbol is independent of the order of rows or columns, that is
it satisfies the symmetry relations
$$ \xx{p}{q}{r}{s}\=\xx{r}{s}{p}{q}\=\xx{s}{r}{q}{p}\=\xx{q}{p}{s}{r}~.$$\ss

Now let \[c_1\] and \[c_2\] be the critical points of the bicritical map \[f\],
let \[v_i=f(c_i)\] be the critical values, and let \[p_i\in
f^{-1}(c_i)\] be any one of the \[n\] preimages of \[c_i\]. For example,
using the usual normal form $$f(z)\=(az^n+b)/(cz^n+d)~,$$
we can take
$$\matrix{	c_1\=\infty~,\quad& v_1\=a/c~,\quad& p_1\in\root n\of{-d/c}~\cr
		c_2\=0~,\quad& v_2\=b/d~,\quad& p_2\in\root n\of{-b/a}~.}$$
Using these notations, we can write 1.6 in the form
$$	\xx{c_1}{v_1}{c_2}{v_2}\={-v_2\over v_1-v_2}\={-b/d\over a/c-b/d}
\={-bc\over ad-bc}\=-X~. $$
Now consider the cross-ratio
$$  \xx{c_1}{v_2}{p_1}{c_2}\={p_1\over v_2}\= {\root n\of{-d/c}\over b/d}~. $$
This is well defined only up to multiplication by an \[n$-the root of unity.
However its \[n$-th power
$$\left(\xx{c_1}{v_2}{p_1}{c_2}\right)^n
\={-d/c\over b^n/d^n}\={-d^{n+1}\over b^nc} $$
is always well defined as an element of \[\widehat\C\]. If \[c\ne 0\]
then we can multiply numerator and denominator by \[c^{n-1}\], thus reducing
this expression to the form
$$	\left(\xx{c_1}{v_2}{p_1}{c_2}\right)^n\=
{-c^{n-1}d^{n+1}\over b^nc^n}\={-c^{n-1}d^{n+1}\over (ad-bc)^n}
\left({bc\over ad-bc}\right)^{-n}\={-Y_2~\over ~X^n}~.\eqno(30) $$
This proves that we can compute \[Y_2\] by a cross-ratio whenever \[X\ne 0\].
Still assuming that \[c\ne 0\], we can also consider the expression
$$ \left(\xx{c_1}{v_1}{p_1}{c_2}\right)^n
 \=(p_1/v_1)^n  \={-d/c\over a^n/c^n}
  \=	{c^{n-1}d^{n+1}\over a^nd^n}\={-Y_2\over(X+1)^n}~.\eqno (30') $$
Note that the left side of $(30')$ is indeterminate if and only if \[c=0\],
or equivalently if and only if
\[c_1=v_1=p_1\]. But whenever this happens, it follows that \[Y_2=0\].
Since \[X\] and \[X+1\] cannot both be zero, this shows that we can
compute \[Y_2\] in all cases. Interchanging the roles of the two critical
points, we find similar expressions for \[Y_1\].\bs\bs

\headline={%
    \ifnum\pageno > 1 {%
      \hss\tenrm\ifodd\pageno D. ENTIRE AND MEROMORPHIC MAPS\hss\folio
        \else\folio\hfill MAPS WITH TWO CRITICAL POINTS\hfill\fi}%
    \else \hss\vbox{}\hss
    \fi}


\cl{\bf Appendix D: Entire and Meromorphic Maps}\ms

Douady has pointed out to me that there is an analogous theory of
entire transcendental or meromorphic maps
which have only two singular values. 

{\QP{\bf Lemma D.1.} \it Suppose that \[f:U\to\widehat\C\] is
a holomorphic map with only two singular values, say \[v_1\] and \[v_2\],
where \[U\] is a connected open subset of \[\widehat\C\]. 
If \[f\] has infinite degree, then \[U\] must be
the complement of a single point \[s\in\widehat\C\], and \[f\]
must map \[\widehat\C\ssm\{s\}\] onto  \[\widehat\C\ssm\{v_1\,,\,v_2\}\] by a
free cyclic covering map.\ss}

Here, by a {\bit singular value\/} \[v\in\widehat\C\] we mean either a critical
value, or an asymptotic value for the map \[f\]. (By definition,
\[z\] is an  {\bit asymptotic value\/} if it can be described as the limit
\[z=\lim_{t\to 1} f(p(t))\] where \[p:[0,1)\to U\] is a path
which eventually leaves any compact subset of \[U\].)
The Lemma asserts that the \[f\] in question
cannot have any critical values, so we could equally well describe it
as a map with two asymptotic values.\ss

For a quite different characterization of this family of maps see [DK2].\ss

{\bf Proof of D.1.} It is not difficult to see
that \[f\] must carry the complement \[U\ssm f^{-1}
\{v_1\,,\,v_2\}\] onto \[\widehat\C\ssm\{v_1\,,\,v_2\}\] by a covering
map. Since the fundamental group of the image space is cyclic, this can
only be a cyclic covering. If it has infinite degree, then it must be a
universal covering. Hence \[U\ssm f^{-1}\{v_1\,,\,v_2\}\] must be conformally
equivalent to the universal covering of \[\widehat\C\ssm\{v_1,v_2\}\], or in
other words to \[\C\] itself. Since \[U\] cannot be the entire Riemann sphere,
this proves that \[U=\widehat\C\ssm\{s\}\] for some unique point \[s\],
and that
$$	f:U~\to~\widehat\C\ssm\{v_1\,,\,v_2\} $$
is a free cyclic covering map. Note that \[f\] has an
essential singularity at \[s\].\QED\ss

If we put the two singular values at zero and infinity, then one example
of a universal covering of \[\widehat\C\ssm\{v_1,v_2\}=\C\ssm\{0\}\] is
given by the exponential map
$$	\exp:\C~\to~\C\ssm\{0\}~.$$
Hence \[f\] can be described as the composition \[f=\exp\circ\mu\], where
\[\mu\] is some conformal isomorphism from \[\widehat\C\ssm\{s\}\] onto
\[\C\]. Evidently we can express any such \[\mu\] as a M\"obius transformation
$$	\mu(z)\={az+b\over cz+d}\qquad{\rm with}\qquad ad-bc\ne 0~, $$
so that
$$ f(z)\=\exp\circ\mu(z)\=\exp\left({az+b\over cz+d}\right)~, \eqno (31) $$
with essential singularity at \[s=-d/c\in\widehat\C\].

In practice, it will be convenient to conjugate \[f=\exp\circ\mu\] by \[\mu\],
so as to obtain a normal form
$$	g\=\mu\circ f\circ\mu^{-1}\=\mu\circ\exp~.$$
Thus 
$$	g(w)\={ae^w+b\over ce^w+d}~, \eqno (32) $$
or equivalently \[
 g(2w)\={ae^w+be^{-w}\over ce^w+de^{-w}}\]. The essential
singularity of \[g\] is at infinity, and its singular values  are
\[v_1=a/c\] and \[v_2=b/d\]. Note the identity
$$	g(w+2\pi i)\=g(w)~. $$
In fact we can describe \[g\] as a composition
$$	\C~\buildrel{\rm projection}\over\longrightarrow~\C/2\pi i\Z~
	\buildrel\cong\over\longrightarrow~\widehat\C\ssm\{a/c\,,\,b/d\}~.$$
Since this description is rather rigid, we are free to change the coordinate
\[w\] only by a translation, or a sign change followed by translation.
If we conjugate \[g\] by a translation, then the resulting map \[w\mapsto
g(w+t)-t\] has the form (32) with matrix of coefficients
$$	\left[\matrix{ e^t(a-tc) & b-td\cr	e^t c & d}\right]\qquad({\rm
or~ equivalently} \qquad\left[\matrix{ e^{t/2}(a-tc) & e^{-t/2}(b-td)\cr
 e^{t/2} c & e^{-t/2}d}\right] $$
if we prefer to work only with unimodular matrices).
Thus the additive group of translations acts on the manifold \[{\bf
 PGL}(2,\C)\] by a correspondence which we can write as
$$	\Phi_t:\left[\matrix{a&b\cr c&d}\right]~\mapsto
\left[\matrix{1&-t\cr 0&1}\right]\left[\matrix{a&b\cr c&d}\right]
\left[\matrix{e^{t/2}&0\cr 0&e^{-t/2}}\right]~. \eqno (33)$$
The difference
$$	v_1-v_2\={a\over c}-{b\over d}\={ad-bc\over cd} $$
remains invariant under this transformation, since \[v_1\] and \[v_2\]
are both translated by the same constant \[-t\]. Let us define
\[X\] to be the reciprocal
$$	X\={1\over v_1-v_2}\={cd\over ad-bc}~\in~\C~. \eqno (34)$$
This number \[X\] is a conformal conjugacy invariant. For if we replace
\[g(w)\] by \[-g(-w)\] then the matrix \[\left[\matrix{a & b\cr c & d}\right]\]
of coefficients will be replaced by
$$\left[\matrix{-1&0\cr 0&1}\right]\left[\matrix{a&b\cr c&d}\right]
 \left[\matrix{0&1\cr 1&0}\right]\=\left[\matrix{-b & -a\cr ~d & ~c}\right]~,
	\eqno (35)$$
so that \[(v_1\,,\,v_2)\] is replaced by \[(-v_2\,,\,-v_1)\], and again
the expression (34) remains invariant. Since
any cooordinate change which carries (32) to another expression of the same
form is a composition of translations and sign changes, it follows
that \[X\] is indeed a conformal conjugacy invariant.\ss

We will prove several statements which are reminiscent of the descriptions
of moduli space
in \S1 and \S2. For each fixed \[X\], we need one further complex number
to give a complete conjugacy class invariant. Furthermore, given \[X(f)\] and
given \[\lambda\ne 0\], there is a unique conjugacy class of maps \[(f)\]
which have a fixed point of multiplier \[\lambda\]. However, unlike the finite
degree situation of \S1, it is necessary to treat the cases \[X=0\] and
\[X\ne 0\] separately in order to obtain reasonable moduli spaces.\ss

{\bf Entire Transcendental Case.} If \[X=0\], or equivalently if
the essential singularity occurs at one of the two singular values,
then putting \[s=v_1\] at infinity and putting \[v_2\] at zero, the normal
form (32) reduces to
$$	g(w)\=k e^w\qquad{\rm for~some~constant}
 \qquad k\in\C\ssm\{0\}~. $$
Here \[k=g'(0)\] is a conjugacy class invariant. 
Thus we obtain the {\bit exponential family\/}, which has
been studied by several authors. (See for example [D].)
We could also ue the normal form
$$ 	f(z)=g(k z)/k\=e^{k z}~, $$
again with \[s=\infty\] and zero as singular values, and again with invariant
\[k=f'(0)\]. If we specify that \[g\] has a fixed point \[w_0=g(w_0)\]
of multiplier \[g'(w_0)= \lambda\], then since
\[g'(w_0)=g(w_0)=w_0\] it follows that \[w_0=\lambda\], and
hence that the invariant \[k=w_0 e^{-w_0}\] is uniquely
determined by \[\lambda\].\ss

{\bf Meromorphic Case.} If  \[X\ne 0\], or equivalently if
the three distinguished points \[s\,,\,v_1\,,\,v_2\] are pairwise distinct,
then \[g:\C\to\widehat\C\ssm\{v_1,v_2\}\] is a meromorphic map. (For
a survey of meromorphic dynamics, see [Be].)
In analogy with \S1 we will construct a conjugacy invariant \[Y=Y_1+Y_2\]. Let
$$	Y_1\={c\over d}\exp{a\over c}~,\qquad Y_2
 \={d\over c}\exp{-b\over d}~, \eqno(36) $$
with product
$$	Y_1Y_2\=\exp\left({a\over c}-{b\over d}\right)\=\exp\left({ad-bc
	\over cd}\right)\=e^{1/X}~. \eqno (37) $$
Then it is not difficult to check that both \[Y_1\] and \[Y_2\] are invariant
under the action (33), and that \[Y_1\] and \[Y_2\] are interchanged by
the transformation (35). Hence the sum \[Y=Y_1+Y_2\] is a conjugacy class
invariant. The two quantities \[Y_1\] and \[Y_2\] individually will
be described as {\bit half-invariants\/}, since only the unordered pair
\[\{Y_1\,,\,Y_2\}\] is actually a conjugacy invariant.

Given \[X\ne 0\] and \[Y\], the conjugacy class is
uniquely determined. In fact we can compute the
unordered pair \[\{Y_1\,,\,Y_2\}\] by solving a quadratic equation. Given
\[Y_1\], we can assume, by using the action (33), that \[c=d\] so that
\[Y_1=e^{a/c}\]. Thus we can compute the ratio \[a/c\], up to a summand
in \[2\pi i\Z\], and then compute the ratio \[b/d\]
from the identity \[a/c-b/d=1/X\]. Thus the function \[g(w)=(ae^w+b)/(ce^w+d)\]
is uniquely determined, up to a summand in \[2\pi i \Z\] which can be
removed by a corresponding translation of the coordinate \[w\]. This
proves the following.

{\QP{\bf Lemma D.2.} \it The moduli space for meromorphic maps with two
singular values is biholomorphic to \[(\C\ssm\{0\})\times\C\], with coordinates
\[X\] and \[Y\].\ms}

{\bf Cross-Ratios.} Just as in Appendix C, we can express the \[Y_i\]
as cross-ratios. For example, using the normal form (31), one easily checks
that:
$$	\xx{v_1}{s}{f(v_1)}{v_2}\= -Y_1~.$$
However, it does not seem possible to compute the invariant \[X\] as a
cross-ratio.\ss

{\bf Bad Topology.} It might seem natural to combine the moduli space
for the exponential family, together with the moduli space
for our meromorphic family, into a larger moduli space, isomorphic to a union
\[(\C\ssm\{0\})\cup(\C\ssm\{0\})\times\C)\]  . However this does not
yield any useful result since the natural topology on the union is not
Hausdorff. To see this, consider the two maps \[z\mapsto e^z\] and
\[z\mapsto 2e^z\] in the exponential family.
Although these are not conjugate, I will show that given any neighborhood
\[N_1\] of the first and any neighborhood \[N_2\] of the second
within the space
of maps with two singular values, some map in \[N_1\] is conjugate to a
map in \[N_2\]. Consider the analytic function \[\xi\mapsto
2\xi^2e^{1/\xi}\], which has an essential singularity at the origin.
Since this function omits the values \[0\] and \[\infty\], it follows from
Picard's Theorem that we can choose a sequence of values of \[\xi\] tending
to zero which satisfy the equation
\[2\xi^2e^{1/\xi}=1\]. For each such \[\xi\], consider the
two meromorphic maps
$$	{e^w\over \xi e^w+1}\qquad{\rm and}\qquad {2e^w\over 2\xi e^w+1}~.$$
A brief computation shows that both have invariants \[X=\xi\] and \[Y=\xi
e^{1/\xi}+1/\xi\]. Hence the two are conjugate, although
these sequences tend to non-conjugate limits as \[\xi\to 0\].\QED\ss

{\bf Comparison with Bicritical Maps.} Our invariants for meromorphic maps
look rather different from the invariants of \S1, but in fact there is
a definite relationship. As in [DGH], we can approximate the
exponential map by the unicritical polynomials \[E_n(w)=(1+w/n)^n\].
Hence we can approximate the meromorphic function \[g\] of (32) by
the bicritical maps
$$	g_n(w)\={aE_n(w)+b\over cE_n(w)+d} $$
throughout any compact subset of \[\C\]. A
straightforward computation then shows that
$$ X(g)\=\lim_{n\to\infty} X(g_n)/n\quad{\rm and~that}\qquad
	Y_i(g)\=\lim_{n_\to\infty} Y_i(g_n)/X(g_n)^n~.$$\ss

{\bf The Symmetry Locus.}
There is a non-trivial M\"obius automorphism which commutes
with \[g\] if and only if \[Y_1(g)=Y_2(g)\]. In view of (37), this means that
$$	Y_1\=Y_2\=(Y/2)\=\pm e^{1/2X}~. $$
If we choose the plus sign, then the most general example is conjugate to
$$	g(w)\=k{e^w-1\over e^w+1}\=k\,\tanh(w/2)~, $$
with a fixed point of multiplier \[k/2\] at the center of symmetry. Here
\[Y_1=Y_2=e^k\] and \[X=1/2k\]. (For a discussion of the
tangent family, see for example [DK1].) If we choose the minus sign,
then we can take
$$	g(w)\=k{e^w+1\over e^w-1}\=k\,\coth(w/2)~. $$
This is just the image of \[k\,\tanh(w/2)\] under the involution
\[g\mapsto J_g\circ g\] of 1.5.\ss

{\bf Fixed Points and the Curves \[\per_1(\lambda)\].} Suppose that \[g\]
has a fixed point \[w_0=g(w_0)\] with multiplier \[\lambda=g'(w_0)\].
We will prove that
the conjugacy class \[(g)\] is uniquely determined by \[X(g)\] together
with \[\lambda\]. We again use the normal form (32), but now we translate
coordinates so that the fixed point is at the origin. Then the fixed
point equation \[0=g(0)\] reduces
to the equation \[a+b=0\]. If \[ad-bc=1\], we have \[X=cd\],
and we see easily that
\[g'(0)=1/(c+d)^2\] so that \[c+d=\pm 1/\sqrt\lambda\] is uniquely determined
up to sign. Hence we can solve for the unordered pair \[\{c,d\}\] up to
sign. We can then solve for \[a\] and \[b=-a\] since \[ad-bc=a(c+d)=1\].
Thus the conjugacy class is uniquely determined.
Multiplying all coefficients by a common constant, this solution can be
written as
$$	\left[\matrix{a&b\cr c&d}\right]\=\left[\matrix{2\lambda&-2\lambda\cr
1+r&1-r}\right] $$
with determinant \[ad-bc=4\lambda\], where \[r^2 = 1 - 4\lambda X\ne 1\].
From this, one can easily write down a precise but somewhat
complicated formula for the invariant \[Y\] as a function
of \[\lambda\] and \[X\]. This function is holomorphic, since it can be
expressed as the sum of a power series, convergent for \[|r|<1\],
in which only even powers of \[r\] appear.\ss

{\bf Real forms.} A meromorphic map with two singular values
commutes with some antiholomorphic involution \[\alpha\]
if and only if the invariants \[X\] and \[Y\] are both real. (Compare \S5.)
If \[Y^2>
4\,e^{1/X}\], then the  half-invariants \[Y_1\,,\,Y_2\] are also real, and
the singular values are fixed by \[\alpha\]. In this case, using the
normal form (32) with real coefficients, the map
\[g\] carries \[\R\] diffeomorphically onto an open interval in
\[\R\cup\{\infty\}\], bounded by the two singular values. On the
other hand, if \[Y^2< 4\,e^{1/X}\], then the half-invariants
\[Y_1\,,\,Y_2\] are complex conjugate, and the
singular values are interchanged by \[\alpha\].
In this case we can use a normal form
$$	z~\mapsto~{a\,\tan(z)+b\over c\,\tan(z)+d} $$
with real coefficients. These maps carry \[\R\] onto \[\R\cup\infty\]
by a composition
$$	\R~\buildrel{\rm projection}\over\longrightarrow~\R/\pi\Z~
	\buildrel\cong\over\longrightarrow~\R\cup\{\infty\}~, $$
with degree \[\pm\infty\]. On the symmetry locus, with
\[Y^2=4\,e^{1/X}\], there are two possible choices for \[\alpha\]
and we can use either normal form.
\bs

\headline={%
    \ifnum\pageno > 1 {%
        \hss\tenrm\ifodd\pageno REFERENCES\hss\folio
        \else\folio\hfill MAPS WITH TWO CRITICAL POINTS\hfill\fi}%
    \else \hss\vbox{}\hss
    \fi}

\cl{\bf References.}\ss

\ref [A] L. Ahlfors, ``Lectures on Quasiconformal Mappings'',
Van Nostrand 19<66, Wadsworth 1987.

\ref [AH] P. Alexandroff and H. Hopf, ``Topologie'', Berlin 1935.

\ref [BB] R. Bam\'on and J. Bobenrieth, The rational maps \[z\mapsto
1/\omega z^d\] have no Herman rings, Proceedings A.M.S. ??

\ref [Be] W. Bergweiler, Iteration of meromorphic functions.
 Bull. Amer. Math. Soc. (N.S.) 29 (1993), 151--188.

\ref [BH] B. Branner and J. H. Hubbard,  The iteration of cubic polynomials,
Part II: patterns and parapatterns, Acta Math. {\bf 169} (1992) 229-325.

\ref [Bo] T. Bousch, Sur quelques probl\`emes de la dynamique holomorphique,
th\`ese, Univ. Paris-Sud (Orsay) 1992.

\ref [CG] L. Carleson and T. Gamelin, ``Complex Dynamics'', Springer 1993.

\ref [D] R. Devaney, \[e^z\]: dynamics and bifurcations. Internat. J. Bifur.
 Chaos Appl. Sci. Engrg. 1 (1991), 287--308.

\ref [DGH] R. Devaney, L. Goldberg and J. Hubbard, A dynamical approximation
to the exponential map by polynomials, Preprint, circa 1986.

\ref [DH1] A. Douady and J. H. Hubbard, \'Etude dynamique des polyn\^omes
complexes II, Publ. Math. Orsay (1984-5).

\ref [DH2] A. Douady and J. H. Hubbard, {\it On the dynamics of
polynomial-like
mappings,\/} Ann. Sci. Ec. Norm. Sup. (Paris) {\bf 18} (1985) 287-343.

\ref [DK1] R. Devaney and L. Keen; Dynamics of tangent. Dynamical
systems (College Park, MD, 1986--87), 105--111, Lecture Notes in Math., 1342,
Springer, Berlin-New York, 1988.

\ref [DK2] R. Devaney and L. Keen, Dynamics of maps with
 constant Schwarzian derivative. Complex analysis,
Joensuu 1987, 93--100, Lecture Notes in Math., 1351, Springer,
Berlin-New York, 1988. (See also:
Dynamics of meromorphic maps: maps with
polynomial Schwarzian derivative. Ann. Sci. \'Ecole Norm. Sup. (4) 22 (1989),
55--79.)

\ref [GK] L. Goldberg and L. Keen, The mapping class group of a generic
quadratic rational map and automorphisms of the 2-shift, Invent. Math.
{\bf 101} (1990) 335-372.

\ref [KK] L. Keen \& J. Kotus, Dynamics of the Family $\lambda\tan z$,
Stony Brook IMS preprint 1995\#8.

\ref [L] M. Lyubich, Some typical properties of the dynamics of rational maps,
Uspekhi Mat. Nauk {\bf 38} (1983) 197-198 (Russ. Math. Surv. {\bf 38} (1983)
154-155). See also: An analysis of the stability of the dynamics of rational
functions, Funk. Anal. i Pril. {\bf 42} (1984) 72-91 (Selecta Math. Sov.
{\bf 9} (1990) 69-90).

\ref [LS1] E. Lau amd D. Schleicher, Internal Addresses in the Mandelbrot
Set and Irreducibility of Polynomials, Stony Brook IMS Preprint 1994/19.

\ref [LS2] E. Lau amd D. Schleicher, Symmetries of fractals revisited,
Math. Intelligencer {\bf 18} (1996).

\ref [M1] J. Milnor, Dynamics in One Complex Variable: Introductory Lectures,
Stony Brook IMS Preprint 1990/5.

\ref [M2] ------, Geometry and dynamics of quadratic rational maps,
Experimental Math. {\bf 2} (1993) 37-83.

\ref [M3] ------, Periodic orbits, external rays and the Mandelbrot set:
an expository account, Stony Brook IMS Preprint 1999/3; to appear
Ast\'erisque.

\ref [Ma], P. Makienko, Totally disconnected Julia sets, MSRI preprint 042-95.


\ref [MSS] R. Ma\~n\'e, P. Sad and D. Sullivan,
On the dynamics of rational maps,
Ann. Sci. \'Ec. Norm. Sup. Paris {\bf 16} (1983) 193-217.

\ref [NS] S. Nakane and D. Schleicher, On multicorns and unicorns: the
dynamics of antiholomorphic polynomials, in preparation.

\ref [P]\  F. Przytycki, Iterations of rational functions: which hyperbolic
components contain polynomials? Fund. Math. 149 (1996), no. 2, 95-118.

\ref [R1] M. Rees, Positive measure sets of ergodic rational maps, Ann. Sci.
\'Ecole Norm. Sup. (4) {\bf 19} (1986) 383-407.

\ref [R2] ------, Components of degree two hyperbolic rational maps, Invent.
Math. {\bf 100} (1990) 357-382.

\ref [R3] ------, Views of parameter space: topographer and resident, in
preparation.

\ref [Sh] M. Shishikura, On the quasiconformal surgery of rational functions,
Ann. Sci.  \'Ec. Norm. Sup. Paris {\bf 20} (1987) 1-29.

\ref [Si] J. Silverman, The space of rational maps on \[{\bf P}^1\], preprint,
Brown University 1996.


\ref [Ste] N. Steinmetz, ``Rational Iteration: Complex Analytic Dynamical
Systems'', de Gruy\-ter 1993.

\ref [Sti] J. Stimson, Degree Two Rational Maps with a Periodic
  Critical Point, thesis, Liverpool 1993.

\ref [Su] D. Sullivan, Quasiconformal homeomorphisms and dynamics I, solution
of the Fatou-Julia problem on wandering domains, Ann. Math. {\bf 122} (1985)
401-418.

\ref [T] Tan Lei, Matings of quadratic polynomials, Ergodic
Th.~\& Dy.~Sy.~{\bf 12} (1992) 589-620.

\ref [TY] ------ and Yongcheng Yin, Local connectivity of the Julia set for
geometrically finite rational maps. Sci. China Ser. A {\bf 39} (1996) 39--47.
(See also Tan Lei's habilitation; Lyon 1997.)

\hskip 3.2in John Milnor

\hskip 3.2in Institute for Mathematical Sciences

\hskip 3.2in State University of New York

\hskip 3.2in Stony Brook, New York 11794-3660\ss

\hskip 3.2in Email: jack@math.sunysb.edu

\hskip 3.2in http://www.math.sunysb.edu/$\sim$jack

\bye